\documentclass[a4paper, 
10pt
]{article}

\usepackage[a4paper,left=3cm,right=3cm,top=2.5cm,bottom=2.5cm]{geometry}

\usepackage{subfig}
\usepackage{pgfplots}
\usepackage{pgfplotstable}
\usepackage{mathtools}
\usepackage{multicol}
\usepackage{comment}
\usepackage{booktabs}
\pgfplotsset{compat=1.5}
\usepackage{amssymb}
\usepackage{url}
\usepackage{bm}

\usepackage{pdfsync}
\usepackage{float}
\usepackage{tabularx}
\usepackage{enumerate}
\usepackage{array}
\usepackage{xspace}
\usepackage{tikz}
\usepackage{tikzsymbols}
\usetikzlibrary{calc,trees,positioning,arrows,chains,shapes.geometric,%
    decorations.pathreplacing,decorations.pathmorphing,shapes,%
    matrix,shapes.symbols, decorations.markings, patterns,fit}

\usepackage{authblk}
\usepackage{siunitx}

\usepackage{hyperref}
\usepackage{cleveref}

\usepackage{color, colortbl}

\newcommand{\RA}[1]{{\color{black}#1}}
\newcommand{\RB}[1]{{\color{black}#1}}

\usepackage{tikz}
\usetikzlibrary{calc,trees,positioning,arrows,chains,shapes.geometric,%
    decorations.pathreplacing,decorations.pathmorphing,shapes,%
    matrix,shapes.symbols}

\tikzset{
  >=stealth',
  punktchain/.style={
    rectangle, 
    rounded corners, 
    draw=black, very thick,
    text width=30em, 
    minimum height=3em, 
    text centered, 
    on chain},
  line/.style={draw, thick, <-},
  element/.style={
    tape,
    top color=white,
    bottom color=blue!50!black!60!,
    minimum width=8em,
    draw=blue!40!black!90, very thick,
    text width=10em, 
    minimum height=3.5em, 
    text centered, 
    on chain},
  every join/.style={->, thick,shorten >=1pt},
  decoration={brace},
  tuborg/.style={decorate},
  tubnode/.style={midway, right=2pt},
}

\DeclareRobustCommand\sampleline[1]{%
  \tikz\draw[#1] (0,0) (0,\the\dimexpr\fontdimen22\textfont2\relax)
  -- (1.5em,\the\dimexpr\fontdimen22\textfont2\relax);%
}

\begin{document}

\title{On the comparison of LES data-driven reduced order approaches for hydroacoustic analysis}

\author[1]{Mahmoud~Gadalla\footnote{mahmoud.gadalla@sissa.it}}
\author[2]{Marta~Cianferra\footnote{marta.cianferra@dia.units.it}}
\author[1]{Marco~Tezzele\footnote{marco.tezzele@sissa.it}}
\author[1]{Giovanni~Stabile\footnote{giovanni.stabile@sissa.it}}
\author[1]{Andrea~Mola\footnote{andrea.mola@sissa.it}}
\author[1]{Gianluigi~Rozza\footnote{gianluigi.rozza@sissa.it}}

\affil[1]{Mathematics Area, mathLab, SISSA, via Bonomea 265, I-34136
  Trieste, Italy}
\affil[2]{Industrial and Environmental Fluid Dynamic Research Group,
  University of Trieste, Trieste, Italy} 

\maketitle

\begin{abstract}
In this work, Dynamic Mode Decomposition (DMD) and Proper Orthogonal Decomposition (POD)
methodologies are applied to hydroacoustic dataset computed using Large Eddy Simulation (LES)
coupled with Ffowcs Williams and Hawkings (FWH) analogy. First, a low-dimensional
description of the flow fields is presented with modal decomposition analysis.
Sensitivity towards the DMD and POD bases truncation rank is discussed,
and extensive dataset is provided to demonstrate the ability of both algorithms
to reconstruct the flow fields with all the spatial and temporal frequencies necessary
to support accurate noise evaluation.
Results show that while DMD is capable to capture finer coherent
structures in the wake region for the same amount of employed modes, reconstructed flow fields using
POD exhibit smaller magnitudes of global spatiotemporal errors compared with DMD counterparts.
Second, a separate set of DMD and POD modes generated using half the snapshots is employed into two data-driven
reduced models respectively, based on DMD mid cast and POD with Interpolation (PODI).
In that regard, results confirm that the predictive character of both reduced approaches on
the flow fields is sufficiently accurate, with a relative superiority of PODI results over
DMD ones. This infers that, discrepancies induced due to interpolation errors in PODI
is relatively low compared with errors induced by integration and linear regression
operations in DMD, for the present setup. Finally, a post processing analysis on the
evaluation of FWH acoustic signals utilizing reduced fluid dynamic fields as input
demonstrates that both DMD and PODI data-driven reduced models are efficient and sufficiently
accurate in predicting acoustic noises.
\end{abstract}

\tableofcontents

\section{Introduction}
\label{sec:intro}
In several engineering fields, there has been recently
a growing need to include fluid dynamic performance evaluation criteria
associated with acoustic emissions. This is, for example, the case in
aircraft engine design, car manufacturing, and ship design
optimization. A particular motivation behind the present research
article is to investigate and propose a methodology that can be successively
used for the noise level prediction of naval propellers since the early design process
~\cite{cianferra2018hydrodynamic, Ianniello2014}.

The need for a reduction of the acoustic emissions through ship design optimization usually involves virtual prototyping and parametric high
fidelity simulations. Thanks to the increase of the available computational
resources, a deeper insight towards the complex physics associated
with hydroacoustic phenomena has become nowadays
affordable with unprecedented spatial and temporal scales (see, for
example, wall-resolving Large Eddy
Simulation (LES)~\cite{PosaBalaras,Kumar2017}).
However, the enormous data sizes resulting from such simulations pose several
challenges on the input/output operations, post-processing, or
the long-term data storage.
On the other hand, hybrid techniques such as, among others,
Detached-Eddy Simulation or wall-layer model LES (WLES),
have allowed obtaining eddy resolving field data with a reasonable use of computational resources.
These techniques are, however, still expensive for an early stage design process,
when a number of different geometric configurations has to be rapidly analysed to restrict the
range of variation of the principal design parameters.
Therefore, seeking a suitable data compression
strategy that allows extracting the most relevant and revealing
information in a reduced order manner, hence providing quick access as
well as efficient data storage, becomes a crucial asset.

Through multidisciplinary scale, efforts have been
made to realize optimal shape design for underwater noise
sources, including ship
hulls~\cite{Campana2006,Tezzele2018Dimension,tezzele2018model} and
propellers~\cite{mola2019marine,valdenazzi2019marine}, using efficient
geometrical parameterization
techniques~\cite{pygem,gadalla19bladex}. On the acoustic side, the
development of new generation noise prediction tools was considered a
major focus in this work. Particularly, based on the Ffowcs Williams
and Hawkings (FWH) analogy~\cite{fwh1969}, several improvements have
been developed. For example, Cianferra et al.~\cite{CianferraIA19}
compared several implementations of the non-linear quadrupole term, highlighting
its significant contribution to the overall hydrodynamic noise emissions in wide
range of frequencies. In a companion paper~\cite{CianferraAI18}, they
showed the effect of shape deformation on the radiated noise for
elementary geometries. On a more engineering level, the generated hydrodynamic noise from a benchmark marine propeller was evaluated in open sea
conditions~\cite{cianferra2018hydrodynamic,Cianferra2019} using FWH coupled with LES.


The aforementioned hydroacoustic models are typically described by a
system of non-linear Partial Differential Equations (PDEs), the
resolution of which results in the fluid dynamic fields necessary to
reproduce the noise source for an acoustic predictions. In
fact, a reliable reconstruction of the noise source is crucially
dependent on the flow structures and the resolved spatial and temporal
scales of the fluid dynamic fields (for a discussion,
see~\cite{Brentner2003}).

Besides LES works~\cite{Seror2001,Balaras2015,Nitzkorski2014, Kumar2017,
CianferraAI18, Cianferra2019, CianferraIA19}, several fluid dynamic models
have been employed in marine hydroacoustics. To name a few,
we mention the potential flow theory~\cite{Kerwin1986}, boundary
element method~\cite{SEOL2002}, Reynolds Averaged Navier-Stokes
equations
(RANS)~\cite{Gloerfelt2005, Ianniello2014, Baek2015}, DES~\cite{DiMascio2014}, and Direct Numerical
Simulation (DNS)~\cite{FosasdePando2014,Seo2008,SANDBERG2008}. Being
regarded as an optimal advance between RANS and DNS, recent literature
has reported LES to be the most suitable model which reproduces the
noise source with high level of realism~\cite{Nitzkorski2014,Balaras2015,Ianniello2013}.
The resolution of the described system of PDEs using
standard discretization methods (finite elements, finite volumes, finite
differences), which we will refer hereafter as the Full Order Model (FOM),
allows for high fidelity acoustic evaluation.

Although hybrid fluid dynamic techniques constitute a good compromise
between accuracy and computational cost for engineering purposes, they
are still expensive in the case of parametric analysis and shape optimization. To overcome this
issue, the development of a Reduced Order Model
(ROM)~\cite{salmoiraghi2016advances,rozza2018advances,hesthaven2016certified,morhandbook2019}
--- which alleviates both
the computational complexity and data capacity --- becomes
essential.
\RB{Moreover, several ROM developments have been realized to account for various effects including subgrid scales~\cite{HiStaMoRo2019}, heat transfer~\cite{GeoStaRoBlu2019} and stabilization of the Galerkin projection~\cite{StaBaZuRo2019}.}

One of the necessary assumptions to construct an efficient ROM
is that, the solution manifold of the underlying problem lies in a low
dimensional space and, therefore, can be expressed in terms of a linear
combination of a few number of global basis functions (reduced basis
functions). Among various techniques to generate such reduced
basis set, the Proper Orthogonal Decomposition (POD) and the
Dynamic Mode Decomposition (DMD) have been widely exploited due to their versatile
properties~\cite{stegeman2015proper}. POD provides a set of orthogonal and optimal basis
functions~\cite{Kunisch2008}, whereas DMD computes a set of
modes with an intrinsic temporal behavior, hence it is particularly
suited for time advancing problems~\cite{schmid2010dynamic}. One 
\RA{of the main goals of this work is to understand whether these
two established techniques are able to generate global basis functions
which include all the wide range of frequencies associated with
hydroacoustic phenomena, and therefore allow for accurate noise
predictions. In addition, the interest is in assessing the
efficiency of the modal decompositions algorithms, evaluating
the amount of modal shapes required to accurately reproduce
all high frequencies relevant to acoustic analysis. This investigation
is of course a first, fundamental step towards developing efficient
reduced order models for hydroacoustic applications. As will be thoroughly
discussed, in the present work we only considered data-driven ROMs,
but the POD and DMD efficiency assessment reported provides valuable
information also for the development of projection based ROMs.}

In literature, POD has been widely used for the past few decades to
identify the coherent structures of turbulent flows (see, among
others,~\cite{Sirovich1987,Berkooz1993}), and has been applied
towards various flow
conditions~\cite{MEYER2007,Rowley2004}. Correspondingly, DMD has been
also exploited~\cite{schmid2010dynamic,Seena2011}. An intuitive
question may arise considering the comparative performance. In that
regard, several works have carried out both DMD and POD on various
flow configurations. For instance, Liu et
al.~\cite{zhang2014identification} conducted comprehensive analysis,
concluding that DMD has the ability to clearly separate the flow
coherent structure in both spatial and spectral senses, whereas POD
was contaminated by other uncorrelated structures. Consistent to the
previous, in~\cite{Bistrian2015} it was noted that DMD is useful when the main
interest is to capture the dominant frequency of the phenomenon, while
the optimality of the POD modes prevails for coherent structure
identification that are energetically ranked. In high-speed train DES,
Muld et al.~\cite{muld2012flow} examined the convergence and reported
that the most dominant DMD mode requires a longer sample time to
converge when compared to the POD counterpart.

As discussed earlier, parametric studies such as shape optimization
\RA{can result in extremely high computational costs. Possible ways
to circumvent such issue could be degrading the high fidelity model,
or restricting the design parameter space sampling. Since the only analysed
parameter in the unsteady fluid dynamic problem is the
time, in this work we decided to rely on data-driven ROMs. A significant
source of error in data-driven ROMs is in fact associated
with the interpolation of modal coefficients based on the problem
parameters. Thus, such methods are an extremely valid alternative
in presence of a one dimensional parameter space, in which
interpolation errors are rather modest. In particular}, POD with
Interpolation (noted as PODI here-after) can be an adequate solution in
such scenarios~\cite{Ly2001}. The basic idea is to exploit POD on
selected ensemble of high fidelity solutions in the design space to
identify the set of optimal basis functions and associated projection
coefficients representing the solution dynamics. Such finite set of
scalar coefficients are then utilized to train a response surface that
allows predictions at parameter values that are not in the original
high-fidelity ensemble. It was demonstrated that PODI can be
efficiently utilized in various events: 1) enhancing the temporal
resolution of experimental measurements~\cite{Bouhoubeiny2009}, 2)
optimal control~\cite{Ly2001}, 3) reconstruction of incomplete
data~\cite{Druault2007}, 4) multi-dimensional parametric
analysis~\cite{Xiao2015,Fossati2013}, 5) inverse
design~\cite{bui2003proper}, 6) variable fidelity
models~\cite{Mifsud2016}, and more. 

Both DMD and PODI are regarded as data-driven reduced models~\cite{brunton2019data}
since they operate on the snapshots
produced from the FOM, and predict the system dynamics in a
non-intrusive manner. In shape optimization context, DMD and PODI have
been successively applied as
in~\cite{demo2018shape,tezzele2018model,demo2018isope,demo2019marine,tezzele2019marine,Franz2014}
in naval
engineering,~\cite{salmoiraghi2018free,dolci2016proper,LeGuennec2018}
in automotive engineering, or~\cite{ripepi2018reduced,Iuliano2013} in
aeronautics. For acoustic analysis, only very few studies have been
reported in literature. This can possibly include the DMD application
as in~\cite{broatch2019dynamic,Jourdain2013} or POD as
in~\cite{GLEGG2001, Mancinelli2017,Shen2020}.

As previously demonstrated, there are very limited
examples of ROMs specifically tailored for acoustic analysis.
To the best of authors' knowledge, the literature is presently
devoid of documents which characterize, in a thorough and systematic
fashion, the performance of DMD and POD/I techniques on the
hydroacoustic flow fields reconstruction or prediction.

The overarching goal of this
study is to investigate the use of DMD and PODI on a hydroacoustic dataset
corresponding to turbulent incompressible flow past
sphere at $\text{Re}=5000$, and computed using
wall-resolving LES for the fluid dynamic fields, and FWH analogy with
direct integration of the nonlinear quadrupole terms for the acoustic
fields. In particular, the objectives of this work are 1) to understand the effect
of DMD/POD modal truncation on the local and global reconstruction accuracy of the
fluid dynamic fields, 2) to compare the efficiency of DMD and PODI on recovering
the flow and spectral information when half the dataset are utilized,
3) to evaluate the performance of DMD and PODI in terms of data
compression and dipole/quadrupole acoustic prediction.


The present article is organized as follows: first, a brief overview
on the FOM (LES, FWH) and ROM (DMD, PODI) formulations and respective specific works are presented in \autoref{sec:methods}. In
\autoref{sec:results_fom}, the FOM results are presented and then followed by
modal analysis with DMD and POD in \autoref{sec:modal}. The
reconstructed and predicted fluid dynamic data obtained from both ROMs are
discussed in \cref{sec:results_recon,sec:results_pred}, while their spectral and acoustic performances
are addressed in \cref{sec:results_spectral,sec:results_ac}. Conclusions are drawn in
 \autoref{sec:the_end}.

\section{Methodology}
\label{sec:methods}
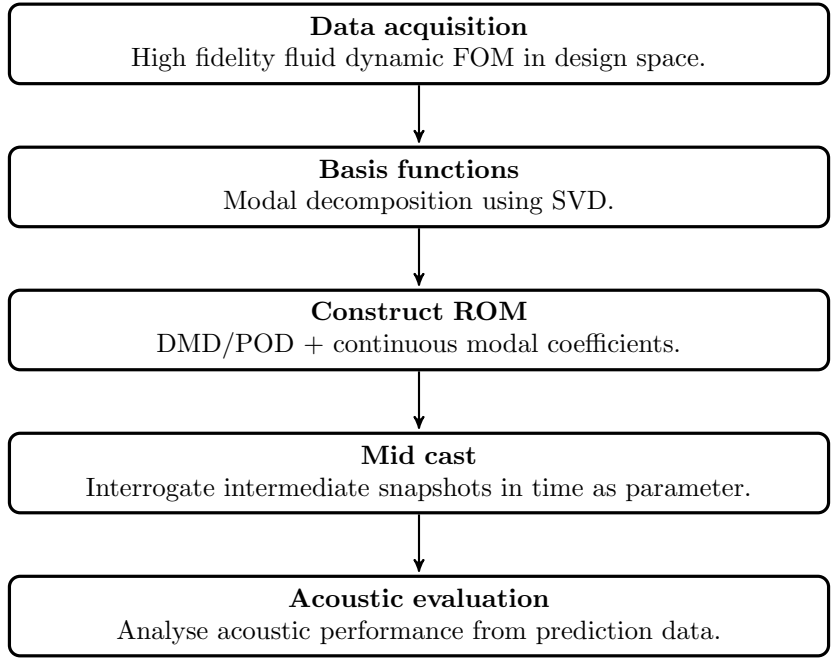
\begin{figure}
\centering
\begin{tikzpicture}
  [node distance=0.8cm, start chain=going below,]
     \node[punktchain, join] (fom) {\textbf{Data acquisition}\\{\text High fidelity fluid dynamic FOM in design space.}};
     \node[punktchain, join] (basis)      {\textbf{Basis functions}\\{\text Modal decomposition using SVD.}};
     \node[punktchain, join] (rom)      {\textbf{Construct ROM}\\{\text DMD/POD + continuous modal coefficients.}};
     \node[punktchain, join] (interpolate) {\textbf{Mid cast}\\{\text Interrogate intermediate snapshots in time as parameter.}};
     \node[punktchain, join] (acoustic) {\textbf{Acoustic evaluation}\\{\text Analyse acoustic performance from prediction data.}};

  \end{tikzpicture}
\caption{Flow chart of the FOM/ROM operations performed in the
  present study} \label{fig:procedure}

\end{figure}

First, the high fidelity data are generated and uniformly sampled in time using high fidelity
simulations with LES turbulence modelling. The resulting
matrix of snapshots is then factorized using Singular Value Decomposition (SVD) which is then used to
construct both the DMD and POD spaces. The scalar coefficients
resulting from the projection of the FOM data onto the POD space are
used to train a continuous representation of the system temporal
dynamics, i.e. PODI approach. In this context, both DMD and PODI are
considered as linear approximation of the dynamical system of the
snapshot matrix and, therefore, are used to predict the data at
intermediate timesteps. The predicted snapshots are finally exploited
to run a post-processing acoustic analogy and validate the accuracy of
the noise generation compared to the FOM data. A summary of the
procedure is presented in the flow chart in \autoref{fig:procedure}. This section is
organized as follows: in \ref{sec:les} the FOM is introduced recalling also the
utilized LES turbulence model, in \ref{sec:math2} the acoustic model used to perform the hydro-acoustic analysis is introduced, in \ref{sec:dmd} and \ref{sec:podi} the non-intrusive pipeline
used to perform model order reduction using DMD and PODI, respectively, is described.

\subsection{Full order model}
\label{sec:les}
The full order model, adopted to provide the snapshots dataset as
input for the reduced order model, is a Large Eddy Simulation.
In LES, the large anisotropic and energy-carrying scales of motion
are directly resolved through an unsteady three dimensional (3D) simulation, whereas
the more isotropic and dissipative small scales of motion are confined
in the sub-grid space. Scale separation is carried out through a
filtering operation of the flow variables.
In literature, the contribution of the Sub-Grid Scales (SGS) of motion
on noise generation and propagation has been found
negligible~\cite{PiomelliCS97,Seror2001}. This means that the
LES model can be considered accurate enough to provide a noise-source flow field, when compared to DNS.
At the same time, the unsteady vortex and coherent structures can be
of extreme importance when computing the noise signature.
Indeed, it has been shown (see among others~\cite{Ianniello2014})
that, in case of complex configurations (e.g. marine propellers),
the RANS methodology may be unsatisfactory in reproducing the flow
fields, adopted as input for the acoustic analysis.

The detailed numerical simulation of the flow around sphere was described in a previous work~\cite{CianferraAI18}, where three different
bluff bodies (a sphere, a cube and a prolate spheroid) were investigated
concerning their noise signature.
\RB{In the mentioned work, validations of the LES solver against experimental and DNS data were presented.
In particular, the flow fields resulting from the employed framework reproduce the thin boundary layer developing along the wall-normal direction, the leading edge, and in the wake. Moreover, they were able to reproduce distributions of the friction coefficients, as well as capturing the underlying flow dynamics with a wide coverage of the spectral content.
}

The filtered Navier-Stokes equations in the incompressible regime are
considered. Within this setting, the SGS stress tensor $\tau^{\text{sgs}}_{ij} = \overline{u_i u_j} -
\bar{u}_i\bar{u}_j$, which represents the effect of the unresolved
fluctuations on the resolved motion, is modeled considering the
Smagorinsky eddy-viscosity closure:
\begin{equation}
\tau^{\text{sgs}}_{ij}-\frac{1}{3} \tau^{\text{sgs}}_{kk} \delta_{ij} = -2 \nu_{t}
\bar{S}_{ij}, \,\,\,\text{with} \,\,\, \bar{S}_{ij} = \frac{1}{2}
\left( \frac{\partial \bar{u}_{j}}{\partial x_{i}} + \frac{\partial
    \bar{u}_{i}}{\partial x_{j}}  \right) ,
\end{equation}
and $\delta_{ij}$ is the isotropic second order tensor, while the overbar denotes the filtering operator.
The SGS eddy viscosity $\nu_t$  is expressed as:
\begin{equation}
 \nu_t = (C_s \Delta)^2 \, |\bar{S}_{ij}|,
\end{equation}
where  \RB{$\Delta=\sqrt[3]{V_c}$ is the
filter width which is defined as the cubic root of the cell volume $V_c$}, and the Smagorinsky constant $C_s$ is computed dynamically using the
Lagrangian procedure of~\cite{MeneveauLC96}, averaging over the
fluid-particle Lagrangian trajectories.

\subsection{Acoustic model}
\label{sec:math2}
The acoustic model herein considered is the one proposed by Ffowcs Williams and
Hawings~\cite{fwh1969}, which is an extension of the
Lighthill theory.

The basic idea behind the acoustic analogies is that pressure perturbation
originates in the flow field and propagates in the far-field where the
medium is assumed quiescent and uniform.
The integral solution of the acoustic wave equation presented by
Ffowcs Williams and Hawkings consists of surface and volume integrals,
meaning that the sources of fluid-dynamic noise can be found as
pressure-velocity fluctuations developing in the fluid region or as
reflected pressure on an immersed solid surface.

We consider the original formulation presented in~\cite{fwh1969}, and modified according to the works of Najafi et
al.~\cite{NajafiBM11} and Cianferra et al.~\cite{CianferraAI18}. The
modification of the original FWH equation takes into account the
advection of acoustic waves.
To account for the surrounding fluid moving at a constant speed (along
the $x$ axis), the advective form of the Green's function must be
considered.
A derivation of the advective FWH equation is reported
in~\cite{NajafiBM11}, where the authors developed an integral solving
formulation for the linear (surface) terms.
The advective formulation of the volume term for the particular case
of the wind tunnel flow is reported in~\cite{CianferraAI18}.

In the present work, as done for the fluid dynamic part, we
describe the formulation without dwelling into details, for which we
refer to the previous works~\cite{CianferraAI18,CianferraIA19}.

The acoustic pressure $\hat{p}$, at any point $\textbf{x}$ and time
$t$, is represented by the sum of surface ($\hat{p}_{\text{2D}}$) and volume
($\hat{p}_{\text{3D}}$) integrals, respectively:
\begin{equation} \label{eq:FWH2D}
4 \pi \hat{p}_{\text{2D}}(\textbf{x},t) \, = \frac{1}{c_{0}}
  \frac{\partial}{\partial t} \int_{S} \left[ \frac{\tilde{p}
      \hat{n}_i \hat{r}_i}{r^*} \right]_{\tau} dS
  + \int_{S} \left[ \frac{\tilde{p} \hat{n}_i
      \hat{r}_i^{*}}{{r^{*}}^2} \right]_{\tau} dS ,
\end{equation}
\begin{equation}\label{eq:FWH3D}
  \begin{aligned}
  4 \pi \hat{p}_{\text{3D}} (\textbf{x}, t) &= \frac{1}{c_{0}^{2}}\frac{\partial^{2}}{\partial t^{2}}\int_{f > 0} \left\{ T_{ij} \left[ \frac{\hat{r}_i \hat{r}_j}{r^{*}} \right] \right\}_{\tau} dV  \\
   & + \frac{1}{c_{0}}\frac{\partial}{\partial t} \int_{f>0} \left\{
       T_{ij} \left[ \frac{ 2 \hat{r}_i \hat{r}^{*}_{j}}{{r^{*}}^2} +
       \frac{\hat{r}^{*}_i \hat{r}^{*}_{j} - R^{*}_{ij}}{\beta^2
       {r^{*}}^2} \right] \right\}_{\tau} dV  \\
   &+ \int_{f>0} \left\{  T_{ij} \left[\frac{3 \hat{r}^{*}_i
       \hat{r}^{*}_j - R^{*}_{ij}}{{r^{*}}^3}  \right] \right\}_{\tau}
       dV.
\end{aligned}
\end{equation}

The pressure perturbation with respect to the
reference value $p_0$ is denoted with $\tilde{p} = p - p_0$, $\hat{n}$
is the (outward) unit normal vector to the surface element $dS$, and
$c_{0}$ is the sound speed. $\hat r$ and $\hat r^*$ are unit radiation
vectors, $r$ and
$r^{*}$ are the module of the radiation vectors $\textbf{r}$ and
$\textbf{r}^{*}$ respectively. Their description is given in detail
in~\cite{CianferraAI18}.

Equation~\eqref{eq:FWH3D} contains two second--order tensors:
$R^*_{ij}$ and the Lighthill stress tensor $T_{ij}$,  the latter
characterizing the FWH quadrupole term. Under the assumption of
negligible viscous effects and iso-entropic transformations for the
fluid in the acoustic field, the Lighthill tensor reads as:
\begin{equation}
T_{ij} = \rho_{0} u_{i} u_{j} + \left( \tilde{p} - c_0^2 \tilde{\rho} \right) \delta_{ij},
\end{equation}
where $\tilde{\rho}$ is the density perturbation of the flow which, in
our case, is equal to zero. The surface integrals in Equation~\eqref{eq:FWH2D}
are referred to as linear terms of the FWH equation and represent the
loading noise term. The volume integrals in Equation~\eqref{eq:FWH3D}
are slightly different from the standard FWH (non-advective)
equation. For their derivation we consider a uniform flow with
velocity $U_0$ along the streamwise direction.

Obviously, the direct integration of the volume terms gives accurate
results, however, this method can be used if the calculation of the
time delays can be omitted, otherwise the computational burden makes
it unfeasible. In fact, the calculation of the time delays requires storing at each
time step the pressure and velocity data related to the entire
(noise-source) volume, in order to perform an interpolation over all
the data.

However, for the case herein investigated, the evaluation of the non-dimensional Maximum Frequency Parameter (MFP)~\cite{CianferraIA19} (which is greater than unity for
every microphone considered)
allows to adopt the assumption of compact noise source.
This means that, in the investigated case, the time delay is very small
and the composition of the signals is not expected to contribute to
the radiated noise.
Since the evaluation of the time delays may be reasonably omitted,
a remarkable saving of the CPU time is achieved, and the direct computation of the quadrupole volume terms becomes feasible.

In the two previous works~\cite{CianferraAI18,CianferraIA19},
as to perform a validation test for the acoustic model, the solution
of the advective FWH equation was compared with the pressure signal
provided by LES, considered as reference data.  This comparison is
useful to verify the ability of the acoustic post processing to
accurately reconstruct the pressure field. Also, it points out the
frequency range which is important to consider.

\subsection{Dynamic mode decomposition}
\label{sec:dmd}

Dynamic Mode Decomposition is a data-driven modal
decomposition technique for analysing the dynamics of nonlinear
systems~\cite{schmid2011application,schmid2011applications}. DMD was
developed in~\cite{schmid2010dynamic}, and since then it has \RA{attracted a lot of
attention} and lead to diverse applications. Its popularity is
also due to its equation-free nature since it relies only on snapshots
of the state of the system at given times. DMD is able to identify
spatiotemporal coherent structures, that evolve linearly in time, in
order to approximate nonlinear time-dependent systems. For a
comprehensive overview we refer to~\cite{kutz2016dynamic}.

In the last years, many variants of the classical algorithm arose such
as, for instance, the multiresolution DMD~\cite{kutz2016multiresolution}, DMD with
control~\cite{proctor2016dynamic}, and compressed
DMD~\cite{erichson2016compressed}. Higher order DMD~\cite{le2017higher} was developed for the cases that show limited spatial complexity but a very
large number of involved frequencies. Sparsity promoting DMD to choose dynamically
important DMD modes, can be found in~\cite{Jovanovi2014}. Kernel DMD
was proposed in~\cite{williams2015data}, where a large number of observables were  able to be incorporated into DMD, exploiting the
kernel trick. Randomized DMD for non-intrusive ROMs can be found
in~\cite{bistrian2017randomized,bistrian2018efficiency}.  We also cite
a new paradigm for data-driven modelling
that uses deep learning and DMD for signal-noise
decomposition~\cite{rudy2018deep}. For
numerical non-intrusive pipelines in applied sciences and industrial
application we mention~\cite{rozza2018advances,tezzele2018ecmi,tezzele2019mortech}. From a
practical point of view, the Python package PyDMD~\cite{demo18pydmd}
contains all the major and most used aforementioned versions stemmed
from the classical DMD in a user friendly way, from which the present
work also exploits.

Here we present a brief overview of the standard algorithm, suppose a
given discrete dataset $\bm x_k \in \mathbb{R}^n$, where $k \in [1,...,m]$
denotes the snapshot number, and $n$ denotes the number of degrees of
freedom. The $m$ snapshots can be arranged by columns into the data
matrix $\bm S = \big[ \bm x_1 , \bm x_2 , \dots , \bm
  x_{m-1} \big]$ and its temporal evolution $\dot{\bm S} = \big[
  \bm x_2 , \bm x_3 , \dots , \bm x_{m} \big]$. Now, the objective is
  to seek a Koopman-like operator~\cite{rowley2009spectral} $\bm A$,
  which is a best fit linear approximation that minimizes $\| \dot{\bm
    S} - \bm A \bm S \|_F$. Such operator is given as $\bm A =
  \dot{\bm S} \bm S^\dagger \in \mathbb{R}^{n\times n}$, where
  $^\dagger$ is the Moore-Penrose pseudoinverse.

For practical reasons, the matrix $\bm A$ is not solved explicitly
since it contains $n^2$ elements, from which $n$ is typically of
several order of magnitudes, and is usually much larger than the
number of snapshots $m$. Instead, a reduced operator $\widetilde{\bm
  A}$ is considered which has much lower dimension, yet it preserves
the spectral information of $\bm A$. The procedure is defined in the
following. First, a singular value decomposition (SVD) is performed on
the snapshots matrix,
\begin{equation}
\label{eq:dmd_svd}
\bm{S} = \bm{U} \boldsymbol{\Sigma} \bm{V}^* \approx
\bm{U}_r \boldsymbol{\Sigma}_r \mathbf{V}^*_r,
\end{equation}
where $r \leq (m-1)$ defines the SVD
truncation rank, $\bm U_r$ is the truncated POD modes as it will be further discussed in the following section, and the three terms $\bm U_r \in \mathbb{C}^{n \times r}$, $\bm \Sigma_r \in
\mathbb{C}^{r \times r}$ and $\bm V^*_r \in \mathbb{C}^{r \times (m-1)}$
denote the truncated SVD components. Second, the
reduced operator $\widetilde{\bm A}$ is
obtained through an $r \times r$ projection of $\bm A$ onto the POD
modes $\bm U_r$, that is
\begin{equation}
\widetilde{\bm A} = \bm U^*_r \bm A \bm U_r = \bm U^*_r \dot{\bm S} \bm V_r
\bm \Sigma^{-1}_r  \in \mathbb{C}^{r \times r}.
\end{equation}

The eigenspace of $\bm A$ is denoted as $\bm \Phi$ (also called the DMD
modes), and it can be revealed through the low-rank projection of the
eigenvectors $\bm W$ of $\widetilde{\bm A}$ onto the POD modes,
\begin{equation}
\bm \Phi = \bm U_r \bm W ,
\end{equation}
where $\bm W$ is obtained from the eigendecomposition of $\widetilde{\bm A}$, i.e.
\begin{equation}
\widetilde{\bm A} \bm W = \bm W \bm \Lambda ,
\end{equation}
and $\bm \Lambda$ is the diagonal matrix of eigenvalues of
$\widetilde{\bm A}$ which is also considered the same matrix obtained from
the eigendecomposition of ${\bm A}$. Finally, the DMD modes, $\bm
\Phi$, and the diagonal matrix of eigenvalues, $\bm \Lambda$, are both
used to provide a linear approximation to the solution vector $\bm
x(t)$ at any time instance, that is
\begin{equation}
\bm x(t) \approx \bm \Phi \exp(\bm \Lambda t) \bm \Phi^\dagger \bm x(0) .
\end{equation}

In this work, we are going to characterize the DMD modes for velocity
and pressure fields, and exploit them to make predictions between time
instants. The reconstructed fields will be used to approximate the
hydroacoustic noise.

\subsection{Proper orthogonal decomposition with interpolation}
\label{sec:podi}

Proper orthogonal decomposition is a linear dimensionality reduction technique
that is widely used to identify the underlying structures within large datasets. In general,
a linear dimensionality reduction assumes that each snapshot in the dataset can be generated as a linear
combination of a properly chosen small set of basis functions. In this regards, similar to the previous section, we
introduce the snapshots matrix $\mathbf{S} \in \mathbb{R}^{n \times m}$ which now contains $m$ snapshots,
thus $\mathbf{S} = [\mathbf{x}_1 \, \mathbf{x}_2 \, \cdots \, \mathbf{x}_m]$. Then, we assume it exists
the linear combination $\mathbf{S} \approx \mathbf{U}_r \mathbf{C}$,
where the matrix $\mathbf{U}_r \in \mathbb{R}^{n \times r}$ contains the
vectors defining the basis functions, $r \leq m$, and $\mathbf{C} \in
\mathbb{R}^{r \times m}$ is the matrix containing the linear
coefficients. This approximation is equivalent to low-rank matrix
approximation~\cite{markovsky2019low}.

If we choose the basis functions that minimize the sum of the squares
of the residual $(\mathbf{S} - \mathbf{U}_r \mathbf{C})$, in the least
square sense, then we recover the Proper Orthogonal Decomposition.
In such case, the matrix $\mathbf{U}_r$ contains the so-called POD
modes, and they can be recovered with the SVD of the snapshots matrix, as defined in \autoref{eq:dmd_svd}, for a low-rank truncation that
retain the first $r$ modes. In fact, the error
introduced by the truncation can measured as~\cite{quarteroni2015reduced}:
\begin{align}
\| \mathbf{S} - \mathbf{U}_r \boldsymbol{\Sigma}_r \mathbf{V}^*_r
  \|_2^2 &= \lambda_{r+1}^2, \\
\| \mathbf{S} - \mathbf{U}_r \boldsymbol{\Sigma}_r \mathbf{V}^*_r \|_F
&= \sqrt{\sum_{i=r+1}^m \lambda_i^2},
\end{align}
where the subscripts 2 and $F$ refer to the Euclidean norm and to the Frobenius
norm respectively, and $\lambda_i$ are the singular values arranged in descending order.

To approximate the solution manifold, after the projection of the
original snapshots onto the POD space spanned by the POD modes, we
interpolate the modal coefficients for new input parameters.

This non-intrusive data-driven approach is called POD with
interpolation~\cite{Ly2001,bui2003proper,bui2004aerodynamic}. A coupling with isogeometric
analysis can be found in~\cite{garotta2018quiet}. For an enhancement
of the method using active subspaces property, in the case of small
computational budget devoted to the offline phase, see~\cite{demo2019cras}.
An open source Python implementation of PODI can be found in the EZyRB
package~\cite{demo2018ezyrb}. For an implementation based on OpenFOAM~\cite{of}
see the freely available ITHACA-FV library~\cite{stabile_stabilized} from which the present work exploits.

The PODI method can be decomposed in two distinct phases: an offline phase,
and an online one. In the first phase, given the solutions snapshots for
given parameters $\mu_k$, we compute the POD modes $\phi_k$ (which construct the matrix $\mathbf{U}$) through SVD, and the
corresponding modal coefficients for all the original snapshot with a
projection. This results in
\begin{equation}
\mathbf{x}_k = \sum_{j=1}^m \alpha_{jk}\phi_j \approx \sum_{j=1}^r
\alpha_{jk}\phi_j,\quad\forall k \in [1, 2, \dotsc, m],
\label{eq:pod}
\end{equation}
where $\alpha_{jk}$ are the elements of the coefficient matrix
$\mathbf{C}$.

After performing the offline phase, the $r$ most energetic modes are now obtained, and the corresponding
modal coefficients $\mathbf \alpha_k$ are used to construct the reduced space. Using the pairs
$(\mu_k, \mathbf \alpha_k)$ we are able to reconstruct the function that
maps the input parameters to the modal coefficients. This function is
interpolatory for the snapshots used in the training phase, hence the name \emph{POD with interpolation}. With this surrogate function, we are able
to reconstruct in real-time the solution fields of interest for an arbitrary
parameter in the so-called online phase.
In the present study, the input parameter is the time while the output fields are the velocity and pressure.
Also, a cubic spline interpolation is considered for the modal coefficients.
\RA{
The degree of the interpolator was rather a modelling choice, given the flexibility provided by the simple, one-dimensional coefficient samples.
In addition, the fairly easy implementation allowed for a robust performance that is comparative to simpler models, while smoother time response is herein preserved due to the higher order polynomial interpolation.
Indeed, a more complicated application involving higher dimensional parameter space would require a careful choice of the hypothesis space, and typically more advanced machine learning (or deep learning) techniques would be more relevant, which is not the case in the present study.
}

\begin{figure}
\centering
\includegraphics[width=0.9\textwidth]{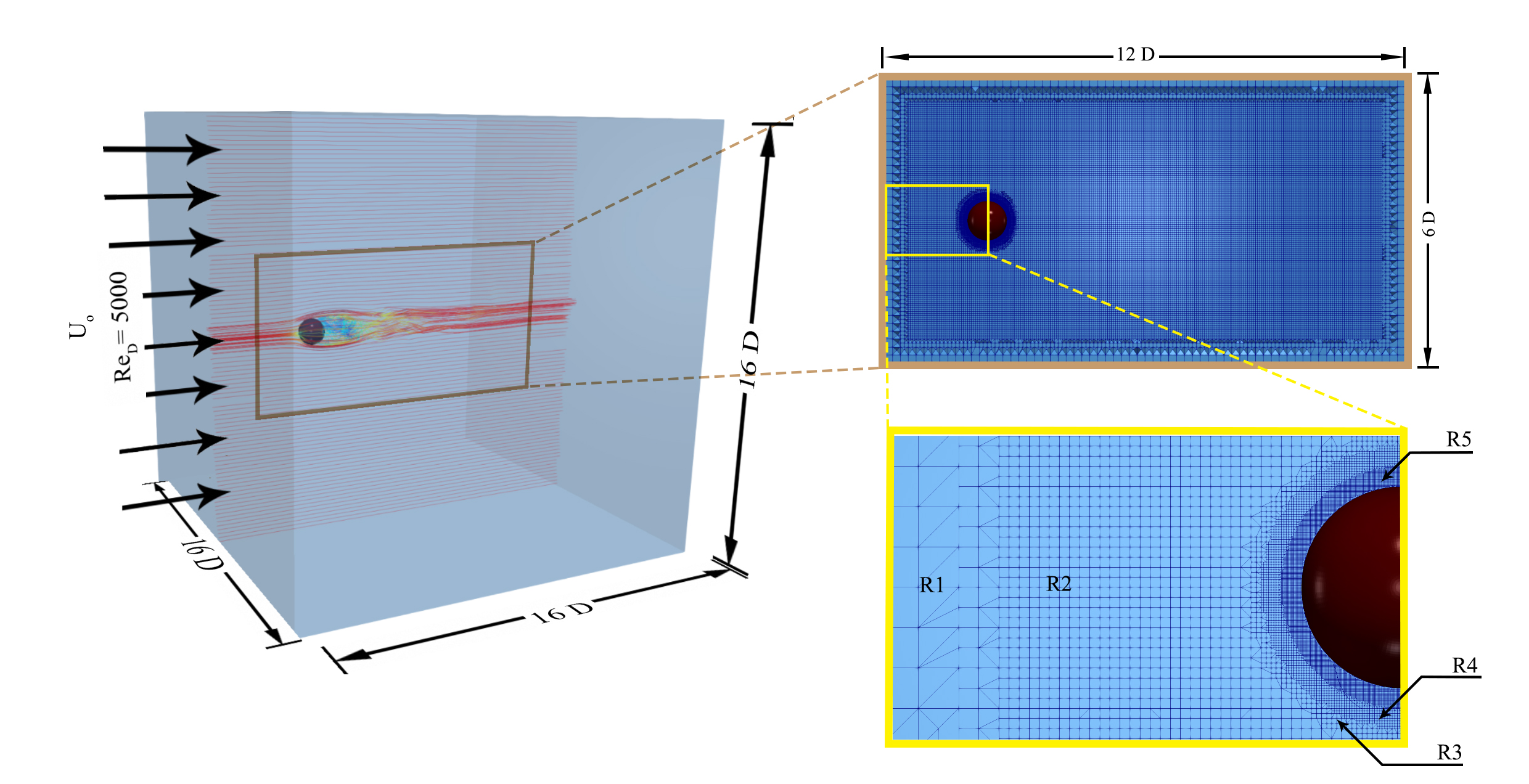}
\caption{Sketch of the computational domain used for the FOM. The sphere diameter is $D=0.01$~\si{m}. Successive mesh refinement layers ($R2$, $R3$, $R4$) are performed through cell splitting approach until reaching the finest grid spacing $0.001D$ in the region $R5$.}
\label{fig:mesh}
\end{figure}

In literature, several works have been devoted to study the accuracy
of the projection coefficients manifold construction and
performance. In~\cite{Mifsud2009} they explored the full factorial and
latin hyper-cube sampling techniques of the parameter space. The
manifold representation was also studied
by~\cite{Karri2009,Mifsud2009} mainly comparing between linear
regression, polynomial, spline, finite difference, and radial basis
functions (RBF). The high-order singular value decomposition (HOSVD)
was demonstrated in~\cite{Lorente2008} as an efficient method for
flows with strong discontinuities. Kriging interpolation was utilized
in~\cite{Fossati2013} for $n$-dimensional approximation of potentially
complex surfaces.

\section{Numerical results}
\label{sec:results}

\subsection{Full order CFD}
\label{sec:results_fd}
\label{sec:results_fom}

The fluid dynamic fields are solved in the framework of the
OpenFOAM library~\cite{of} which is based on the Finite
Volume Method (FVM). The filtered Navier-Stokes equations are
solved using the PISO pressure-velocity coupling algorithm implemented
in the pisoFoam solver. The spatial derivatives are discretized through
second-order central differences.
Implicit time advancement runs according to the Euler scheme.
The numerical algorithm, including the SGS closure, has been customized at the laboratory of
Industrial and Environmental Fluid Mechanics (IE-Fluids) of the
University of Trieste, and more details can be found in~\cite{CintolesiPA15}.

The fluid dynamic full order model simulates a sphere of diameter $D=0.01$~\si{m}, immersed in a water
stream with constant streamwise velocity $U_0=1$~\si{m/s}. The kinematic viscosity is
$\nu=2.0 \times 10^{-6}$~\si{m^2/s}, so that the Reynolds number based on
the sphere diameter is $\text{Re}_D=5000$.

The computational domain, depicted in \autoref{fig:mesh}, is a box with dimensions $16D \times 16D
\times 16D$ along the $x$, $y$ and $z$ axes respectively. The sphere is located such that a distance
of $12D$ is attained downstream, along the $x$-axis, while it is centered with respect to the other axes.
A zero-gradient
condition is set for the pressure at the domain boundaries, except for the outlet where
pressure is set to zero. The velocity is set to $U_0$ at the inlet, stress-free
condition is set at the lateral boundaries, and zero-gradient condition
is set for the velocity components at the outlet.

The grid, unstructured, and body-fitted, consists of about $5$ millions
of cells. It is created using the OpenFOAM \textit{snappyHexMesh} utility.
The grid spacing normal to the wall for the densest
layer of cells (indicated as $R5$) is such to have first cell center within a
wall unit $y^+$ ($y^+=u_\tau y/\nu$ with $u_\tau= \sqrt{\tau_w/\rho_0}$ and $\tau_w$
the mean shear stress).
An \textit{A posteriori} analysis showed that about $5$ grid points are placed within
$10$ wall units off the wall. The grid spacing is obtained through successive
transition refinements (indicated as $R3$ and $R4$ in \autoref{fig:mesh}).
A refinement box around the body (named $R2$ in \autoref{fig:mesh}) is considered
so as to obtain, in the wake region, a grid size of less than $0.1D$ at a distance of $8D$.
Out of the region of interest, a coarser grid (indicated as $R1$) allows for possible extension of the domain dimensions, and reducing possible
disturbance effects coming from the boundaries.

For the time integration, a constant time step is set to \RB{$\Delta t = 10^{-5}$~\si{sec}}
in order to keep the Courant number under the threshold of $0.5$. The flow
around the sphere is completely developed after about 80 characteristic times $D/U_0$.

\subsection{Modal decomposition}
\label{sec:modal}
In this section, we report the numerical results concerning the modal decomposition
of the full order snapshots. This phase is particularly useful in order to have an insight onto the
dominant structures and on the frequencies hidden in the full order dynamical system.
We report both an analysis on the eigenvalue decay which is associated with the Kolmogorov width,
and the modal representation which permits to visualize the turbulent structures associated
with each mode. Moreover, we analyse the time evolution
of the temporal coefficients in order to identify the time frequencies associated with each mode.
All the computations have been carried out using the
PyDMD Python package~\cite{demo18pydmd} and
the ITHACA-FV library~\cite{stabile_stabilized, StaHiMoLoRo2017}.

\subsubsection{Singular values decay}
\begin{figure}
\centering
\includegraphics[width=0.6\textwidth]{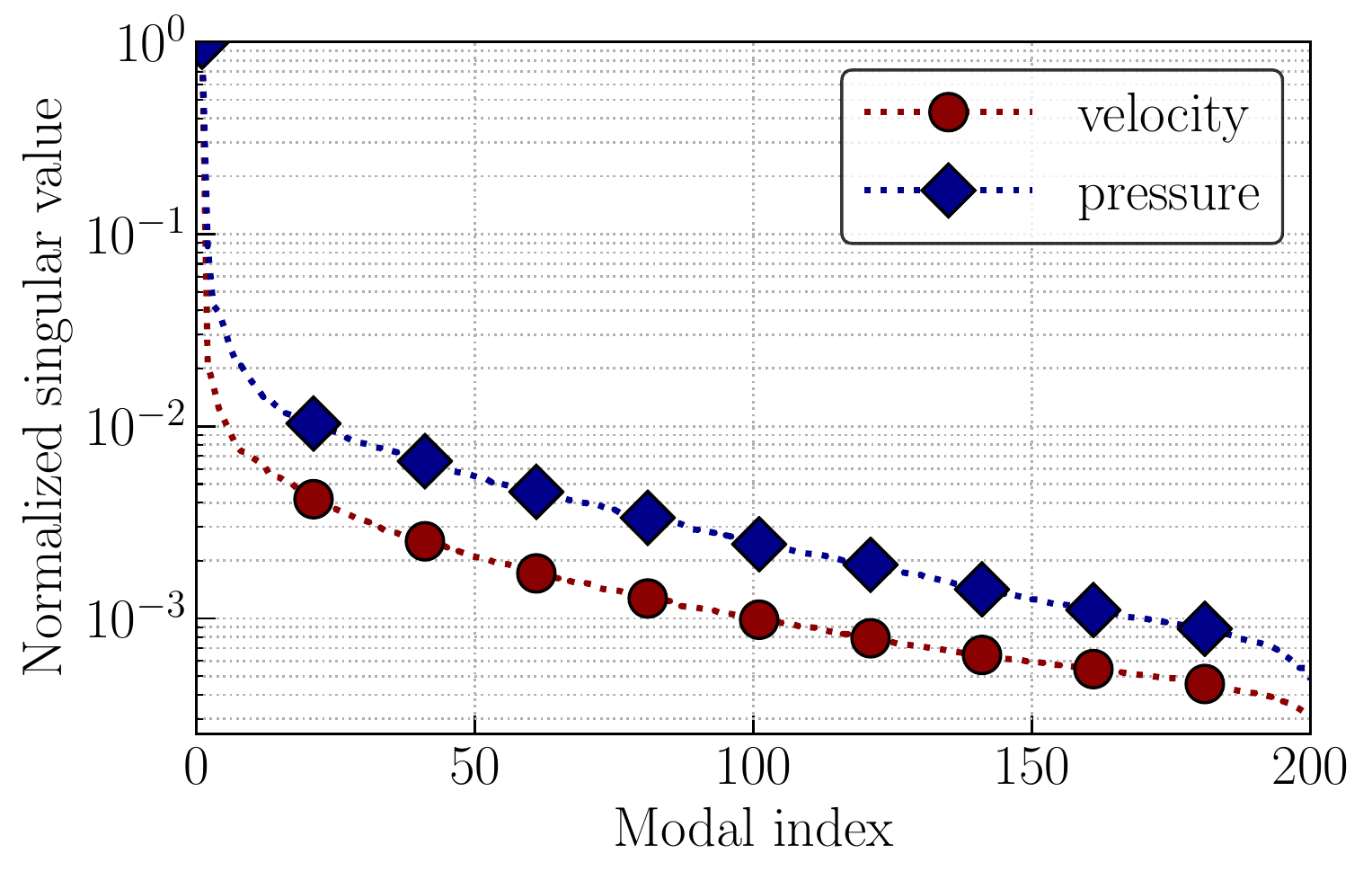}
\caption{Normalized singular values (SV) for streamwise
  velocity and pressure snapshots. SV are arranged in descending order. The
  absence of any gaps in the plots suggests no specific truncation
  rank for the SV-based reduced models.}
\label{fig:svd}
\end{figure}

\autoref{fig:svd} depicts the first step of the modal analysis applied to the
Navier--Stokes fluid dynamic problem considered. The plot shows the magnitude of
the normalized singular values (SV) obtained from the SVD factorization of the
snapshot matrix, obtained both for the streamwise velocity field component
and for the pressure field. In the diagram --- and in all the following discussion ---
the modes are arranged in descending order according to the corresponding SV magnitude.
Typically, a presence of SV magnitude gaps in such plot provides an indication of a
convenient truncation rank for the modal analysis. In the present case, a steady and
continuous decay is observed after a steep slope corresponding to the first $15$ to
$20$ modes. Thus, the absence of SV magnitude gaps in the higher frequencies region suggests that, for both fluid dynamic
fields considered, there is no specific truncation rank for the modal analysis.

\subsubsection{Modal representation}
As shown, the SV magnitude observation is not resulting in an obvious indication
of the modal truncation rank. To start understanding the effect of modal truncation rank
on the fluid dynamic solution accuracy, we then resort to considerations based on
the spatial and time frequencies which need to be reproduced in the hydro-acoustic simulations.

\begin{figure*}
\centering
\includegraphics[width=0.8\textwidth]{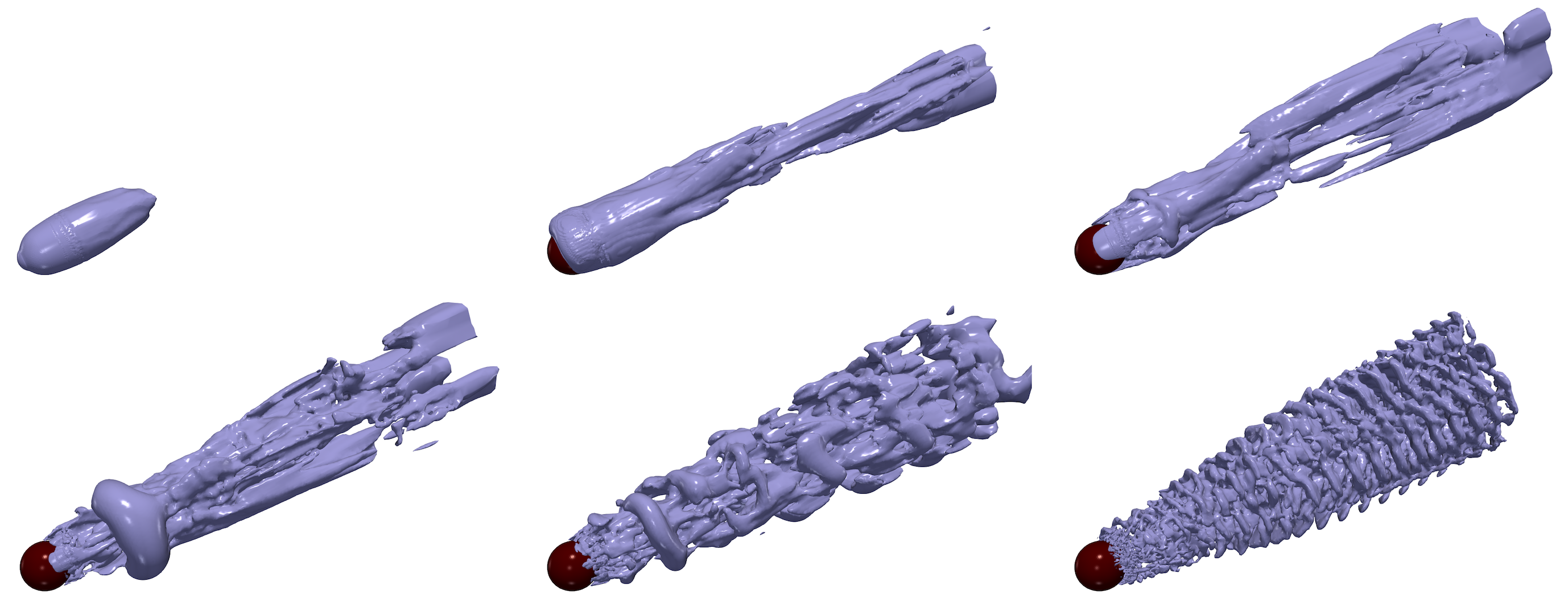}
\put(-420,170){\color{black}\fontsize{8}{8}\fbox{$\text{Velocity (DMD)}$}}

\includegraphics[width=0.8\textwidth]{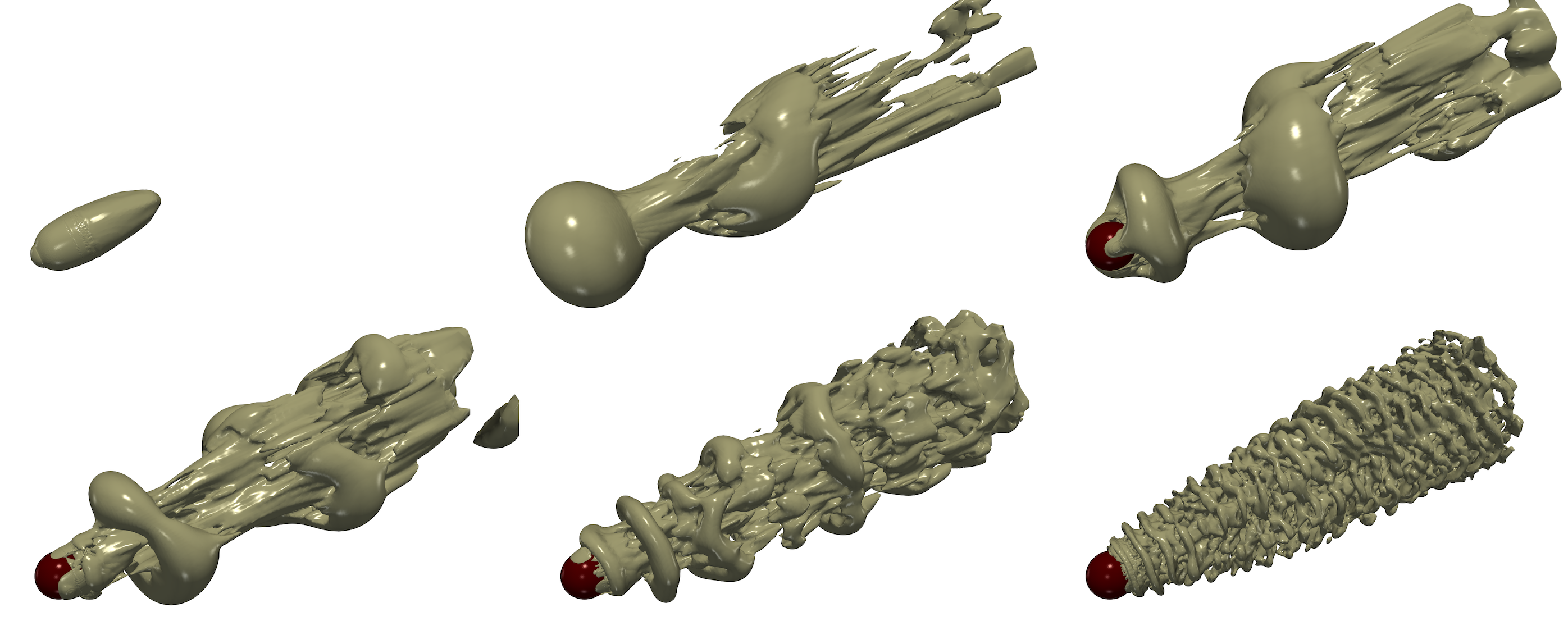}
\put(-420,170){\color{black}\fontsize{8}{8}\fbox{$\text{Velocity (POD)}$}}
\caption{Isosurfaces of the DMD modes (purple color) at value of
  $-15\times 10^{-5}$~\si{m/s} versus the isosurfaces of
  POD modes (olive color) at value of $6.5$~\si{m/s} for
  the streamwise velocity field. Modes: $1$, $2$, $4$, $8$, $36$, and $128$.}
\label{fig:Umodes}
\end{figure*}

\autoref{fig:Umodes} depicts a set of three dimensional modal shapes
resulting from the longitudinal velocity field DMD and POD modal decomposition, respectively.
The purple diagrams represent isosurfaces passing through the field data points of the DMD modes at value $-15\times 10^{-5}$ \si{m/s}, while
the olive color plots refer to isosurfaces of POD modes at value $6.5$ \si{m/s}. For each modal
decomposition methodology, the images in the Figure are arranged in tabular fashion and refer
to modes $1$, $2$, $4$ on the first row, and $8$, $36$, $128$ on the second one.
The plots suggest that the DMD modal shapes present a more pronounced tendency to be organized according to
spatial frequencies with respect to the POD modes. In fact wider turbulent structures are only found in
the very first DMD modes, while in the case of POD modes they can be identified also in higher rank
modes, among higher frequency patterns. Along with this tendency, the turbulent structures associated
with the low rank DMD modes also appear organized in longitudinal streaks, and gradually become more
isotropic for higher rank modes. It should be pointed out though that despite the fact that POD modes
are in general not designed to separate the contributes of single harmonic components, the modal shapes
obtained for the longitudinal velocity do appear to be at least qualitatively correlated to spatial
frequencies. In fact, by a qualitative standpoint, the spatial frequencies appearing in the plots
corresponding to higher modes are in general higher with respect to those associated to the first modes.

\begin{figure*}
\centering
\includegraphics[width=0.8\textwidth]{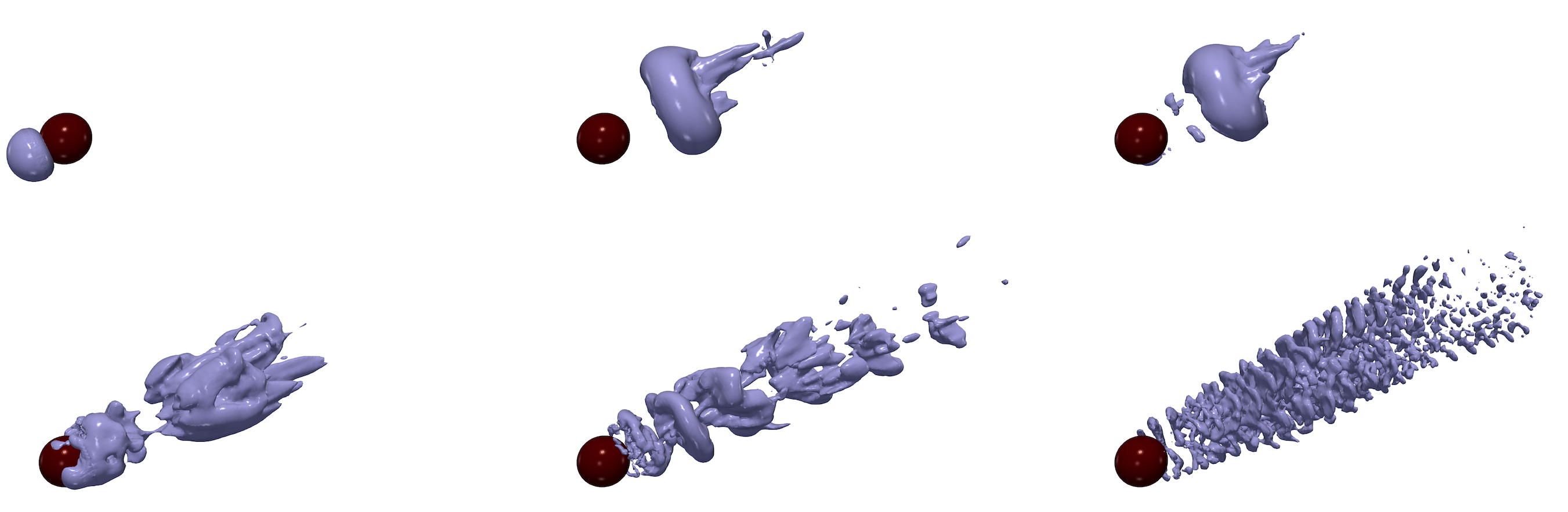}
\put(-420,150){\color{black}\fontsize{8}{8}\fbox{$\text{Pressure (DMD)}$}}

\includegraphics[width=0.8\textwidth]{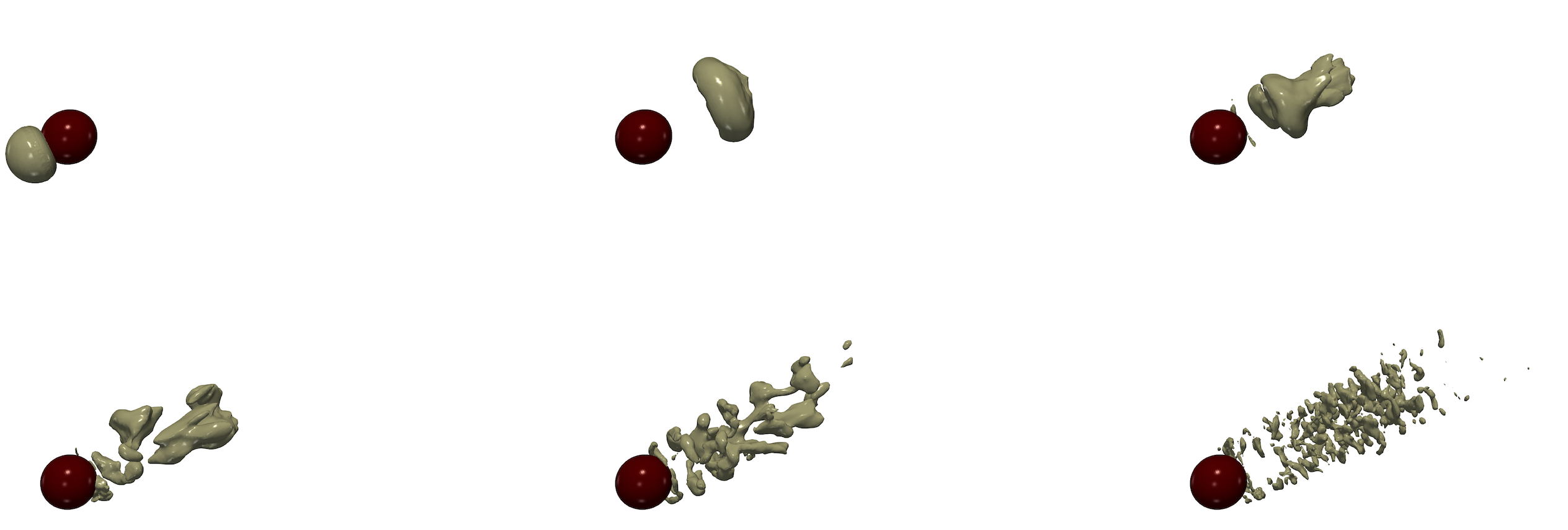}
\put(-420,150){\color{black}\fontsize{8}{8}\fbox{$\text{Pressure (POD)}$}}
\caption{Isosurfaces of the DMD modes (purple color) at value of
  $-5\times 10^{-4}$~\si{Pa} versus the isosurfaces of
  POD modes (olive color) at value of $300$~\si{Pa} for
  the pressure field. Modes: $1$, $2$, $4$, $8$, $36$, and $128$.}
\label{fig:Pmodes}
\end{figure*}

Similar considerations can been drawn from the observation of similar plots corresponding to
the modal decomposition of the pressure fields, presented in \autoref{fig:Pmodes}.
Also in this case, the purple plots refer to isosurfaces passing through all data points of value $-5\times 10^{-4}$ \si{m^2/s^2}
of DMD modal shape functions, while the olive diagrams refer to isosurfaces of POD modes at the
value $300$ \si{m^2/s^2}. The plots are again arranged, for each decomposition methodology considered, in
tabular fashion portraying modes $1$, $2$, $4$ on the first row, and $8$, $36$, $128$ on the second one.
The pressure modes associated with both methods seem again qualitatively arranged according to spatial
frequency content. As previously observed for the longitudinal velocity modes, the DMD modes
appear more closely correlated to spatial frequencies, and present a less isotropic appearance
with respect to their POD counterparts.

\subsubsection{Associated coefficients}
After having characterized the spatial frequency content of different
POD and DMD modes, we now want to analyse the time frequencies associated to
each mode. To do this, we decompose each snapshot into its modal components
and observe the time evolution of the modal coefficients. In fact, as the POD
and DMD modes generated from the snapshots are constant in time, the corresponding
modal coefficients must depend on time to allow for the reconstructed solution
to reproduce the correct time variation.

The four plots presented in \autoref{fig:coeffs} show, the temporal evolution of the
modal coefficients associated, respectively, to the DMD decomposition of the longitudinal
velocity field (top left), to the POD decomposition of the longitudinal velocity field
(top right), to the DMD decomposition of the pressure field (bottom left) and to the POD
decomposition of the pressure field (bottom right). Each diagram reports four lines referring
to the coefficients of modes $2$, $4$, $8$, $36$.
\begin{figure}
\centering
\includegraphics[width=0.9\textwidth]{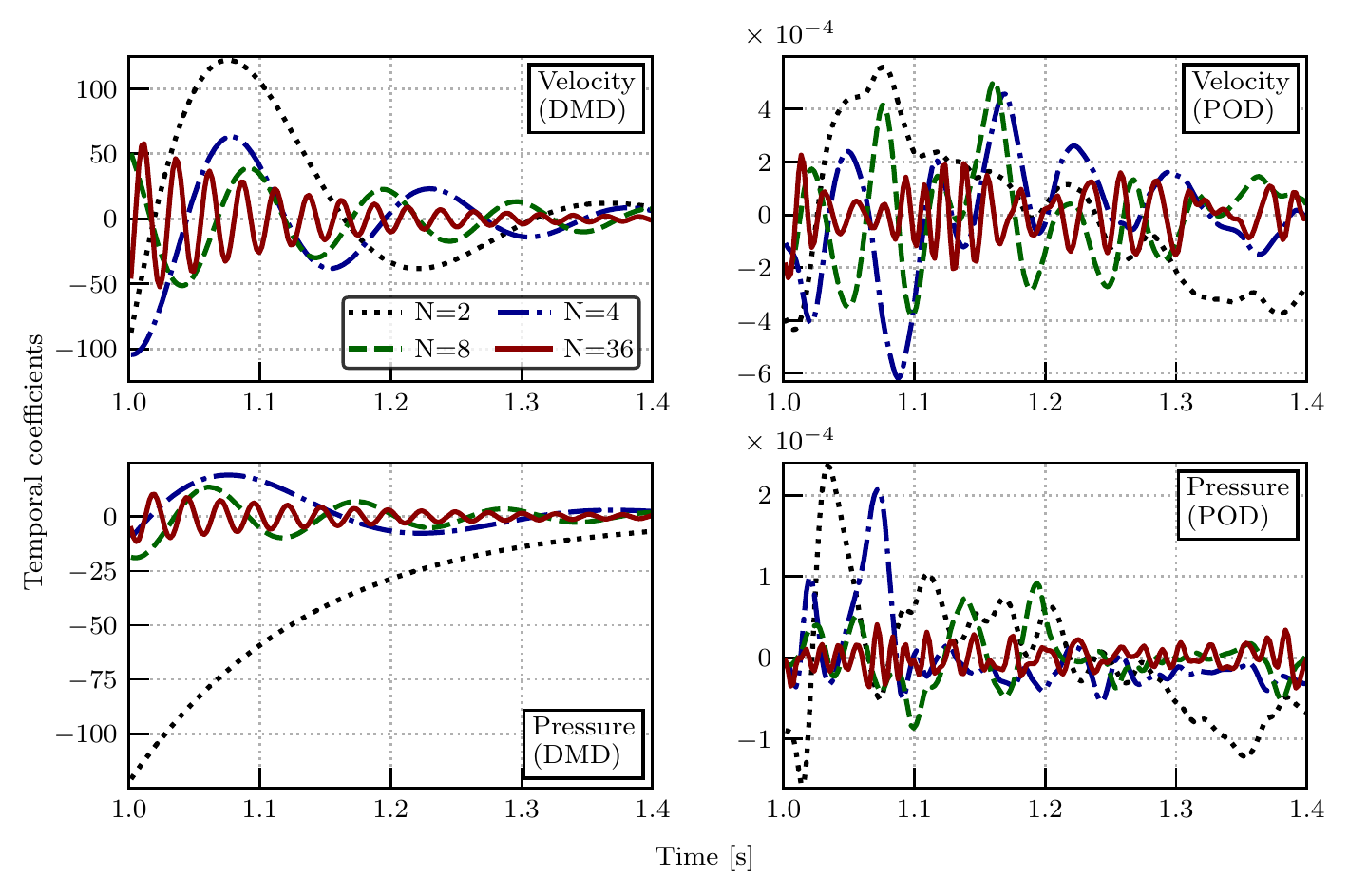}
\caption{Temporal coefficients associated to the DMD and POD (left to
  right) modes number $2$, $4$, $8$, and $36$ of the streamwise velocity and
  pressure data (top to bottom).}
\label{fig:coeffs}
\end{figure}

By a qualitative perspective, the diagrams in \autoref{fig:coeffs} suggest that DMD and POD
modal coefficients are strongly correlated with the time frequencies. In fact, higher frequency harmonics
appear in the time evolution of higher order modal coefficients, which are not observed in lower ones.
As expected, also in this case the frequency-mode association is definitely stronger for DMD modes,
in which a single dominant harmonic can be identified in correspondence with each modal coefficient. As for
POD, the respective frequency content seems to cover wider set of harmonics, associated
with higher frequencies as higher modes are utilized.


\subsection{Flow fields reconstruction}
\label{sec:results_recon}
A first aim of this work is to assess whether the proposed DMD
and POD algorithms, are able to accurately reproduce the full order
model solutions. A first step in such assessment will be that of
checking the effectiveness of the SVD based modal decomposition strategies
of the model reduction algorithm considered. In particular, the reconstructed
solution convergence to the snapshots considered will
be discussed both through the visualization of single snapshots flow fields and
by presenting convergence plots of error averaged among snapshots. Finally,
we will present similar plots for the solution predicted by means of both
DMD and PODI.

\begin{figure*}
\centering
\includegraphics[width=0.9\textwidth]{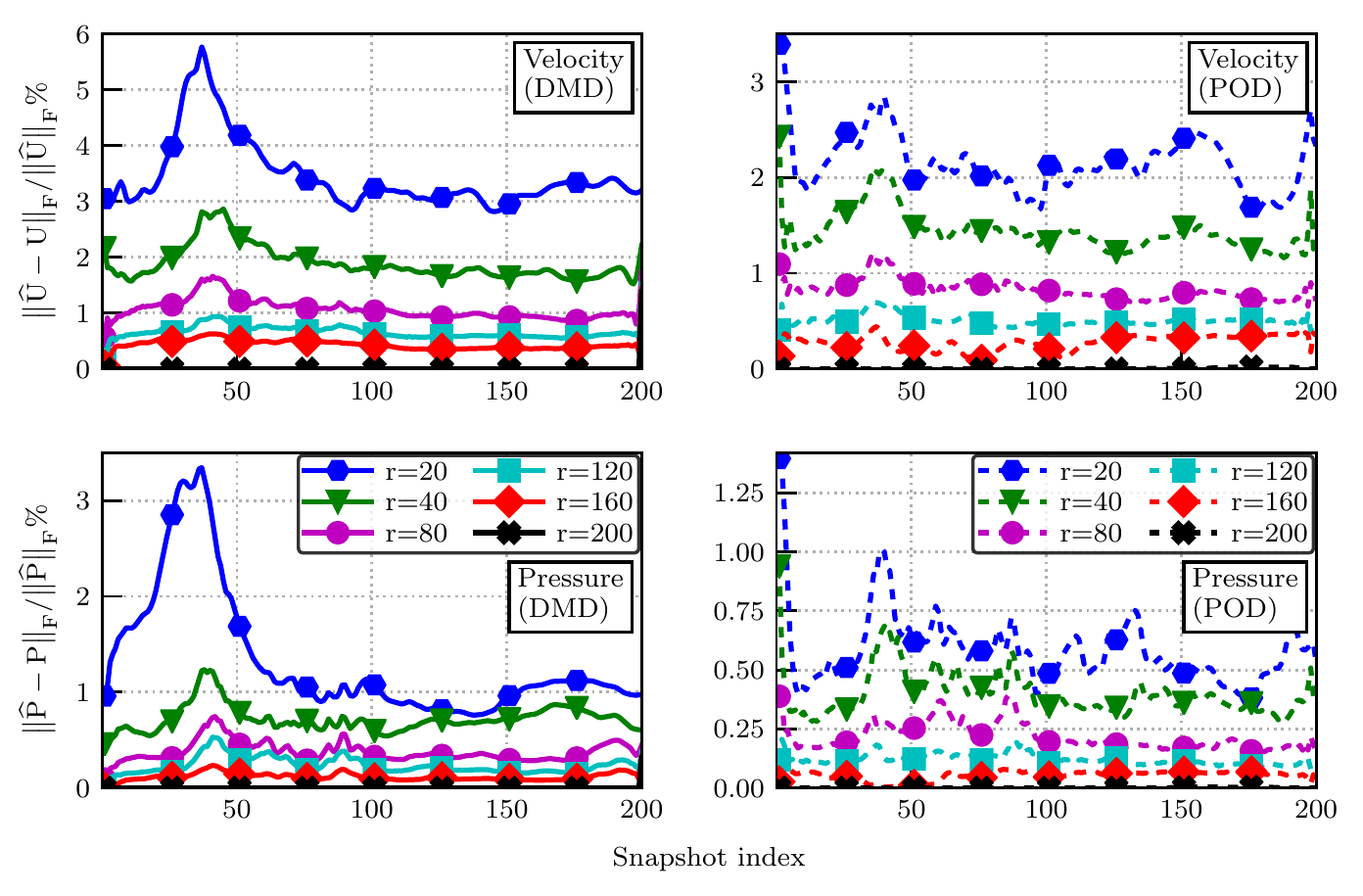}
\caption{Relative error percentage in the Frobenius norm for the
    velocity and pressure (top and bottom, resp.) reconstructed fields
    using the DMD and POD (left and right, resp.)
    modes. \RB{As expected, reconstruction accuracy converges with higher
    SVD truncation rank (r).}}
\label{fig:reconstrConvergence}
\end{figure*}

\autoref{fig:reconstrConvergence} presents a first evaluation of
the effectiveness of the DMD and POD modal decomposition algorithm in reducing
the number of degrees of freedom of the fluid dynamic problem. The
results refer to a reconstruction exercise in which the LES solutions
at all the time steps have been used for the modal decomposition. The
curves in the plots indicate the relative reconstruction error at each
time step for both velocity (top plots) and pressure (bottom plots) when a growing
number of modes are considered. Such relative  reconstruction error is
computed as the Frobenius norm of the difference between the LES
solution vector and the reconstructed one, divided by the  Frobenius norm
of the LES solution. We here remark that to make the velocity and pressure
fields error values comparable, the gauge atmospheric pressure value in the
simulations has been set to one. A null value would in fact result in
lower LES solution norm, leading in turn to large pressure relative
errors compared to the velocity ones, even in presence of comparable
absolute errors. In \autoref{fig:reconstrConvergence}, the two plots on the left refer to DMD
reconstruction results, while the ones on the right present the POD
reconstruction error. As can be appreciated, for both modal decomposition methods
the errors presented follow the
expected behavior, and reduce as a growing number of modes is used in
the reconstruction, until machine precision error is obtained when all
the $200$ modes available are used to reconstruct the $200$
snapshots. More interestingly, the data indicate that both for DMD and POD,
a number of modes between $80$ and $120$ leads to velocity and pressure reconstruction
errors which fall under $1$\% across all the time interval considered. It is worth
pointing out that, compared to most typical low Reynolds and RANS
flows applications of DMD and POD methodologies, such convergence rate is rather
slow, as higher number of modes are needed to obtain comparable
accuracy. This should not surprise, as LES resolves more turbulent
structures than the aforementioned models, resulting in higher spatial
frequencies which in turn require a higher number of modal shapes to
be accurately reproduced. \RB{The peak we see for smaller rank
  truncations around the $40$-th snapshot for the DMD reconstruction
  is justified by the fact that those snapshots present a wider range
  of frequencies thus we need more DMD modes to properly reconstruct
  them.} Finally, the plots suggest that the reconstruction with
POD modes leads to errors that are slightly lower to the corresponding DMD
errors. In fact, the results consistently show that for both the pressure and the velocity
fields, the POD reconstruction error is approximately half of the DMD error obtained
with the same amount of modes.

\subsubsection{Global error}
\begin{figure}
\centering
\includegraphics[width=0.6\textwidth]{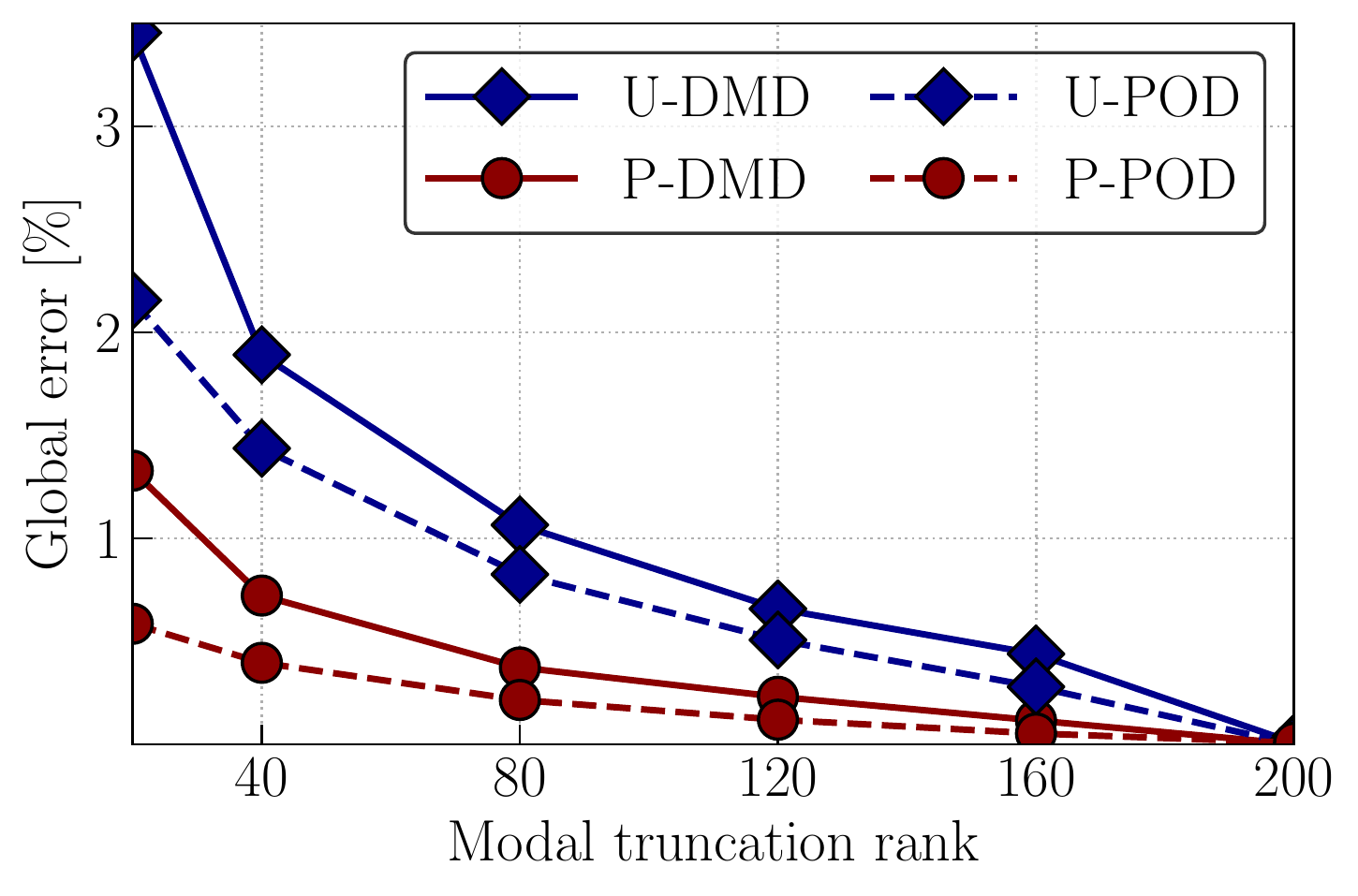}
\caption{Global error of the reconstructed velocity and pressure fields versus SVD truncation rank.
The error denotes spatio-temporal averaging of the flow field data, first by evaluating relative error
in the Frobenius norm, then by averaging over all the snapshots.}
\label{fig:modal_accuracy_recon}
\end{figure}

\autoref{fig:modal_accuracy_recon} presents further verification
of the modal reconstruction accuracy. The curves in the plot represent
the percentage modal reconstruction error --- computed in Frobenius norm
and averaged among all snapshots --- as a function of the number of
modes considered in the reconstruction. The black and red curves confirm
that the reconstructed velocity and pressure fields, respectively, converge
to the corresponding LES fields as the number of modes is gradually increased.
Moreover, these results indicate that an efficient reconstruction,
characterized for instance by a $1$\% relative error, would require more
than $80$ modes for the velocity fields, and more than $30$ for the pressure fields.
These values are higher than those typically observed for RANS and
low Reynolds simulations, probably due to the higher spatial
frequencies typically found in eddy-resolving solution fields.

\subsubsection{Fields visualization}
The error indicators considered in the previous sections are extremely useful
in confirming that the reconstructed solution is globally converging to the LES one
as the number of modes is increased. Yet, they offer little information on
the error distribution in the flow field, and the impact of the reconstruction on
local flow characteristics of possible interest. In particular, for the test case
considered in the present work, it is quite important to assess whether the
reconstruction error does not alter the flow in proximity and in the wake of
the sphere, as such regions are crucial both to the evaluation of the fluid dynamic
forces on the sphere and to the acoustic analysis. To this end, in the present section
we present a series of visualization of the reconstructed flow fields, which
are compared to their LES counterparts.

\begin{figure}
\centering
\includegraphics[width=0.45\textwidth]{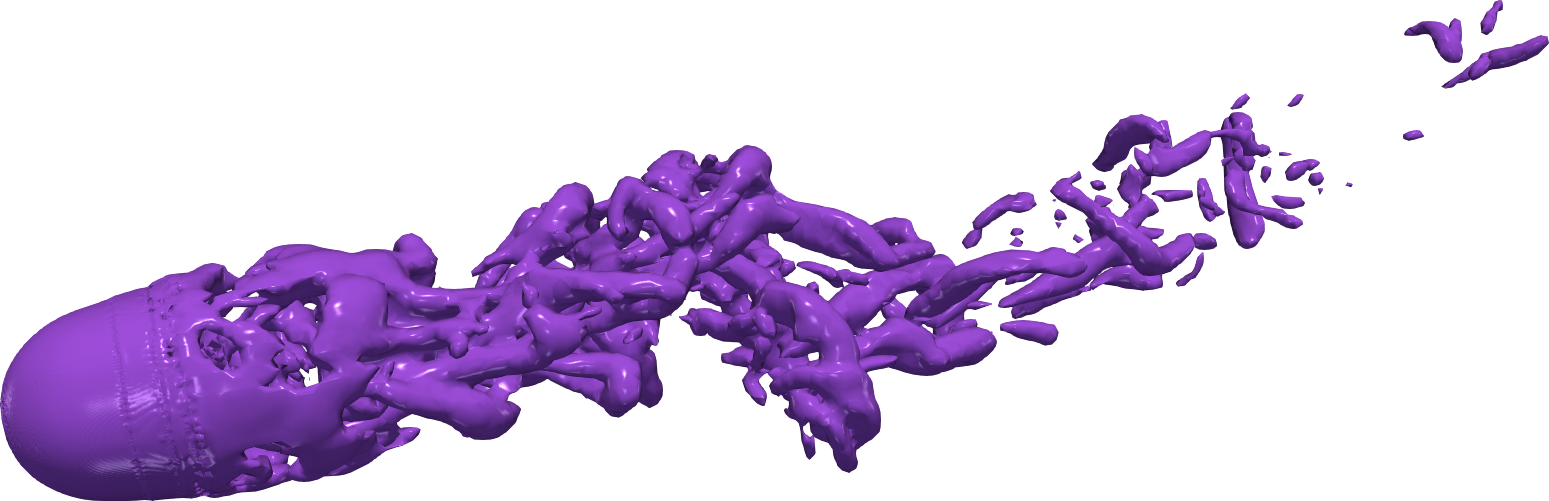} \hfill
\includegraphics[width=0.45\textwidth]{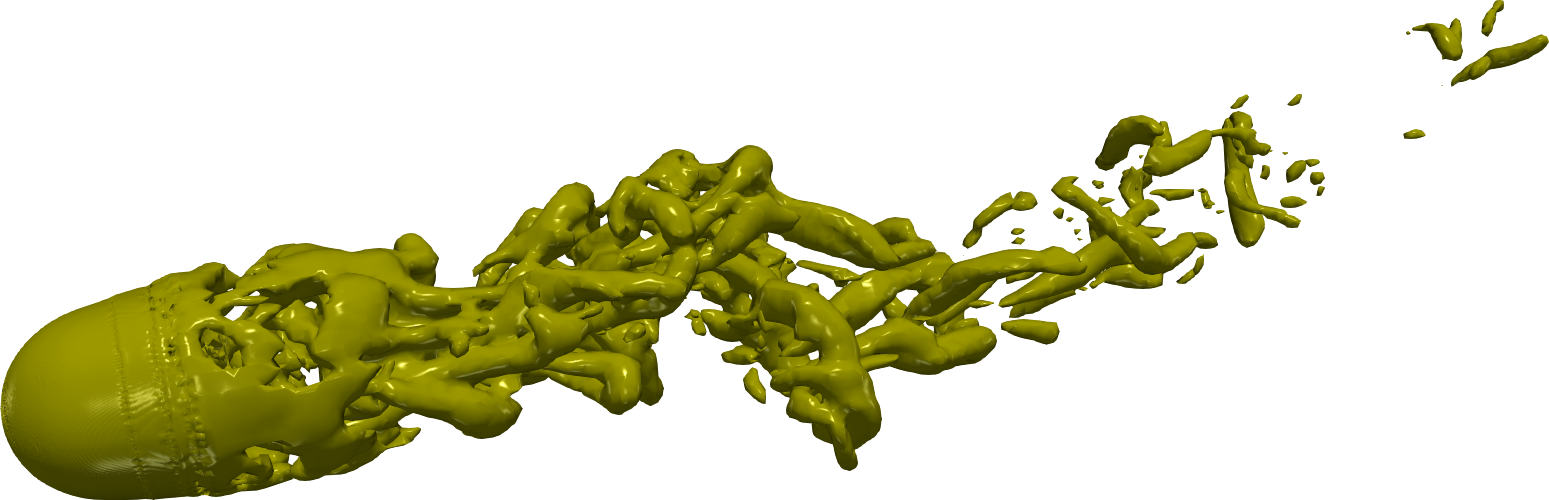} \\
\includegraphics[width=0.45\textwidth]{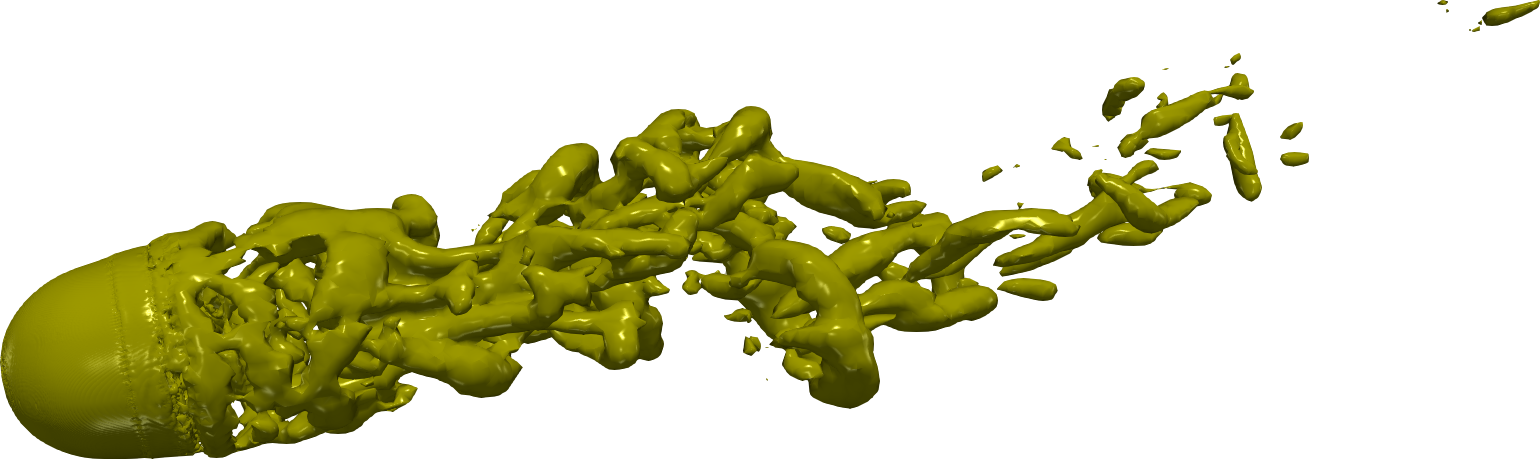}
\caption{Coherent structures represented by the iso-contours of the Q-criterion
         at non-dimensional value $Q D^2 / U_0^2 = 4$ for
         the last snapshot (top), compared with corresponding
         ones obtained from modal reconstruction with DMD (centered) and POD (bottom) using $160$ modes.}
\label{fig:cfd_reconstruct}
\end{figure}

Making use of the Q-criterion, defined as $Q=0.5 (\|\mathbf{\Omega}\|^2-\|\mathbf{S}\|^2)$ with $\mathbf{\Omega}$
and $\mathbf{S}$ denoting the vorticity and strain rate tensors respectively, \autoref{fig:cfd_reconstruct} depicts
the turbulent structures characterizing the flow in the wake region past the sphere. By definition,
a positive value of $Q$ implies relative dominance of the vorticity magnitude over the strain rate~\cite{HALLER2005}.
Here, a positive value of $10^4$ \si{1/s^2} is chosen to generate isosurfaces which pass through all data points holding
this value. In the figure, the top plot refers to the original LES solution obtained at the last snapshot of the dataset, while the centered and
bottom plots refer to the corresponding DMD and POD reconstructions, respectively, utilizing $160$ modes.
The images show that both POD and DMD reconstruction algorithms lead to fairly accurate representation of
the turbulent structures shape past the sphere. In fact, the configuration of the wider vortical
structures detaching from the sphere appears to be correctly reproduced in the reconstructed solution. As for
finer details associated with the smaller turbulent scales, the DMD reconstruction is observed to be in closer
agreement with the original LES solution than the POD reconstructed field. Such observation is consistent with findings
from~\cite{zhang2014identification,Bistrian2015} in which they demonstrated the superiority of DMD to accurately determine spectral
and convective information of the vortical structures in wake regions.

\begin{figure}
\centering
\includegraphics[width=1.\textwidth]{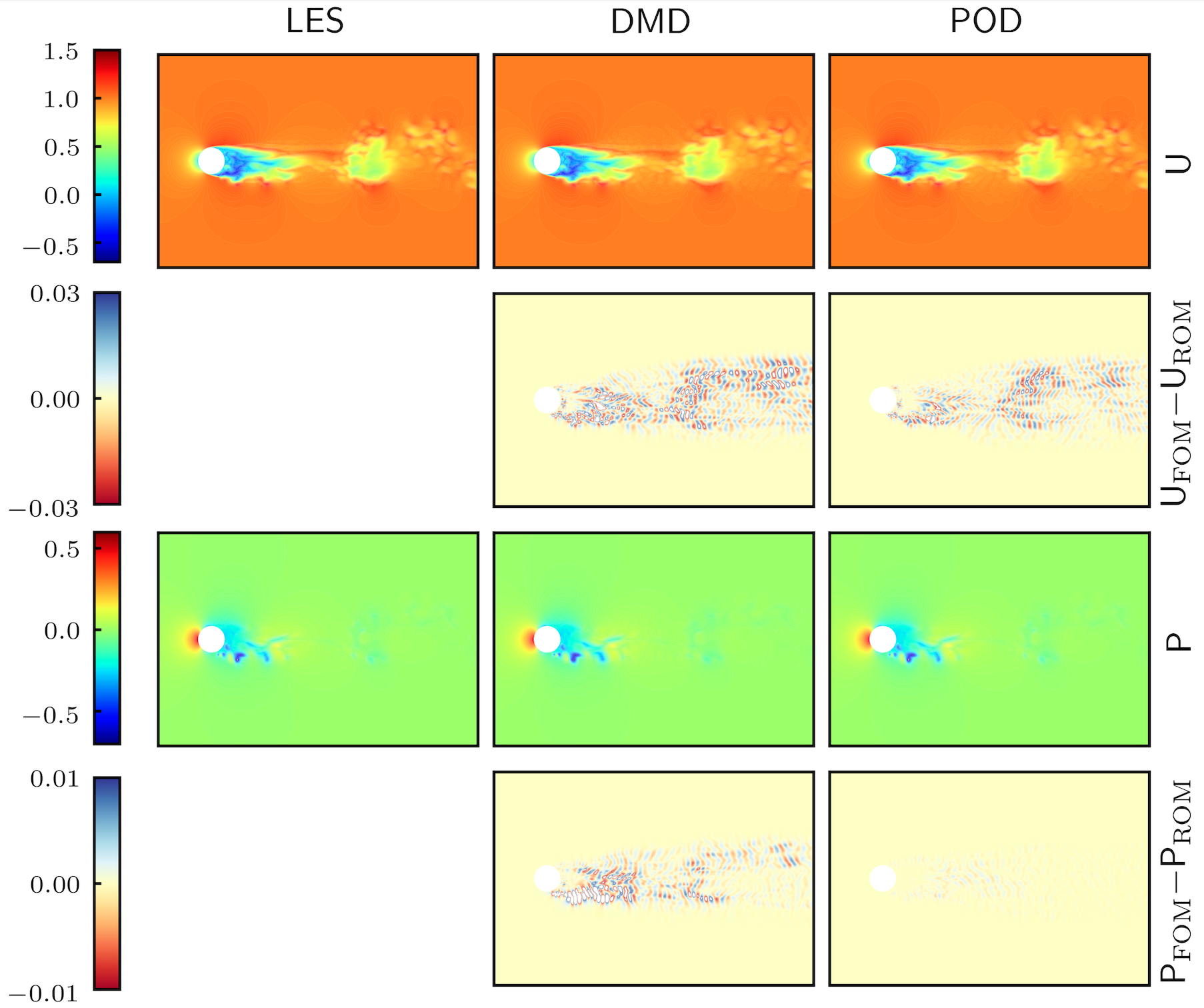}
\caption{\RB{Left to right: LES, DMD and
 POD reconstructions with $N=160$ modes. Top to bottom: streamwise
  velocity, corresponding error fields, pressure, corresponding error
  fields. Instantaneous snapshot is selected based on errors peak.}}
\label{fig:recon_vis}
\end{figure}

For a more significant and quantitative assessment, \autoref{fig:recon_vis} includes a series of
of contour plots representing the instantaneous flow field at time instance corresponding to a maximized relative error, \RB{cf.} \autoref{fig:reconstrConvergence}. 
The three plots in the first row represent contours of the LES streamwise velocity component
field, and of its DMD and POD reconstructed counterparts, respectively, utilizing $160$ modes. At a first
glance, the reconstructed fields seem to reproduce
the main features of the LES flow. In
particular, both the stagnation region ahead of the sphere and the flow detachment
past it appear to be correctly reproduced by both modal reconstruction strategies. In addition,
the detached vortex, located downstream with respect to the
sphere in this particular time instant, is also correctly reproduced. For a better assessment, the two images of the second row
represent contours of the local error for the DMD and POD streamwise velocity reconstruction, respectively.
Both for the DMD and POD reconstruction, higher local error values are found in the wake region
downstream with respect to the sphere. In particular, it observed that the local error peaks
in the DMD reconstruction is larger than that for POD. Additionally, the high frequency error pattern and the
elevated local error values located in the wake region seem to indicate that, as expected, the modes
disregarded in the reconstruction are associated with high spatial frequencies.
A similar comparison is presented for the pressure field in the following rows of the figure. The
three plots in the third row represent the instantaneous pressure field obtained with LES, 
and its reconstructions computed with DMD and POD, respectively, for the same snapshot. Again, the full order field appears well
reproduced by both DMD and POD reconstruction algorithms, as features like the peak pressure in the stagnation
region and the pressure minimum within the vortex detaching past the sphere are correctly reproduced.
The plots of pressure reconstruction error for DMD and POD, respectively, presented in the last row confirm
that both methods are capable to adequately represent the LES solution. Also here, the disregarded
high spatial frequency modes are likely responsible for the high frequency error pattern observed.

\subsubsection{Error statistics in the wake region}
\begin{figure*}
\centering
\includegraphics[width=0.9\textwidth]{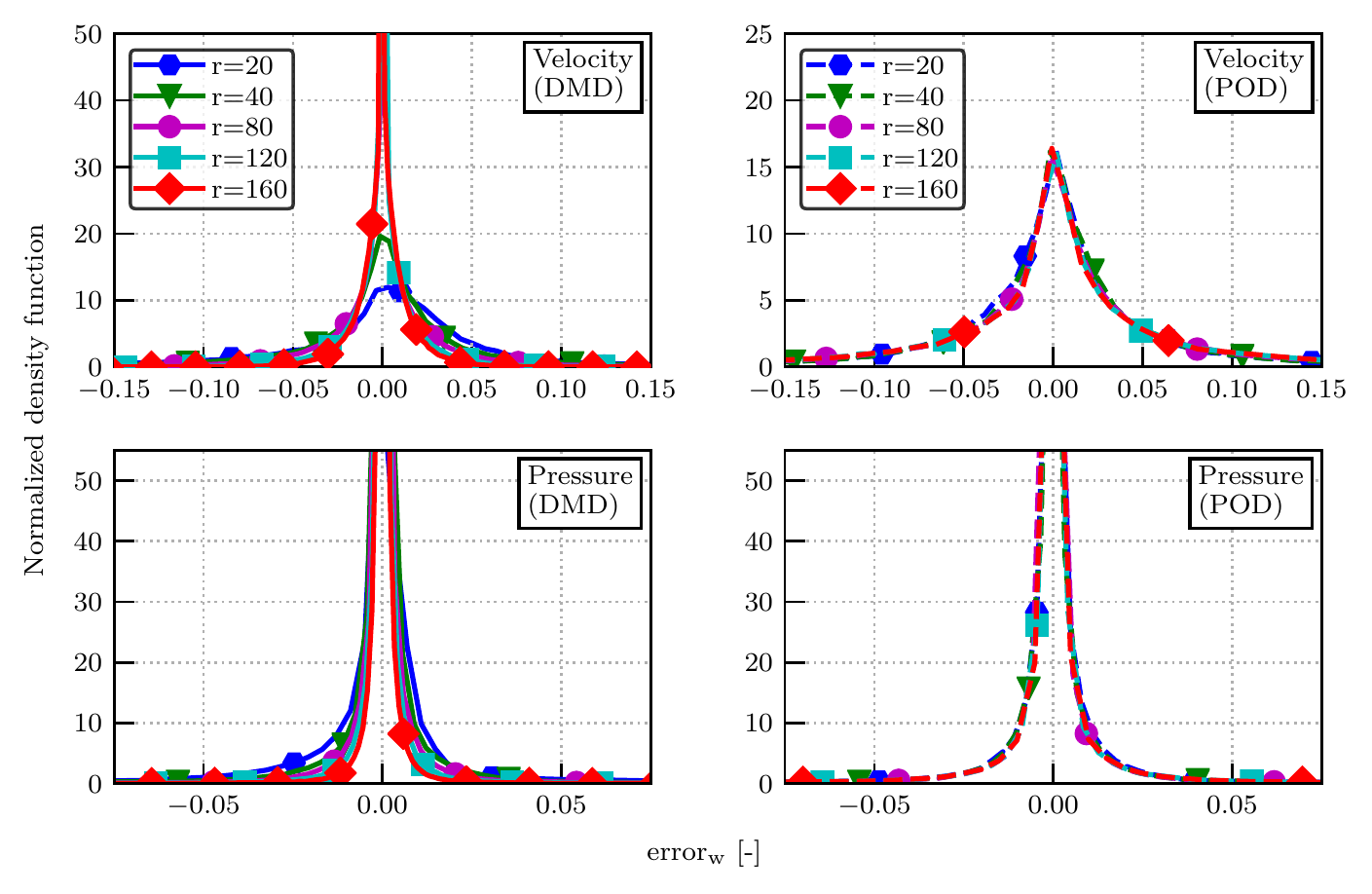}
\caption{Error statistics of velocity and pressure field data of
  selected snapshot corresponding to maximized relative error. Sampled data points are
  conditioned by $\| \bm \omega_x \| = \| \nabla_x \times \bm U \|
  > 1.0$ to identify the wake region. Resulting density function is weighted
  by normalized cell volumes. \RB{Presented plots correspond to various modal truncation ranks (r).}}
\label{fig:stats_recon}
\end{figure*}
As \autoref{fig:recon_vis} shows, the highest reconstruction errors in
the velocity and pressure fields reconstruction are mostly
located in the wake region. Therefore, a follow up analysis is
required to quantify the spatial error distribution within such
region and to assess the local convergence behavior of both the DMD and
POD modal decomposition methodologies. \autoref{fig:stats_recon} presents a normalized density function plot
based on the reconstruction error of both the velocity and pressure
fields corresponding to the snapshot illustrated in
\autoref{fig:recon_vis}.
\RB{To generate the plot, the following steps are performed. First, the computational cells in the wake region are
identified according to the condition $\| \bm \omega_x \| = \| \nabla_x
\times \bm U \| > 1.0$ with $\bm \omega_x$ denoting the streamwise vorticity component.
Second, the arrays corresponding to the velocity and pressure error fields in the wake region
are computed, i.e. $e_w^U=(\hat{U}_w-U_w)$ and $e_w^P=(\hat{P}_w-P_w)$ wherein the subscript denotes the wake cells.
Third, to account for variations in the grid resolution, the mentioned error fields within the wake region are weighted by the corresponding normalized cell volumes $V_c^\textrm{norm}$, i.e.
$\bar{e}_w = (\sum_{i=1}^n V_{c,i}^\textrm{norm} e_{w,i} / \sum_{i=1}^n V_{c,i}^\textrm{norm})$ wherein $n$ is the cell count in the wake.
Finally, the error interval for both the velocity and pressure data are uniformly divided into
equal-width bins and the density function (normalized histogram) is plotted for each bin.}

The scattered plots highlight differences in the
behavior of the two modal reconstruction strategies considered. The
diagrams on the left, which refer to the DMD results for velocity
(top) and pressure (bottom), show in fact a clear error reduction as
the modal truncation order is increased. The curves corresponding to
growing truncation orders tend to get closer to the vertical
axis as the error, displayed on the horizontal axis, is progressively
reduced. The same convergence rate cannot be observed in the POD plots
on the right, as both the velocity (top) and pressure (bottom)
reconstruction error statistical distribution curves appear less
affected by an increase of the modal truncation order. Again, this
observation can be explained in the light of the DMD theory and the
ability of its modes to be organized according to the field spatial
frequencies, unlike POD ones which can instead become contaminated by
uncorrelated structures, according to the claim reported in~\cite{zhang2014identification}.

It is worth pointing out that for image definition purposes, the
left and right tails of \autoref{fig:stats_recon} have not been
reported. Yet, a cross comparison with \autoref{fig:recon_vis} readily
suggests that POD results showed lower error margin --- hence narrower
tails --- in this regard. Therefore, we could infer that higher DMD modes
are capable to capture more frequencies in the wake, resulting in
lower mean error but higher peak errors compared to POD modes.

\subsection{Fields mid cast using DMD and PODI}
\label{sec:results_pred}
The previous sections are focused on assessing the accuracy of the
POD and DMD modal decomposition strategies. The reconstruction results
confirmed that both methods can be considered effective tools for the reduction
of the degrees of freedom of the fluid dynamic problem. We now want to
analyse the ability of the data-driven DMD and PODI reduced order
models considered in this work, in predicting the LES solution at time
steps that are not included in the original snapshot set. In particular,
throughout the remaining analysis of this work, the original $200$ LES snapshots
are decomposed into two sets, one set comprises the $100$ odd snapshots from which they are used
to train the ROMs and are called the \emph{train} dataset hereafter, while the remaining
even snapshots within the temporal interval of the train dataset are contained in a
\emph{test} dataset, and are used to compare the full order solution with the DMD and PODI
model prediction results (i.e. ROM \emph{prediction} dataset).

Before progressing with analysis, a summary on the POD energetic content, as well as the data compression level and computational speedup for the considered DMD and PODI data-driven models at various modal truncation ranks is listed in \autoref{tab:speedup}. Here, the kinetic energy content within a ROM is described by cumulative sum of the POD eigenvalues. The compression level, taken as the arithmetic mean between DMD and PODI in the table, considers the size ratio of the FOM data (train and test datasets, totalling $74$ \si{GB}) to the ROM data (spatial modes of the train dataset, and continuous representation of the temporal dynamics computed via time integration or cubic spline interpolation for DMD or PODI, respectively). Speedup is defined as ($\textrm{CPU}_\textrm{FOM}/\textrm{CPU}_\textrm{ROM}$) with \RB{$\textrm{CPU}_\textrm{FOM}=2.4\times 10^5\si{sec}$} corresponding to the time required to solve for the $200$ snapshots and to write out the train dataset, while $\textrm{CPU}_\textrm{ROM}$ is the time needed to 1) perform modal decomposition on the train dataset, 2) extract the associated dynamics with a continuous representation, and 3) write out the prediction dataset in OpenFOAM format.
\definecolor{lightgray}{rgb}{.8,.8,.8}
\begin{table}
\centering
\caption{ROM performance at various modal truncation rank. Data
  compression level is defined as
  (\si{GB}$_\textrm{FOM}$/\si{GB}$_\textrm{ROM}$) with
  \si{GB}$_\textrm{FOM}$=$74$. Speedup is defined as
  ($\textrm{CPU}_\textrm{FOM}/\textrm{CPU}_\textrm{ROM}$) with
  \RB{$\textrm{CPU}_\textrm{FOM}=2.4\times 10^5$\si{sec}.} Compression level
  is averaged between DMD and PODI.} \label{tab:speedup}
\footnotesize
\begin{tabular}{ c c c c c c }
\hline
\hline
Rank & \multicolumn{2}{c}{POD cumulative energy} & \multicolumn{2}{c}{Speedup} & Compress. \\
(r) & U & P &  DMD & PODI & level \\
\hline
\hline
\rowcolor{lightgray}
10 & 0.999992 & 0.994093 & 29.29 & 35.67 & 19.995\\
20 & 0.999996 & 0.997076 & 28.46 & 32.46 & 9.9973 \\
\rowcolor{lightgray}
40 & 0.999998 & 0.998735 & 26.41 & 27.93 & 4.9986\\
60 & 0.999999 & 0.999295 & 25.41 & 24.68 & 3.3324\\
\rowcolor{lightgray}
80 & 0.999999 & 0.999571 & 23.58 & 22.48 & 2.4993 \\
100 & 1.000000 & 1.000000 & 23.18 & 20.46 & 1.9995 \\
\hline
\hline
\end{tabular}
\end{table}

It is worth noting that, since time is the considered parameter in the present ROM procedure, an offline phase still requires solving for the same temporal window as in the FOM solution, in order to obtain the train dataset. Nevertheless, a computational gain in the generated ROM can be still attained, since the less amount of stored data (i.e. few modes and associated dynamics) compared with FOM are utilized for a swift construction of the fluid dynamic fields at arbitrary sample points. Such procedure is considered a lot faster and more economic than recomputing the simulation to write out the fluid dynamic field data at those sample points. Indeed, an extension of the present procedure to consider additional parameters while performing, for instance, multiparameter interpolation using PODI, would result in a further speedup since the offline phase would solve a full order solution only for a subset of the parameter combinations in the parameter space, hence saving up complete FOM computations from being performed. Such multiparameter investigation is considered a follow up study of this work. To this end, now we progress the analysis by considering the temporal evolution of the error due to ROMs prediction.

\subsubsection{Prediction error analysis}
\label{sec:Snapshot_wise_errors}
\begin{figure*}
\centering
\includegraphics[width=0.9\textwidth]{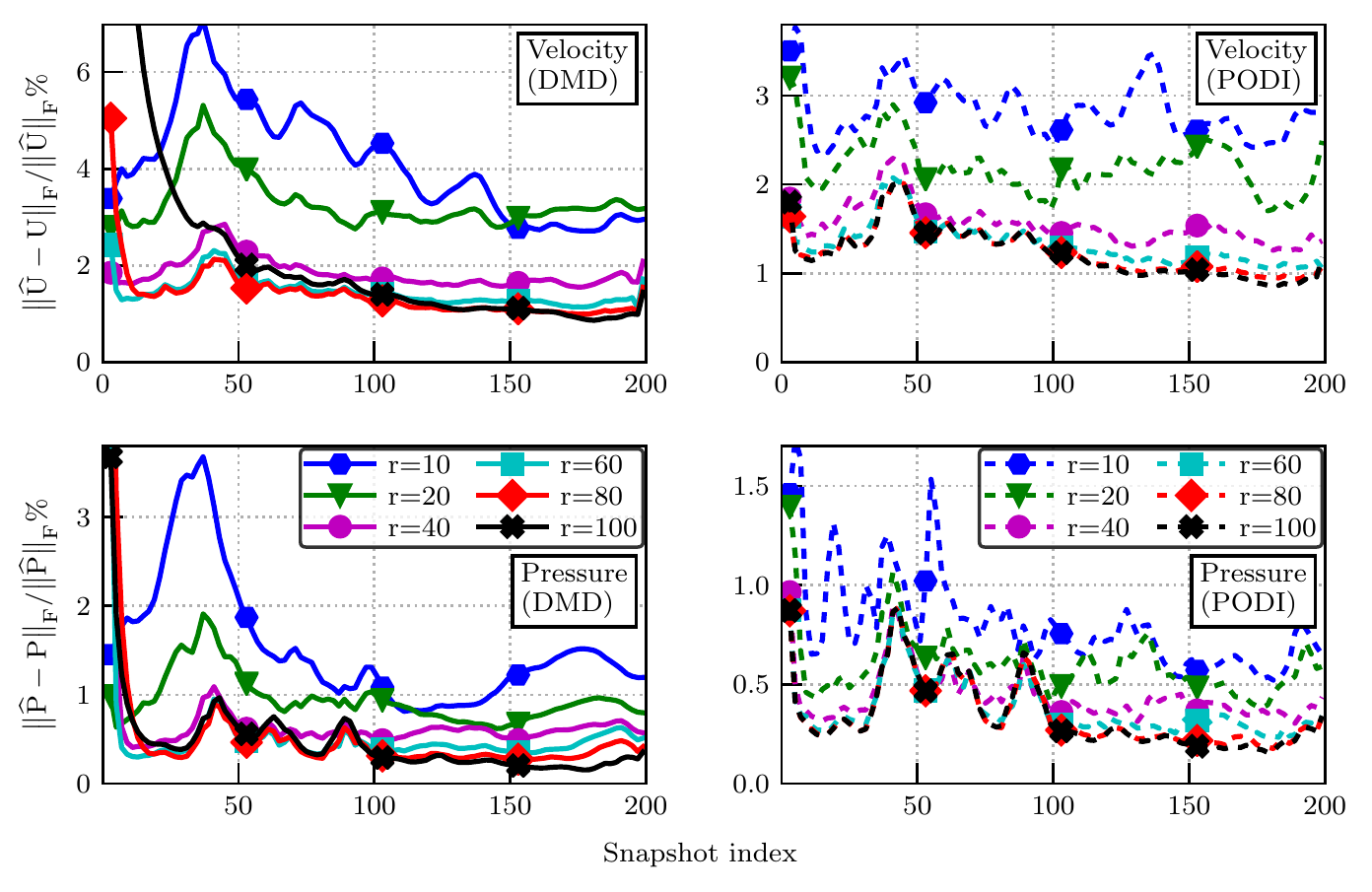}
\caption{Relative error percentage in the Frobenius norm for velocity and
  pressure (top and bottom, resp.) fields, where half sample rate snapshots are used to
  train the reduced model, DMD (left) or PODI (right), \RB{for particular truncation (r)} while predicting the
  intermediate snapshots. Plots show only the
  snapshot-wise prediction error, while disregarding errors in the
  training set, which is almost null.}
\label{fig:error_inter}
\end{figure*}

The first test presented is designed to assess the predictive accuracy of the DMD and PODI methodologies
described in this work. \autoref{fig:error_inter} depicts the results obtained with
the ROM methodologies applied in such alternate snapshots arrangements as previously described. The top plots
report the percentage Frobenius norm error between the DMD and PODI predictions with respect to the LES
solution for the longitudinal component of the velocity, while the bottom diagrams present
similar error norms corresponding to the pressure field. The left plots
refer to the DMD results, while the right ones depict the PODI errors. We point out
that, the --- significantly lower --- error values corresponding to the train dataset
have been omitted for clarity. Finally, the different colors in the plots indicate growing
number of modes considered in the DMD and PODI prediction, up to a maximum of $100$.

The results show a substantial convergence of the data-driven reduced model solutions to the full order one.
In particular, selecting a number of modes between 80 and 90 results in errors lower than $2$\% on the velocity
field for most time steps in both DMD and PODI methods. As for the pressure field, the bottom left diagram in \autoref{fig:error_inter}
shows that $1$\% error goal can be obtained with an even slightly lower amount of DMD modes. It must be pointed out though,
that the DMD solution in the very first time steps is not converging to the LES one, and the error grows
as the number of modes is increased. This behavior, which could be related to the high frequencies
introduced by the higher DMD modes added to the solution, is currently under
further investigation. Aside from the first few time steps, it is
generally observed that PODI errors are consistently
lower by a factor two with respect to the DMD ones, especially in the cases utilizing fewer modes.
It is worth pointing out that a direct comparison between the PODI plots in
\autoref{fig:modal_accuracy_recon} and those in \autoref{fig:error_inter}
can indirectly result in a possible estimate of the interpolation error associated with the
PODI strategy. The plots suggest that the effect of interpolation on the global spatial error at each
time step is rather low, as the errors in \autoref{fig:error_inter} present the same behavior and
are not significantly higher than the ones obtained with pure reconstruction.


\begin{figure}
\centering
\includegraphics[width=0.6\textwidth]{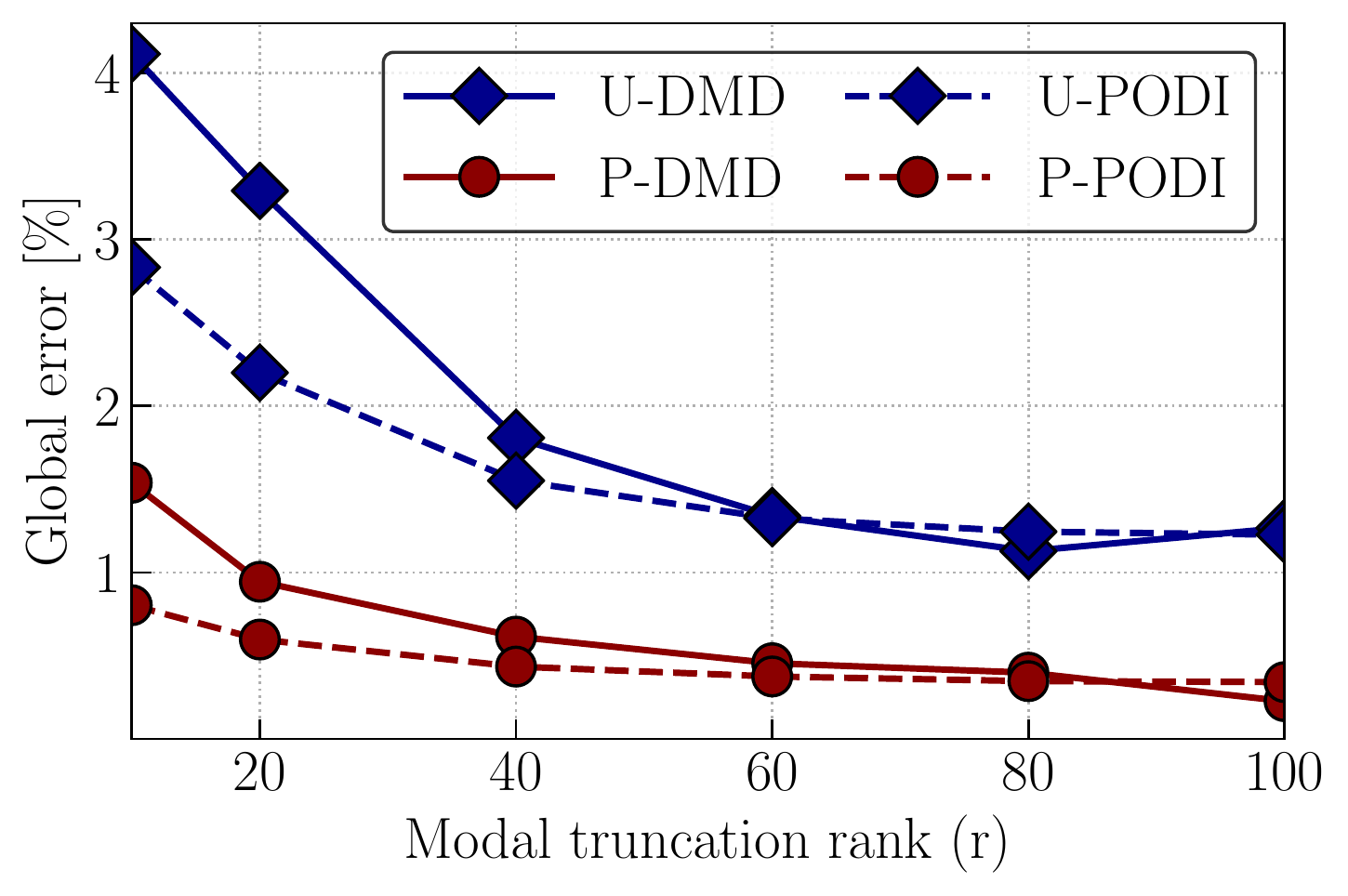}
\caption{Global error of the prediction dataset versus SVD truncation rank.
The error denotes spatio-temporal averaging of the flow field data, first by evaluating relative error
in the Frobenius norm, then by averaging over all the snapshots.}
\label{fig:error_global_inter}
\end{figure}
Now, to summarize the performance of each ROM with respect
to the modal truncation level, a global spatio-temporal error is measured against growing number of modes,
as depicted in \autoref{fig:error_global_inter}. In particular, similar to
the previous analysis in \autoref{fig:modal_accuracy_recon}, the mean absolute
error is evaluated for the percentage Frobenius norm
spatial error for all the prediction dataset, cf. \autoref{fig:error_inter}.
Here, it is observed that, for both the velocity and pressure predictions,
a significant reduction in the global error is achieved in the PODI models with respect
to DMD up to a utilization of $60$ modes, after which the global error becomes
comparable between both ROMs. Additionally, a finite
global error is noted in the figure even when a full modal rank is utilized.
Indeed, such observation should not be surprising since half the dataset is only employed to train both ROMs
while interrogating the remaining sample points. This is not
the case in pure field reconstructions, cf. \autoref{fig:modal_accuracy_recon},
where the global error vanishes with a full modal rank since dynamic and spectral
information become entirely recovered.

\subsubsection{Coherent structures}
The prediction error indicators considered in the previous
section indicate whether the reduced models solution is globally
converging to the LES one as the number of modes is increased. We now
resort to flow visualizations to obtain better information on the
error distribution in the flow field and on the reduced models
performance in reproducing local flow characteristics of possible
interest.

\begin{figure}
\centering
\includegraphics[width=0.45\textwidth]{fig/Q_LES.png} \hfill
\includegraphics[width=0.45\textwidth]{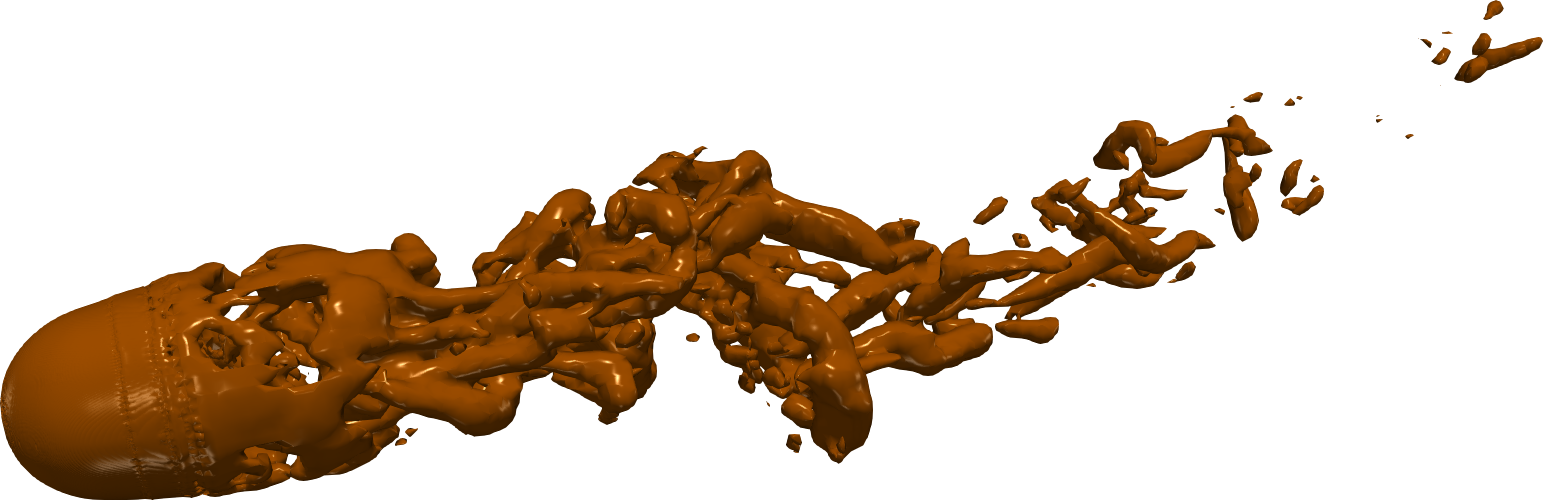} \\
\includegraphics[width=0.45\textwidth]{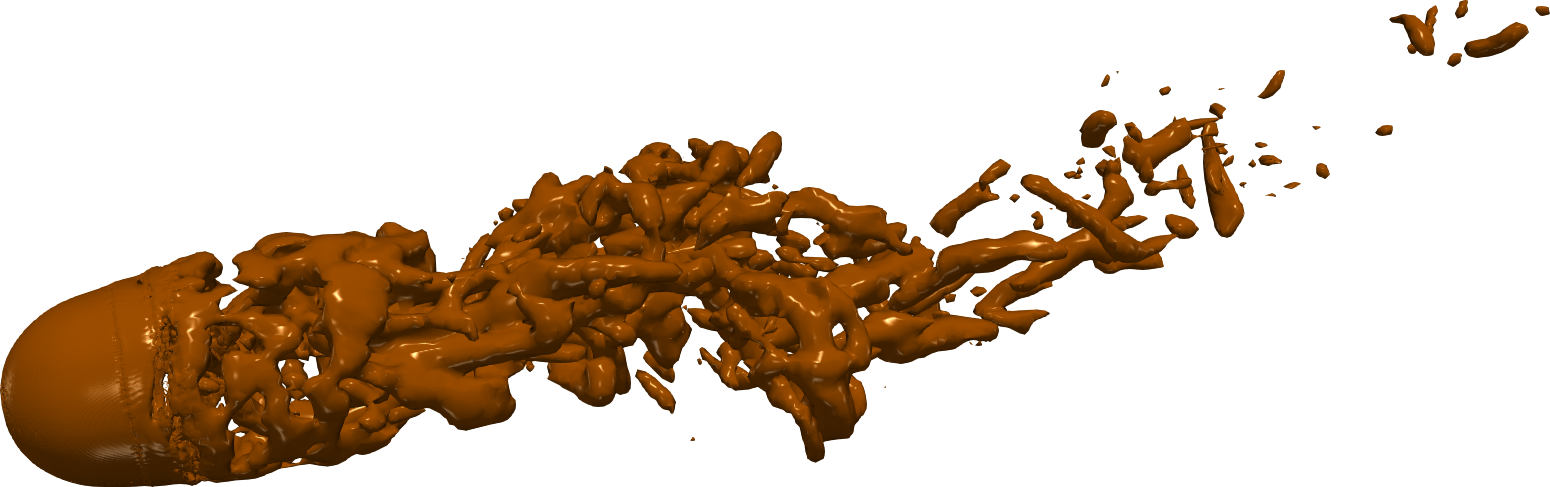}
\caption{Coherent structures of the last snapshot in the prediction dataset, represented by the iso-contours of the
  Q-criterion at non-dimensional value of $Q D^2 / U_0^2 = 4$,
  compared with corresponding ones obtained from DMD-$100$ (centered),
  and PODI-$100$ (bottom).}
\label{fig:cfd}
\end{figure}

Again, making use of the Q-criterion isosurfaces at value $Q D^2 / U_0^2 = 4$,
\autoref{fig:cfd} portrays the turbulent structures
characterizing the wake flow past the sphere. The top plot
refers to the original LES solution obtained at the last
snapshot of the prediction dataset, while the centered and bottom plots refer to the
respective DMD and PODI predicted solution utilizing $100$ modes.
Also in this case, the images show that both DMD and PODI reduced
order models allow for rather accurate reproduction of the turbulent
structures past the sphere. For both ROM solutions, the main vortical structures detaching
from the sphere appear in fact very similar to those of the original
LES flow field. Finer details associated with smaller
turbulent scales are also in good agreement, hence suggesting
that errors in the PODI time interpolation and DMD time integration
are not significantly higher with respect to the reconstruction error analysed earlier.

\RB{
\subsubsection{Drag and lift coefficients}
\label{sec:lift_drag}

Besides coherent structures, it is important to assess the performance of data-driven ROMs in capturing hydrodynamic phenomena that could find particular interest in the engineering community. One of these phenomena is the drag and lift forces on the sphere surface. First, pressure-induced drag and lift forces are evaluated by integrating the pressure in the streamwise and the flow-normal directions, respectively, over the sphere surface. Then, the pressure forces are normalized by the dynamic pressure acting on the sphere, with a sphere projected cross-sectional area being considered for the reference area, to compute coefficients of drag ($C_D$) and lift ($C_L$). The temporal evolution of both coefficients on the sphere surface is depicted in \autoref{fig:coeffs} for various ROM predictions compared with the FOM.
The presented plots show that both $C_D$ and $C_L$ are accurately predicted while only half the dataset is employed.
It is noted that the accuracy of $C_D$ predictions is higher in comparison to $C_L$ predictions, and that employing as low as $40$ modes are able to predict the forces coefficients with a sufficient accuracy. More importantly, the predictive character of PODI models are observed to be superior over DMD models for a particular modal truncation level. The observed discrepancies in the early temporal window of $C_L$ predictions via DMD at $r=80$ are related to the poor DMD predictions of the flow fields at early snapshots, as previously noted in the relative error plots, cf. \autoref{fig:error_inter}. On a general note, the presented data-driven ROMs are able to predict the pressure dynamics on the surface of the immersed sphere while employing considerably lower levels for modal truncation.
We would like to highlight that the framework that we are proposing is
solely based on input and output quantities. Therefore, instead of
computing engineering quantities from the reconstructed flow fields, one could
also construct a reduced order model to directly approximate them. In such a
case we would expect that the so generated ROM would produce even more accurate
results. 
\begin{figure}
\centering
\includegraphics[width=0.9\textwidth]{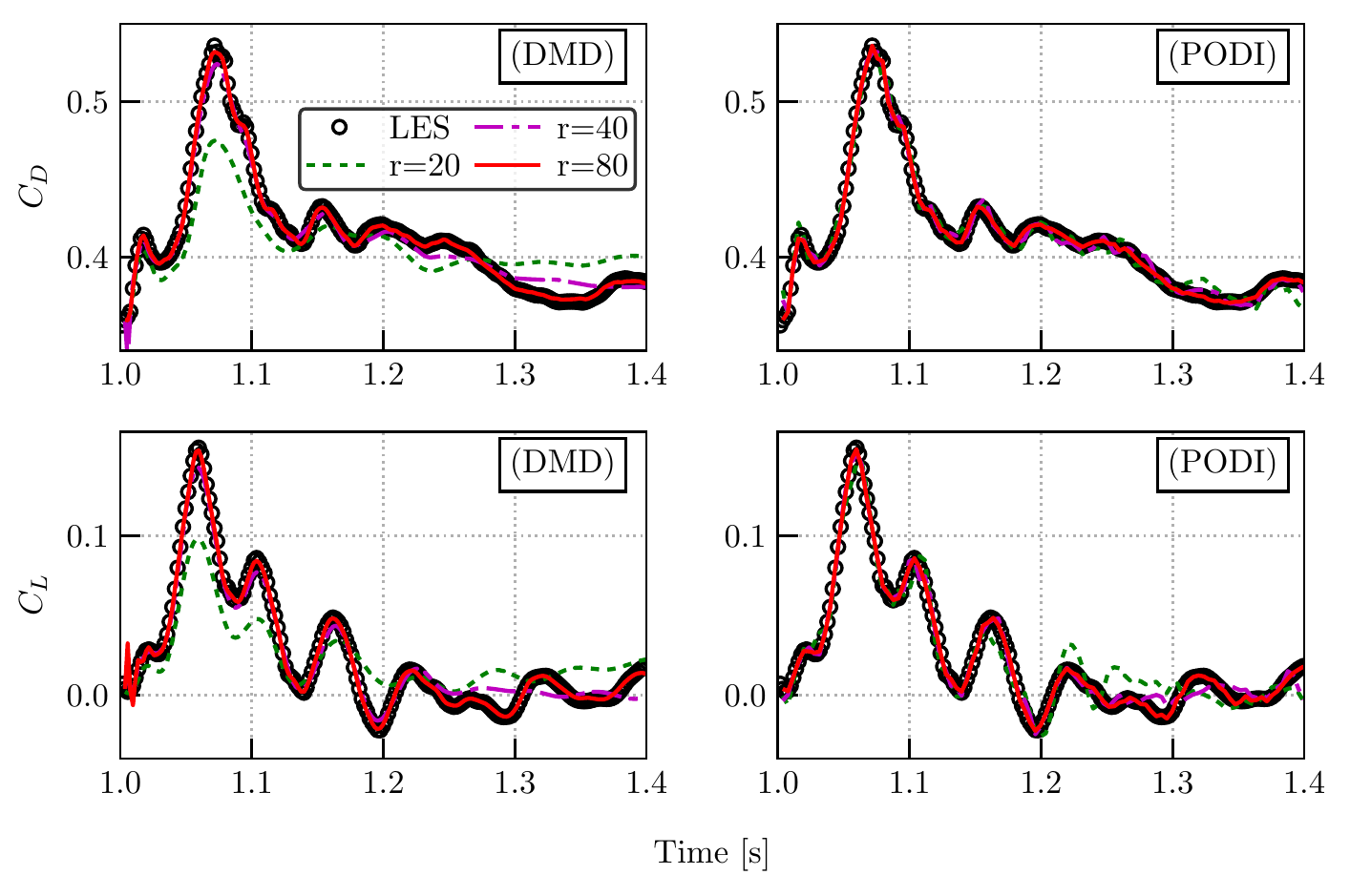}
\caption{Temporal evolution of the drag and lift coefficients corresponding to pressure forces on the immersed sphere surface. The force coefficients are computed for various DMD and PODI models in comparison to FOM data.
Presented plots show the predictive performance of both data-driven ROMs on estimating the sphere pressure forces, with noted superiority of PODI predictions over DMD.}
\label{fig:coeffs}
\end{figure}
}

\subsection{Spectral analysis}
\label{sec:results_spectral}
\label{sec:results_ac}

\subsubsection{Data probes inside and outside the wake region}
\begin{figure}
\centering
\includegraphics[width=0.9\textwidth]{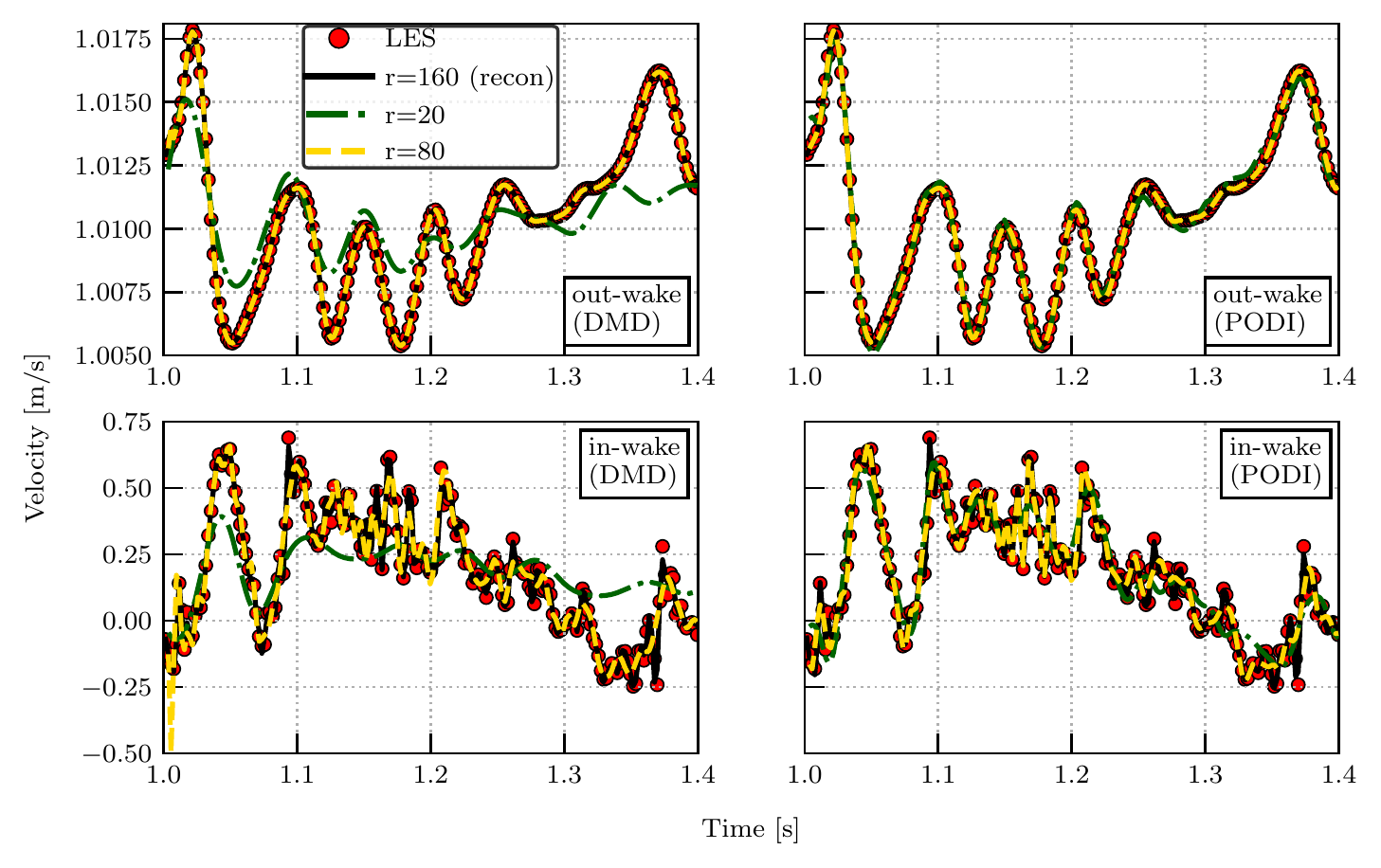}
\caption{Time history of streamwise velocity signals at fixed
  points located outside the wake at $\mathbf{x}_A=\{2D, 0,
  2D\}$ (top) and inside the wake at $\mathbf{x}_B=\{2D, 0,
  0\}$ (bottom) for DMD (left) and PODI (right) models. Red symbols
  represent LES results, while dashed lines represent the ROM results for
  different SVD truncation ranks \RB{(r)}. Black continuous line represents the reconstructed
  field utilizing $160$ modes.}
\label{fig:U_signals}
\end{figure}

A first necessary step in order to evaluate how the frequency content
of the FOM solution is reproduced at the ROM level, consists in
observing the time dependency of the solution at fixed points in the
computational domain. To this end, \autoref{fig:U_signals} presents
the reduced models results for the time evolution of the streamwise
velocity component in correspondence with two different
locations of the flow field. In the $2\times2$ tabular arrangement of
the figure, the top plots correspond to a point located outside of the sphere
turbulent wake,  while the bottom plots refer to a point inside the turbulent
wake. In addition, the two diagrams on the left refer to the DMD results,
and the ones on the right report the PODI result. In the plots, the red filled
circles refer to the LES result at prediction sample points, while the dash-dotted lines represent
the ROM results, which have been drawn for a growing number of
modes. Finally, the black continuous line represents the pure DMD and PODI
reconstruction result utilizing $160$ modes.

As expected, the plots show that the velocity variation over time becomes more
chaotic in the wake region, due to the high frequency fluctuations
associated with the turbulent structures typically observed in such
part of the flow domain. All the plots show that, as the number of DMD
and PODI modes used is increased, the reduced solution approaches the
full order model one. At a first glance, the PODI convergence at
the point located outside of the wake appears faster than that obtained with DMD.
The reduced PODI solution obtained with $20$ modes is, in fact, already rather
close to the original time signal, while the corresponding DMD solution
is still far from the LES one. The behavior of both DMD and PODI reduced models
is clearly more accurate when a higher number of modes is selected. The
$80$ modes curve obtained with both ROMs is practically indistinguishable
from the LES one. Given the small estimate of the time interpolation error, as
discussed in \autoref{sec:Snapshot_wise_errors}, the faster PODI convergence
observed should depend on the fact that single PODI modes are richer in spatial
frequency content than the DMD ones. The low order PODI modes might then
already include higher frequencies not contained in the corresponding DMD modes.
Such spatial frequency content should finally reverberate in the time evolution
of the solution including some higher frequencies also when the only lower order PODI
modes are used. Despite these favorable characteristics, when the full order
model solution presents even higher frequencies, also the PODI solution requires a
higher number of modes to obtain satisfactory accuracy.
This is clear by the plots of
the time history of the velocity at the point within
the wake region. Here, both PODI and DMD  reduced result curves become
sufficiently close to the original
solution one only when $80$ modes are considered. This should not come
as a surprise, as the in-wake signal presents higher frequency
fluctuations due to the wake turbulent structures, which can be fully
represented only by including high order modes.

\begin{figure}
\centering
\includegraphics[width=0.9\textwidth]{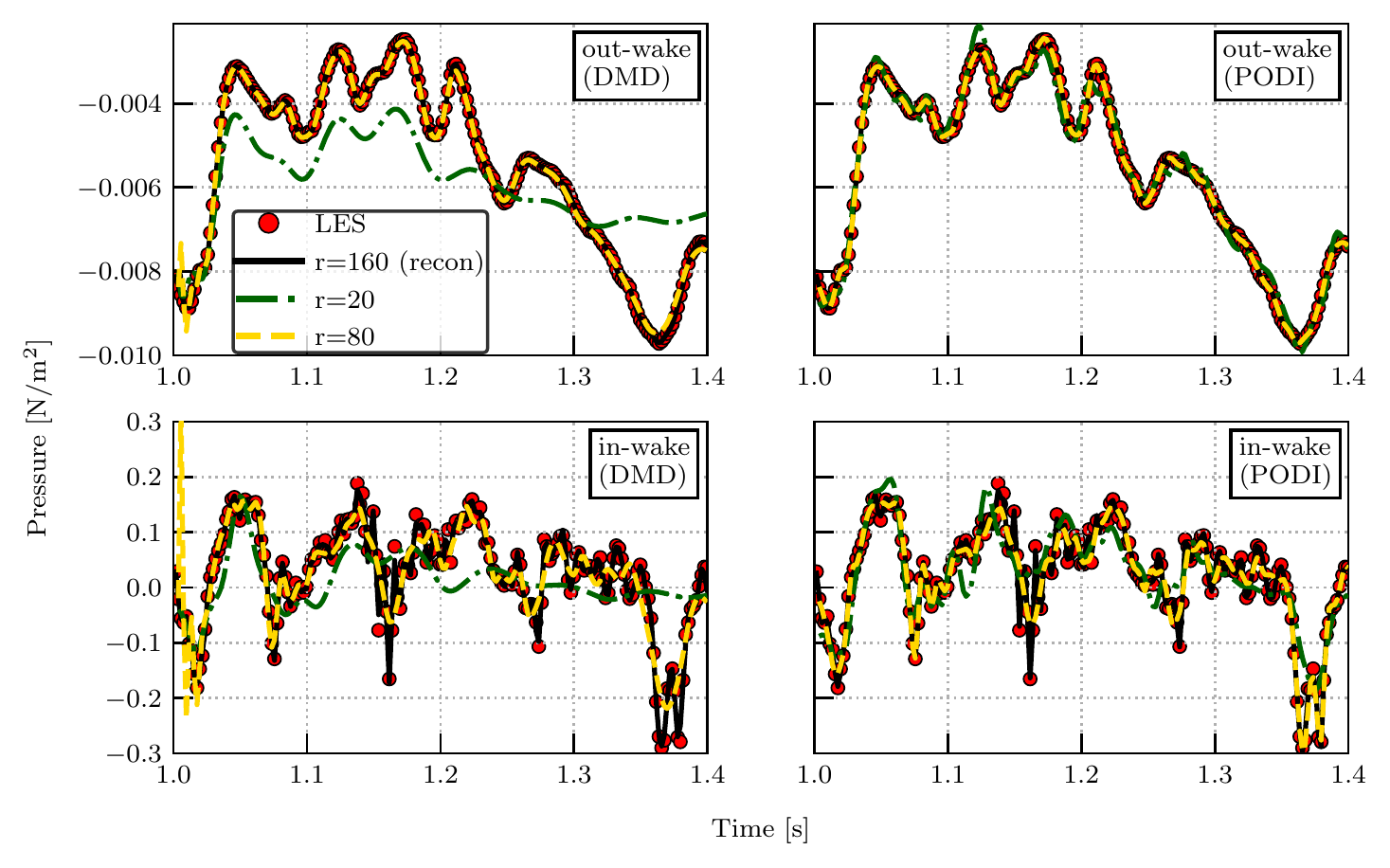}
\caption{Time history of pressure signals at fixed
  points located outside the wake at $\mathbf{x}_A=\{2D, 0,
  2D\}$ (top) and inside the wake at $\mathbf{x}_B=\{2D, 0,
  0\}$ (bottom) for DMD (left) and PODI (right) models. Red symbols
  represent LES results, while dashed lines represent the ROM results for
  different SVD truncation ranks \RB{(r)}. Black continuous line represents the reconstructed
  field utilizing $160$ modes.}
\label{fig:P_signals}
\end{figure}

The reduced pressure local time evolution prediction results presented
in \autoref{fig:P_signals} exhibit the same behavior observed for the
reduced velocity field. Also this Figure, in which the plots and the curve
colors are arranged as described in \autoref{fig:U_signals}, suggests in fact
that both inside and outside the wake region, employing $80$ modes allows
for both ROMs to obtain pressure predictions that are sufficiently
close to the corresponding LES time signals. Again, the
plots also suggest that PODI in the region outside the wake is able
to obtain viable pressure predictions with fewer modes. Thus, these results indicate that
in presence of fully attached flows, in which wake effects are less dominant,
the choice of POD as modal decomposition methodology could result in more
economic reduced models.

\subsubsection{Fast Fourier transform}

\begin{figure*}
\centering
\includegraphics[width=0.9\textwidth]{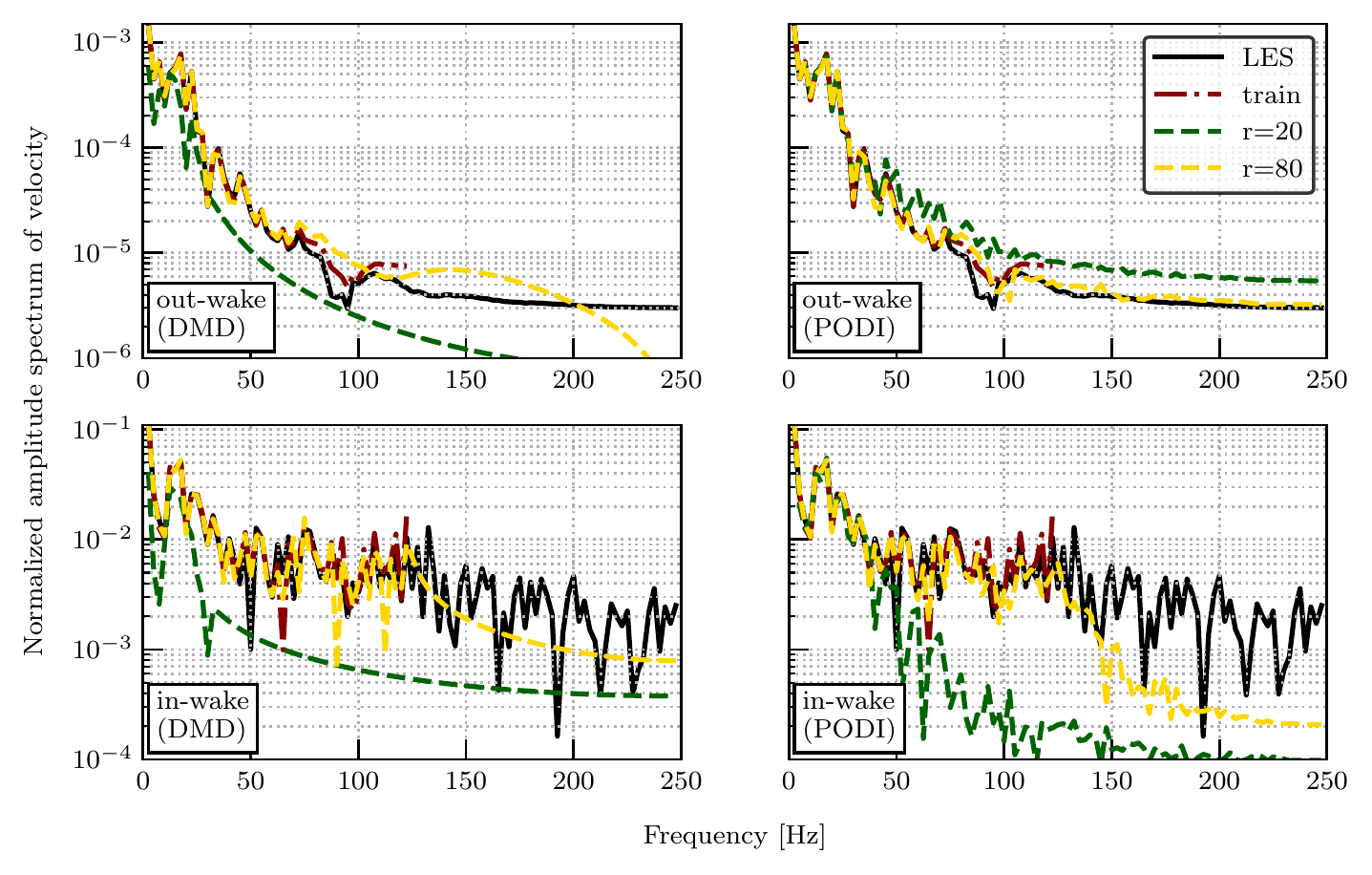}
\caption{Normalized amplitude spectrum of the streamwise velocity
  component, computed as ($|
  \text{fft}(U-\overline{U})|/N_\textrm{freq}$) with
  $\overline{U}$ denoting mean velocity, in a fixed point located
  outside the wake (top) and inside the wake (bottom) for DMD (left)
  and PODI (right) \RB{at different truncation ranks (r)}. Black continuous and red dash-dotted lines refer to LES
  results of the full dataset and train dataset, respectively.
  Green and yellow dashed lines mark ROM surrogates employing $20$ and $80$ modes, respectively.}
\label{fig:velocity_spectra}
\end{figure*}
\begin{figure*}
\centering
\includegraphics[width=0.9\textwidth]{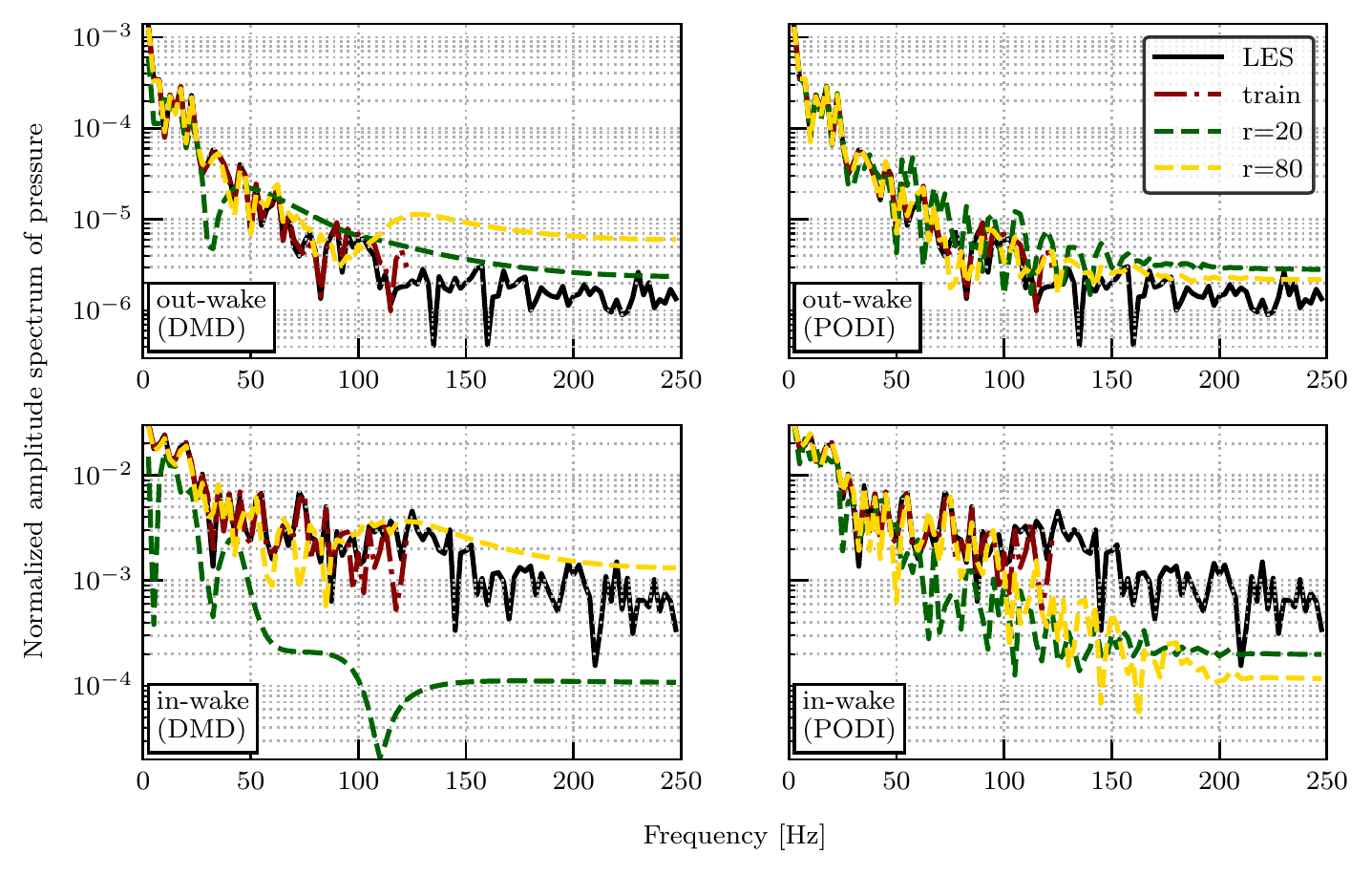}
\caption{Normalized amplitude spectrum of the pressure signal, computed as ($|
  \text{fft}(p-\overline{p})|/N_\textrm{freq}$) with
  $\overline{p}$ denoting the mean pressure, in a fixed point located
  outside the wake (top) and inside the wake (bottom) for DMD (left)
  and PODI (right) \RB{at different truncation ranks (r)}. Color arrangements are same as in \autoref{fig:velocity_spectra}.}
\label{fig:pressure_spectra}
\end{figure*}

The plots in \autoref{fig:U_signals} and \autoref{fig:P_signals}
seem to indicate that, on a qualitative level, the PODI and DMD
solutions are able to recover a growing portion of the full order
model frequency content as the number of modes is gradually
increased. This aspect should be, of course, investigated in a more
accurate way, as the presence of high frequency components in the ROM
solution results in their ability to be effective surrogates in
acoustic analysis. Thus, we make use of the FFT spectra with the aim
of obtaining a quantitative assessment of the impact of the POD and
DMD modes considered on the frequency content of the
solution. \autoref{fig:velocity_spectra} presents the magnitude of the
FFT of the local streamwise velocity signals previously
presented. Also in this case, the top plot refers to the point
located outside of the wake, while the bottom plot refers to the
in-wake point spectrum. The black, continuous lines refer to the FFT
of the local velocity signal obtained with the LES full order
model. The red dash-dotted lines represent the spectrum of the signal
composed only by the train dataset which, of course, is
characterized by half the sampling frequency with respect to the full
LES signal. Finally, the green and yellow dashed lines respectively
denote the corresponding plots obtained by means of truncated DMD (left panel) and
PODI (right panel) models employing $20$ and $80$ modes each. \autoref{fig:pressure_spectra}
--- which also employs the line color arrangement just described ---
displays analogous spectral results obtained when the pressure field
is considered. \RB{We point out that, based on experimental results
presented in \cite{SakamotoVortShedding}, a Re = 5000 flow past a
sphere results in the Strouhal number $S_t=0.2$ associated with
vortex shedding. This corresponds to $f = S_t U / D = 20 \, Hz$
being $U = 1.0 \, m/s$ and $D = 0.01 \, m$.
In \autoref{fig:velocity_spectra} and \autoref{fig:pressure_spectra},
it is clearly visible that a  $20 \, Hz$ peak appears in the LES
signal spectra when the velocity and pressure probe is located within
the cylinder wake (lower plots). Such peak is accurately reproduced by
both PODI and DMD results.}


The top plots in both \Cref{fig:velocity_spectra,fig:pressure_spectra} show that all the reduced models
considered allow for a sufficiently good reconstruction of the
solution spectra for a point outside of the wake. In fact,
all the yellow dashed curves appear very close to the black continuous one
representing the LES solution spectrum. Thus, this confirms that
considering $80$ modes or more leads to a spectrum that is
indistinguishable from the original, especially for the PODI results in which a relative superiority
is again noted in comparison with the spectra from DMD models. Yet, the higher frequency
turbulent structures occur in the wake. In such region, as suggested by the bottom plots in
\Cref{fig:velocity_spectra,fig:pressure_spectra}, the frequency content of the
streamwise velocity and the pressure solutions is definitely richer, and the accuracy of
the reduced PODI and DMD solutions is clearly lower when higher
frequencies are considered. Here, for velocity signals, the plot suggests that considering
$80$ modes both PODI and DMD algorithms lead to good quantitative
spectral reconstructions of the velocity signal for frequencies up to
approximately $110$~\si{Hz}. Interestingly, being based on the train signal yellow
dashed line, in this case the PODI and DMD algorithms seem able to
improve the behavior of the signal interpolated at half the sampling
frequency, and somewhat extend its accuracy at frequencies very close
to the Nyquist one.

\RB{
It is worth commenting on the ROM sensitivity against grid resolution. When a coarser grid is utilized, larger length scales are modeled. As it is expected to note deviations in the corresponding FOM performance with respect to DNS data, it is also doubtful to expect recovering these missing details via a ROM. Furthermore, ROM predictions closely correlate with the energetic content and the induced frequencies from the employed dataset and respective modes. This implies that, while shifting from LES towards RANS grid resolutions, the underlying field dynamics are expected to be sufficiently recovered with less snapshots and less modes for the ROM.
}

\subsection{Acoustic analysis}
The acoustic analysis carried out in this work is aimed at comparing
FWH signal obtained using the full order model data and the
corresponding signal obtained with both the PODI and the DMD reduced order models.
By a practical standpoint, the FWH post processing
described in~\autoref{sec:math2} is applied to the pressure and velocity
fields obtained from the LES full order simulation, and from both the
PODI and DMD reduced models. The FWH integrals are here computed
considering two microphones located at $\mathbf{x}_{\text{mic A}}=\{0,\, 2D,\, 0\}$
(microphone A) and $\mathbf{x}_{\text{mic B}}=\{2D,\, 2D,\, 0\}$ (microphone B). As
mentioned, the sphere is centered at the origin $\mathbf{O}=\{0,\,
0,\, 0\}$ and has diameter $D$. The acoustic pressure time-history is
converted to sound Spectrum Level $\text{SpL} = 20\,\log(p/p_{\text{ref}})$,
considering $p_{\text{ref}}=1$~\si{\mu \si{Pa}}.

\begin{figure*}
\centering
\includegraphics[width=0.9\textwidth]{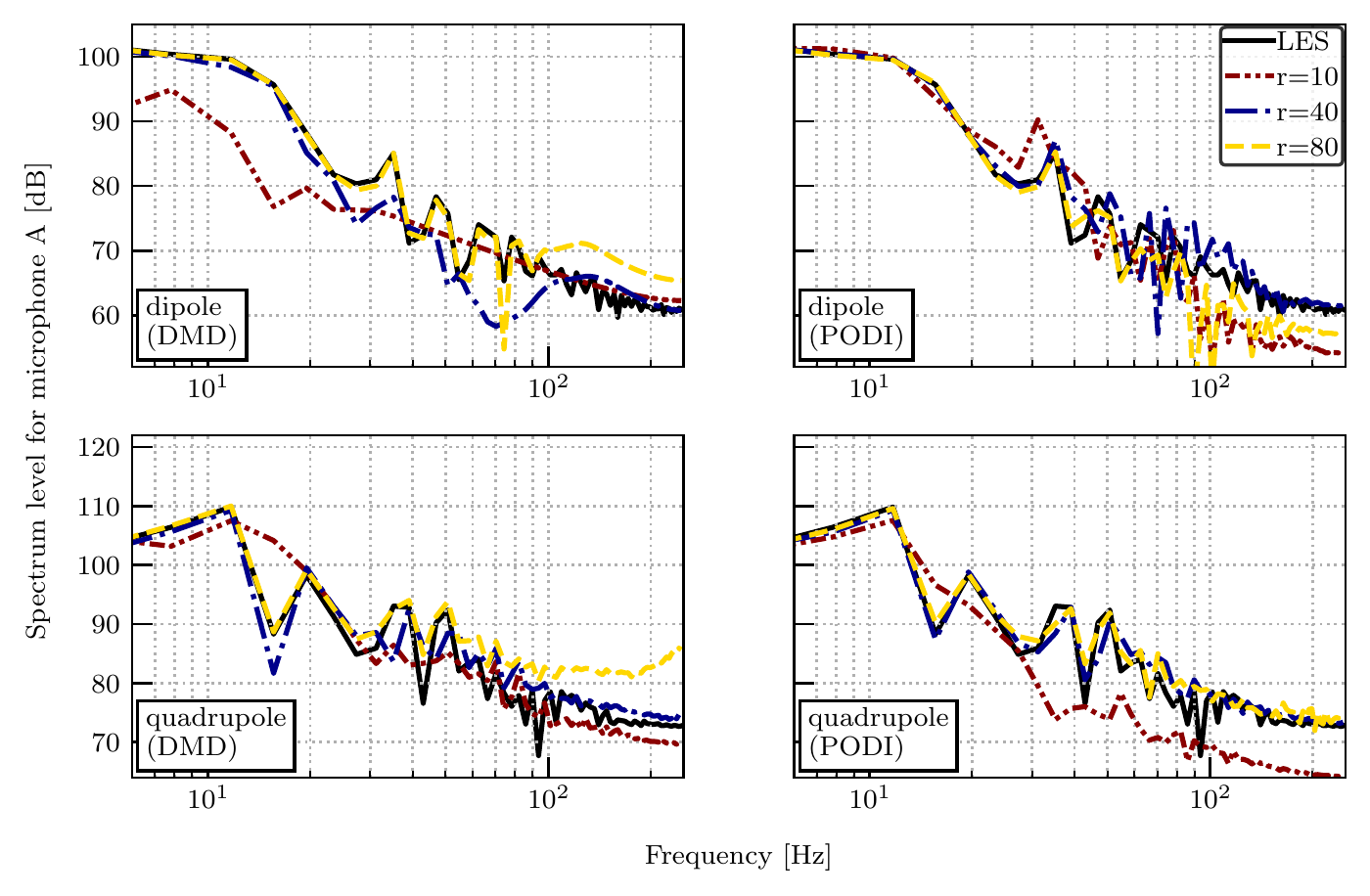}
\caption{Linear dipole (top) and nonlinear quadrupole (bottom) terms
  of FWH equation evaluated from LES data and compared to
  corresponding DMD (left) and PODI (right) reduced models \RB{at different truncation ranks (r)}. Microphone
  A.}
\label{fig:acoA}
\end{figure*}
\begin{figure*}
\centering
\includegraphics[width=0.9\textwidth]{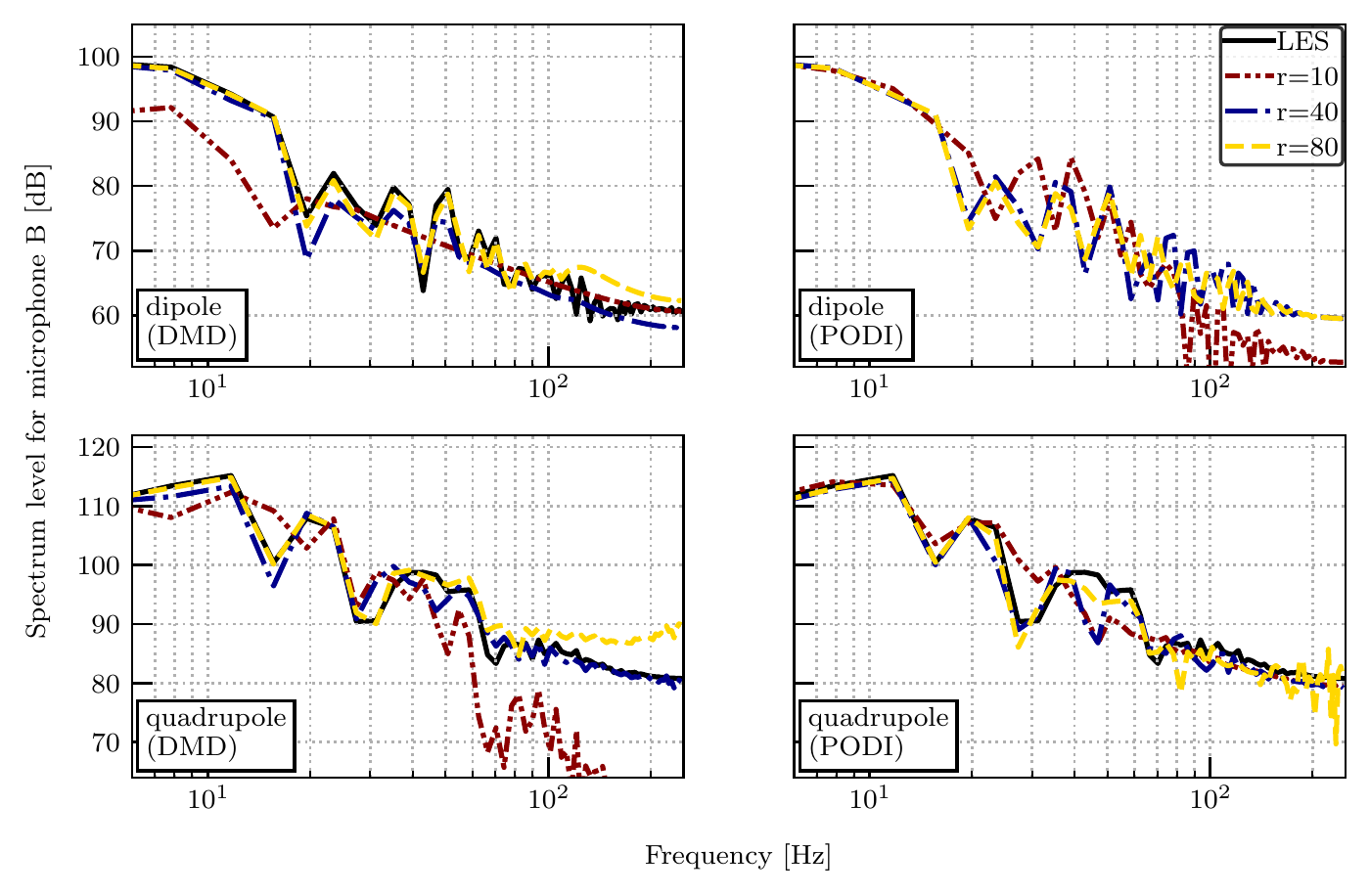}
\caption{Linear dipole (top) and nonlinear quadrupole (bottom) terms
  of FWH equation evaluated from LES data and compared to
  corresponding DMD (left) and PODI (right) reduced models \RB{at different truncation ranks (r)}. Microphone
  B.}
\label{fig:acoB}
\end{figure*}

We will first focus on the results of the PODI model.
The corresponding SpL from both microphones is reported in \Cref{fig:acoA,fig:acoB}, right panels.
In the figures, the different FWH signal curves denote the different
number of modes employed in the PODI reconstruction. Also, as seen
in previous works where the FWH formulation is used, it is convenient
to separate the linear contribution to the acoustic pressure, obtained
from the surface integrals in \autoref{eq:FWH2D}, from the nonlinear
part obtained from the volume integrals in \autoref{eq:FWH3D}.
Such contributions are referred to as dipole and quadrupole terms, respectively.
This distinction is particularly relevant, as the former contribution only
requires the knowledge of the pressure field on the body surface, while
the latter component takes into consideration the evolution of the pressure
and velocity field in the wake region. More specifically, we emphasize that to
obtain an adequate reconstruction of the non-linear acoustic signal, an accurate
reconstruction of the vorticity field is necessary.

In general, the PODI results seem satisfactory
especially for what concerns utilizing $80$ modes. \Cref{fig:acoA,fig:acoB}
(right top plots) show in fact good agreement between the
linear acoustic pressure contribution based on the LES pressure
field and the corresponding linear contribution computed with PODI
employing $80$ modes (yellow lines in the plots).
As the plots suggest, this trend is appreciable
both at microphones A and B. As for the nonlinear contribution,
\Cref{fig:acoA,fig:acoB} (right bottom panels) also show a
satisfactory agreement between the acoustic pressure signals obtained
from the LES flow fields and those computed based on
PODI. In this case, in the plots referring to both
microphones, the blue and yellow lines ($40$ and $80$ modes respectively)
on the right bottom panels are sufficiently close to their reference FOM counterpart,
represented by the continuous black line.

At a closer look, close agreement between PODI and FOM acoustic
pressures, both for the linear and non-linear contributions, can be
established at least up to frequencies in the range $80$--$100$~\si{Hz},
which may be considered the meaningful frequency interval related to
the considered test case. In fact, we point out that, since the reference
pressure is set to $1$~\si{\mu Pa}, the ambient noise may be considered in
the range of $60$--$100$~\si{dB}. At higher frequencies though, we note
that some spurious oscillations appear in the PODI results, especially
in the nonlinear noise terms corresponding to microphone B
(\autoref{fig:acoB}), at which the value of the integrals in the
FWH formulation are strongly conditioned by the presence of the
wake. The observed spurious oscillations could be due to an imperfect
reconstruction of the vorticity field, although it must be emphasized
that the sampling frequency of the original dataset is not very high, $f_s=500$~\si{Hz}, thereby
for the train dataset it becomes $f_s=250$~\si{Hz} for such an alternative arrangement. Therefore, it is
reasonable to obtain less accurate results at frequencies
higher than $125$~\si{Hz}. Further work will be devoted to investigating
whether longer temporal records or higher sampling frequencies will eliminate
or mitigate the problem.

However, considering the $80$ modes PODI reconstruction, the maximum
error is of the order of $5$~\si{dB}, observed at microphone B, for the non
linear terms (\autoref{fig:acoB} bottom right panel).
In addition, the main peaks observed at very
low frequencies, up to $20$~\si{Hz}, are well captured by the PODI
signal. In such case, the agreement is also verified for the $40$
modes signals. As these low frequency peaks describe the main and most
energetic features of the acoustic signal, this indicates that a
moderate amount of modes might still be a good compromise when the
acoustic analysis is only aiming at a general characterization of the
principal noise features. Finally, we must point out that
the errors observed in the non-linear noise contributions are not in general
higher than those introduced in the linear contributions, as might be
expected given the interpolation procedure involved in PODI. In some
occurrence (\autoref{fig:acoA} top right plot), the linear contribution errors
are even quite surprisingly higher than their nonlinear counterparts. This might
be related to the fact that the linear noise contributions are only based on
pressure evaluations on the body surface, while the non-linear contributions
are based on volume integrals. The POD procedure used selects the modal
shapes so as to minimize the error in a norm based on the whole volumetric
solution, rather than only on a surface restriction. For such reason the
modal shape selected might result in non optimal results in the computation
of surface integrals. Possible gains might then be obtained, if needed,
adding a separate PODI only built on the body surface degrees of freedom,
which would result in a modest increase of the computational cost.

As regards the acoustic signals provided by the DMD method, a good agreement is
observed with respect to the FOM spectrum, up to about $80$~\si{Hz}, for both microphones,
and both linear and non-linear terms, see \Cref{fig:acoA,fig:acoB}
(left panels).
The most accurate signals are obtained using $40$ and $80$ modes,
while for the signal related to 10 modes we observe, as expected,
a considerable discrepancy.
In general, the reconstruction of the vorticity field, obtained by both POD and DMD,
is found to provide an adequate input field for the acoustic model.
Particularly, the volume integral of the FWH equation
involves both the velocity and the pressure field in the wake region.
Having observed a good match of the ROM signals with respect to the
FOM reference signals, we may conclude that the entire spectrum
associated with the vortex wake has been adequately reconstructed.

\section{Conclusions and perspectives}
\label{sec:the_end}

This work discussed the application of data-driven dimensionality reduction algorithms on a
hydroacoustic dataset which was numerically measured using Large Eddy Simulation (LES) fluid dynamic turbulence model
and the Ffowcs Williams and Hawkings (FWH) acoustic analogy. An extensive set of
data was presented to fully characterize the ability of both Dynamic Mode Decomposition (DMD) and
Proper Orthogonal Decomposition (POD) to reconstruct the flow fields, based on their spectral and energetic contents,
with all spatial and temporal frequencies needed to support accurate noise predictions.

First, Singular Value Decomposition (SVD) analysis did not indicate the
presence of significant constraints on the modal truncation rank for such flows. In fact,
no significant singular value energy gaps were observed. SVD was
then used to extract DMD and POD modes with associated coefficients, and to employ them for flow fields reconstruction.
In general, both DMD and POD algorithms showed efficient reconstruction accuracy. Spatial and temporal error
analyses indicated relatively lower error magnitudes in POD-based reconstructed fields. On the other hand,
statistical analysis and vortical structures identification methods demonstrated the ability of DMD to
capture finest wake scales by employing higher modes, compared with POD.

Second, two data-driven reduced models based on DMD mid cast and on POD with interpolation (PODI)
were created utilizing half the LES original dataset. Both DMD
and POD based reduced models showed good efficiency and accurate flow reconstruction.
In addition, the spectral analysis of the reduced flow solutions at selected points
inside and outside vortical wake regions indicated that both models were able
to recover most of the model flow frequencies in the ranges of interest for the
acoustic analysis. In particular, PODI showed notable capability to capture
additional frequencies which were present in the original dataset, but not in
the subset employed to train the model. As a consequence, both data-driven reduced
models developed proved efficient and sufficiently accurate in predicting
acoustic noises.

Given the results obtained, introducing bluff body geometrical parameterization as
additional dimension for the PODI analysis is an interesting possibility for future works.
A further future perspective which is currently being explored, is
represented by the development of a reduced model employing the POD modal decomposition
for the Galerkin projection of the continuity and momentum equations for the fluid dynamic
variables. Other studies could be devoted to better understand how the
choice of the interpolator for PODI affects the results.

\appendix

\section*{Appendix}

\section{ROMs cross-comparison}

\begin{figure}
\centering
\includegraphics[width=0.9\textwidth]{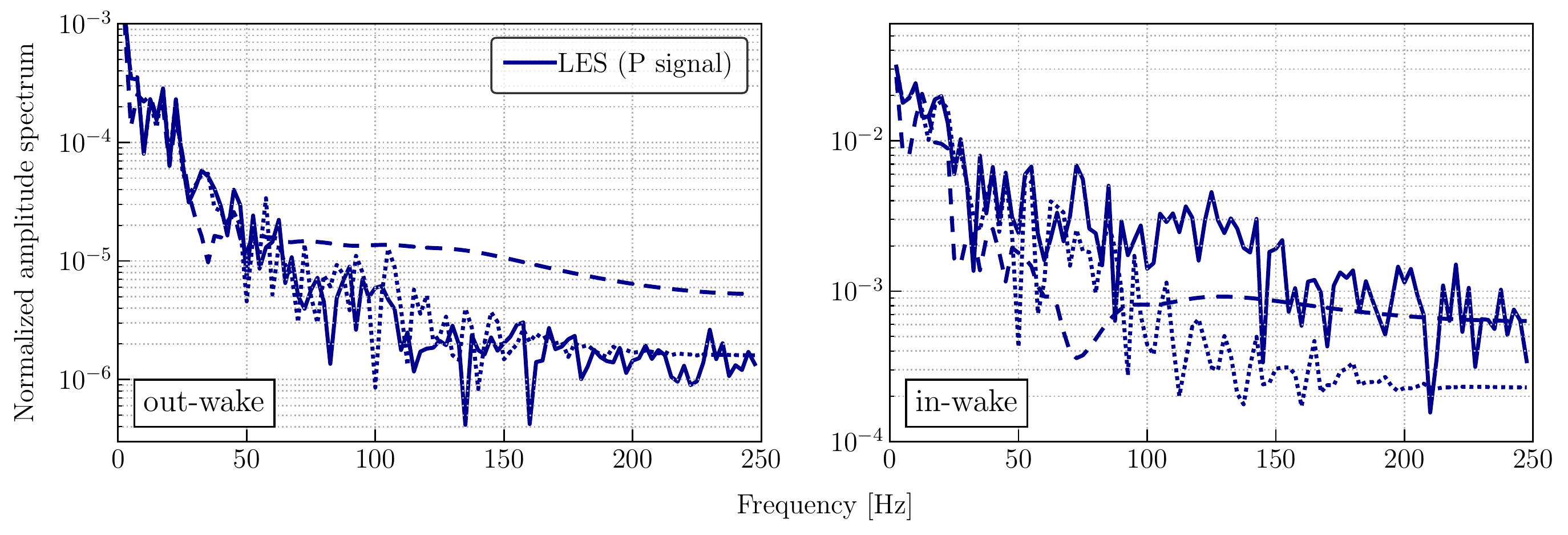}
\includegraphics[width=0.8825\textwidth]{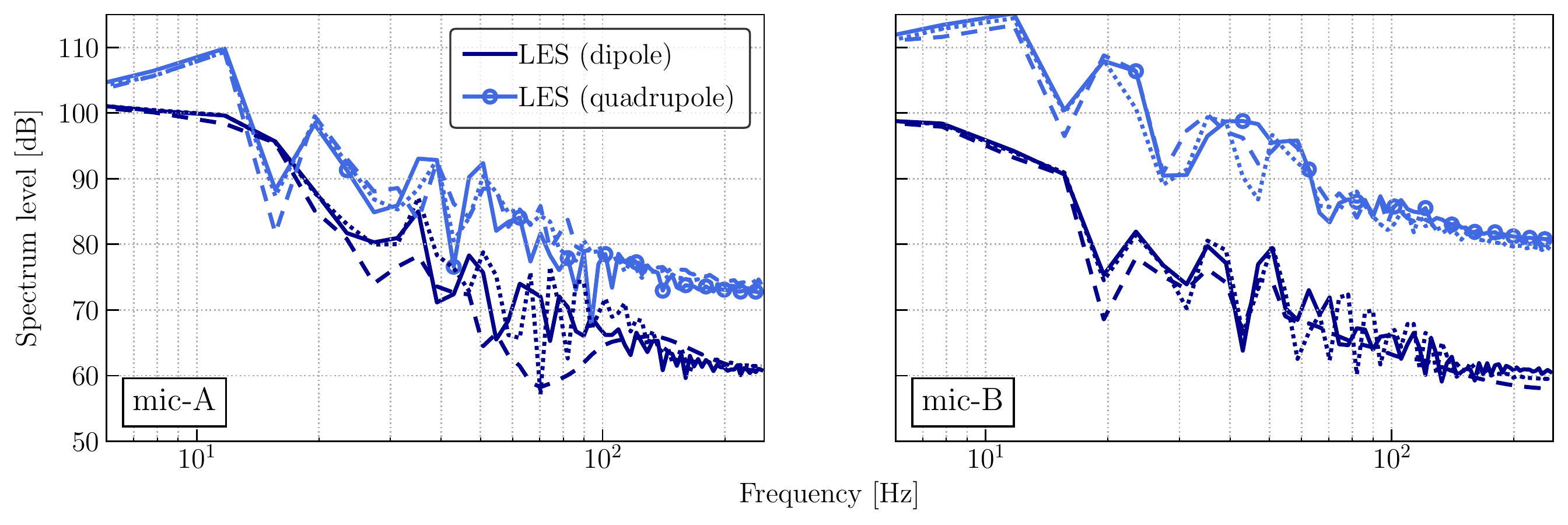}
\caption{Top: FFT spectrum of the pressure signal at points located out of
         wake (left) and inside the wake (right). Bottom: dipole and
         quadrupole acoustic noise spectrum measured at mic-A (left)
         and mic-B (right) locations. PODI
         spectra (\sampleline{line width=1.5 pt, dash pattern=on .15em off .15em on .15em})
         demonstrate better agreement with LES data (\sampleline{line width=1.5 pt}) compared 
         with DMD (\sampleline{line width=1.5 pt, dash pattern=on .5em off .3em on .5em}) ones.
         Both PODI and DMD results in the graphs have been obtained employing
         $r=40$ modes.}
\label{fig:acoMics_r40}
\end{figure}

\RA{For additional clarity, the spectral and acoustic analyses discussed in the
previous sections, cf.\Cref{fig:velocity_spectra,fig:pressure_spectra,fig:acoA,fig:acoB},
are here reported in single graphs directly comparing both DMD and PODI ROMs results against LES data.
In particular, we are considering both DMD and PODI ROMs truncated at $r=40$. The aim is to
provide an easier comparison of the ROMs performance. \autoref{fig:acoMics_r40} depicts FFT spectrum
of the pressure signal (inside and outside the wake, top plots) and the dipole and quadrupole
acoustic noises (microphones A and B, bottom plots). In all plots
DMD (\sampleline{line width=1.5 pt, dash pattern=on .5em off .3em on .5em}) and
PODI (\sampleline{line width=1.5 pt, dash pattern=on .15em off .15em on .15em}) results
are reported alongside LES (\sampleline{line width=1.5 pt}) results. The plots clearly
indicate that, given equal modal truncation level, PODI models appear to accurately recover
a wider range of frequencies compared to DMD models.}

\section*{Acknowledgements}
This work was partially supported by the project PRELICA,
``Advanced methodologies for hydro-acoustic design of naval
propulsion'', supported by Regione FVG, POR-FESR 2014-2020, Piano
Operativo Regionale Fondo Europeo per lo Sviluppo Regionale,
and partially funded by European Union Funding for
Research and Innovation --- Horizon 2020 Program --- in the framework
of European Research Council Executive Agency: H2020 ERC CoG 2015
AROMA-CFD project 681447 ``Advanced Reduced Order Methods with
Applications in Computational Fluid Dynamics'' P.I. Gianluigi Rozza.
We gratefully thank professor Vincenzo Armenio, University of Trieste,
for the general supervision, his precious comments and suggestions,
during the development of the presented work.


\begin{thebibliography}{100}

\bibitem{Baek2015}
D.-G. Baek, H.-S. Yoon, J.-H. Jung, K.-S. Kim, and B.-G. Paik.
\newblock Effects of the advance ratio on the evolution of a propeller wake.
\newblock {\em Computers {\&} Fluids}, 118:32--43, Sept. 2015.

\bibitem{Balaras2015}
E.~Balaras, S.~Schroeder, and A.~Posa.
\newblock {Large-Eddy Simulations of Submarine Propellers}.
\newblock {\em Journal of Ship Research}, 59(4):227--237, Dec. 2015.

\bibitem{Berkooz1993}
G.~Berkooz, P.~Holmes, and J.~L. Lumley.
\newblock {The Proper Orthogonal Decomposition in the Analysis of Turbulent
  Flows}.
\newblock {\em Annual Review of Fluid Mechanics}, 25(1):539--575, Jan. 1993.

\bibitem{Bistrian2015}
D.~A. Bistrian and I.~M. Navon.
\newblock An improved algorithm for the shallow water equations model
  reduction: {D}ynamic {M}ode {D}ecomposition vs {POD}.
\newblock {\em International Journal for Numerical Methods in Fluids},
  78(9):552--580, Apr. 2015.

\bibitem{bistrian2017randomized}
D.~A. Bistrian and I.~M. Navon.
\newblock Randomized dynamic mode decomposition for nonintrusive reduced order
  modelling.
\newblock {\em International Journal for Numerical Methods in Engineering},
  112(1):3--25, 2017.

\bibitem{bistrian2018efficiency}
D.~A. Bistrian and I.~M. Navon.
\newblock Efficiency of randomised dynamic mode decomposition for reduced order
  modelling.
\newblock {\em International Journal of Computational Fluid Dynamics},
  32(2-3):88--103, 2018.

\bibitem{Bouhoubeiny2009}
E.~Bouhoubeiny and P.~Druault.
\newblock Note on the {POD}-based time interpolation from successive {PIV}
  images.
\newblock {\em Comptes Rendus M{\'{e}}canique}, 337(11-12):776--780, Nov. 2009.

\bibitem{Brentner2003}
K.~S. Brentner and F.~Farassat.
\newblock Modeling aerodynamically generated sound of helicopter rotors.
\newblock {\em Progress in Aerospace Sciences}, 39(2-3):83--120, Feb. 2003.

\bibitem{broatch2019dynamic}
A.~Broatch, J.~Garc{\'\i}a-T{\'\i}scar, F.~Roig, and S.~Sharma.
\newblock Dynamic mode decomposition of the acoustic field in radial
  compressors.
\newblock {\em Aerospace Science and Technology}, 90:388--400, 2019.

\bibitem{brunton2019data}
S.~L. Brunton and J.~N. Kutz.
\newblock {\em {Data-Driven Science and Engineering: Machine Learning,
  Dynamical Systems, and Control}}.
\newblock Cambridge University Press, 2019.

\bibitem{bui2003proper}
T.~Bui-Thanh, M.~Damodaran, and K.~Willcox.
\newblock {Proper Orthogonal Decomposition Extensions for Parametric
  Applications in Compressible Aerodynamics}.
\newblock In {\em 21st AIAA Applied Aerodynamics Conference}, page 4213, 2003.

\bibitem{bui2004aerodynamic}
T.~Bui-Thanh, M.~Damodaran, and K.~E. Willcox.
\newblock {Aerodynamic Data Reconstruction and Inverse Design Using Proper
  Orthogonal Decomposition}.
\newblock {\em AIAA journal}, 42(8):1505--1516, 2004.

\bibitem{Campana2006}
E.~F. Campana, D.~Peri, Y.~Tahara, and F.~Stern.
\newblock Shape optimization in ship hydrodynamics using computational fluid
  dynamics.
\newblock {\em Computer Methods in Applied Mechanics and Engineering},
  196(1-3):634--651, Dec. 2006.

\bibitem{CianferraAI18}
M.~Cianferra, V.~Armenio, and S.~Ianniello.
\newblock Hydroacoustic noise from different geometries.
\newblock {\em International Journal Heat Fluid Flow}, 70:348--362, 2018.

\bibitem{CianferraIA19}
M.~Cianferra, S.~Ianniello, and V.~Armenio.
\newblock Assessment of methodologies for the solution of the {Ffowcs}
  {Williams} and {Hawkings} equation using {LES} of incompressible single-phase
  flow around a finite-size square cylinder.
\newblock {\em Journal of Sound and Vibration}, 453:1--24, 2019.

\bibitem{cianferra2018hydrodynamic}
M.~Cianferra, A.~Petronio, and V.~Armenio.
\newblock {Hydrodynamic Noise from a Propeller in Open Sea Condition}.
\newblock In {\em Technology and Science for the Ships of the Future:
  Proceedings of NAV 2018: 19th International Conference on Ship \& Maritime
  Research}, pages 149--156. IOS Press, 2018.

\bibitem{Cianferra2019}
M.~Cianferra, A.~Petronio, and V.~Armenio.
\newblock Non-linear noise from a ship propeller in open sea condition.
\newblock {\em Ocean Engineering}, 191:106474, Nov. 2019.

\bibitem{CintolesiPA15}
C.~Cintolesi, A.~Petronio, and V.~Armenio.
\newblock Large eddy simulation of turbulent buoyant flow in a confined cavity
  with conjugate heat transfer.
\newblock {\em Phys fluids}, 27, 2015.

\bibitem{FosasdePando2014}
M.~F. de~Pando, P.~J. Schmid, and D.~Sipp.
\newblock A global analysis of tonal noise in flows around aerofoils.
\newblock {\em Journal of Fluid Mechanics}, 754:5--38, July 2014.

\bibitem{demo2018shape}
N.~Demo, M.~Tezzele, G.~Gustin, G.~Lavini, and G.~Rozza.
\newblock {Shape Optimization by Means of Proper Orthogonal Decomposition and
  Dynamic Mode Decomposition}.
\newblock In {\em Technology and Science for the Ships of the Future:
  Proceedings of NAV 2018: 19th International Conference on Ship \& Maritime
  Research}, pages 212--219. IOS Press, 2018.

\bibitem{demo2018isope}
N.~Demo, M.~Tezzele, A.~Mola, and G.~Rozza.
\newblock {An Efficient Shape Parametrisation by Free-Form Deformation Enhanced
  by Active Subspace for Hull Hydrodynamic Ship Design Problems in Open Source
  Environment}.
\newblock In {\em Proceedings of ISOPE 2018: The 28th International Ocean and
  Polar Engineering Conference}, volume~3, pages 565--572, 2018.

\bibitem{demo2019marine}
N.~Demo, M.~Tezzele, A.~Mola, and G.~Rozza.
\newblock A complete data-driven framework for the efficient solution of
  parametric shape design and optimisation in naval engineering problems.
\newblock In R.~Bensow and J.~Ringsberg, editors, {\em Proceedings of MARINE
  2019: VIII International Conference on Computational Methods in Marine
  Engineering}, pages 111--121, 2019.

\bibitem{demo2018ezyrb}
N.~Demo, M.~Tezzele, and G.~Rozza.
\newblock {EZyRB: Easy Reduced Basis method}.
\newblock {\em The Journal of Open Source Software}, 3(24):661, 2018.

\bibitem{demo18pydmd}
N.~Demo, M.~Tezzele, and G.~Rozza.
\newblock {PyDMD: Python Dynamic Mode Decomposition}.
\newblock {\em The Journal of Open Source Software}, 3(22):530, 2018.

\bibitem{demo2019cras}
N.~Demo, M.~Tezzele, and G.~Rozza.
\newblock {A non-intrusive approach for reconstruction of POD modal
  coefficients through active subspaces}.
\newblock {\em Comptes Rendus M\'ecanique de l'Acad\'emie des Sciences,
  DataBEST 2019 Special Issue}, 347(11):873--881, November 2019.

\bibitem{dolci2016proper}
V.~Dolci and R.~Arina.
\newblock {Proper Orthogonal Decomposition as Surrogate Model for Aerodynamic
  Optimization}.
\newblock {\em International Journal of Aerospace Engineering}, 2016.

\bibitem{Druault2007}
P.~Druault and C.~Chaillou.
\newblock {Use of Proper Orthogonal Decomposition for reconstructing the 3D
  in-cylinder mean-flow field from PIV data}.
\newblock {\em Comptes Rendus M{\'{e}}canique}, 335(1):42--47, Jan. 2007.

\bibitem{erichson2016compressed}
N.~B. Erichson, S.~L. Brunton, and J.~N. Kutz.
\newblock Compressed dynamic mode decomposition for background modeling.
\newblock {\em Journal of Real-Time Image Processing}, 16(5):1479--1492, Oct
  2019.

\bibitem{fwh1969}
J.~E. Ffowcs~Williams and D.~L. Hawkings.
\newblock Sound generation by turbulence and surfaces in arbitrary motion.
\newblock {\em Philosophical Transaction of Royal Society}, 264:321--342, 1969.

\bibitem{Fossati2013}
M.~Fossati and W.~G. Habashi.
\newblock {Multiparameter Analysis of Aero-Icing Problems Using Proper
  Orthogonal Decomposition and Multidimensional Interpolation}.
\newblock {\em {AIAA} Journal}, 51(4):946--960, Apr. 2013.

\bibitem{Franz2014}
T.~Franz, R.~Zimmermann, S.~G\"{o}rtz, and N.~Karcher.
\newblock Interpolation-based reduced-order modelling for steady transonic
  flows via manifold learning.
\newblock {\em International Journal of Computational Fluid Dynamics},
  28(3-4):106--121, Mar. 2014.

\bibitem{gadalla19bladex}
M.~Gadalla, M.~Tezzele, A.~Mola, and G.~Rozza.
\newblock {BladeX: Python Blade Morphing}.
\newblock {\em The Journal of Open Source Software}, 4(34):1203, 2019.

\bibitem{garotta2018quiet}
F.~Garotta, N.~Demo, M.~Tezzele, M.~Carraturo, A.~Reali, and G.~Rozza.
\newblock {Reduced Order Isogeometric Analysis Approach for PDEs in
  Parametrized Domains}.
\newblock In M.~D'Elia, M.~Gunzburger, and G.~Rozza, editors, {\em
  Quantification of Uncertainty: Improving Efficiency and Technology: QUIET
  selected contributions}, volume 137 of {\em Lecture Notes in Computational
  Science and Engineering}, pages 153--170. Springer International Publishing,
  Cham, 2020.

\bibitem{GeoStaRoBlu2019}
S.~Georgaka, G.~Stabile, K.~Star, G.~Rozza, and M.~J. Bluck.
\newblock {A hybrid reduced order method for modelling turbulent heat transfer
  problems}.
\newblock {\em Computers {\&} Fluids}, 208:104615, 2020.

\bibitem{GLEGG2001}
S.~A. Glegg and W.~J. Devenport.
\newblock Proper orthogonal decomposition of turbulent flows for aeroacoustic
  and hydroacoustic applications.
\newblock {\em Journal of Sound and Vibration}, 239(4):767--784, 2001.

\bibitem{Gloerfelt2005}
X.~Gloerfelt, F.~P{\'{e}}rot, C.~Bailly, and D.~Juv{\'{e}}.
\newblock Flow-induced cylinder noise formulated as a diffraction problem for
  low {Mach} numbers.
\newblock {\em Journal of Sound and Vibration}, 287(1-2):129--151, Oct. 2005.

\bibitem{LeGuennec2018}
Y.~L. Guennec, J.-P. Brunet, F.-Z. Daim, M.~Chau, and Y.~Tourbier.
\newblock A parametric and non-intrusive reduced order model of car crash
  simulation.
\newblock {\em Computer Methods in Applied Mechanics and Engineering},
  338:186--207, Aug. 2018.

\bibitem{HALLER2005}
G.~Haller.
\newblock An objective definition of a vortex.
\newblock {\em Journal of Fluid Mechanics}, 525:1--26, Feb. 2005.

\bibitem{hesthaven2016certified}
J.~S. Hesthaven, G.~Rozza, B.~Stamm, et~al.
\newblock {\em {Certified Reduced Basis Methods for Parametrized Partial
  Differential Equations}}.
\newblock Springer, 2016.

\bibitem{HiStaMoRo2019}
S.~Hijazi, G.~Stabile, A.~Mola, and G.~Rozza.
\newblock {Data-Driven POD--Galerkin reduced order model for turbulent flows}.
\newblock {\em Journal of Computational Physics}, 416:109513, 2020.

\bibitem{Ianniello2013}
S.~Ianniello, R.~Muscari, and A.~D. Mascio.
\newblock {Ship underwater noise assessment by the acoustic analogy. Part I:
  nonlinear analysis of a marine propeller in a uniform flow}.
\newblock {\em Journal of Marine Science and Technology}, 18(4):547--570, July
  2013.

\bibitem{Ianniello2014}
S.~Ianniello, R.~Muscari, and A.~D. Mascio.
\newblock {Ship underwater noise assessment by the acoustic analogy. Part II:
  hydroacoustic analysis of a ship scaled model}.
\newblock {\em J Mar Sci Technol}, 18(4):547--570, 2013.

\bibitem{Iuliano2013}
E.~Iuliano and D.~Quagliarella.
\newblock {Proper Orthogonal Decomposition, surrogate modelling and
  evolutionary optimization in aerodynamic design}.
\newblock {\em Computers {\&} Fluids}, 84:327--350, Sept. 2013.

\bibitem{Jourdain2013}
G.~Jourdain, L.-E. Eriksson, S.~H. Kim, and C.~H. Sohn.
\newblock Application of dynamic mode decomposition to acoustic-modes
  identification and damping in a 3-dimensional chamber with baffled injectors.
\newblock {\em Journal of Sound and Vibration}, 332(18):4308--4323, Sept. 2013.

\bibitem{Jovanovi2014}
M.~R. Jovanovi{\'c}, P.~J. Schmid, and J.~W. Nichols.
\newblock Sparsity-promoting dynamic mode decomposition.
\newblock {\em Physics of Fluids}, 26(2):024103, Feb 2014.

\bibitem{Karri2009}
S.~Karri, J.~Charonko, and P.~P. Vlachos.
\newblock Robust wall gradient estimation using radial basis functions and
  proper orthogonal decomposition ({POD}) for particle image velocimetry
  ({PIV}) measured fields.
\newblock {\em Measurement Science and Technology}, 20(4):045401, Feb. 2009.

\bibitem{Kerwin1986}
J.~E. Kerwin.
\newblock {Marine Propellers}.
\newblock {\em Annual Review of Fluid Mechanics}, 18(1):367--403, Jan. 1986.

\bibitem{Kumar2017}
P.~Kumar and K.~Mahesh.
\newblock Large eddy simulation of propeller wake~instabilities.
\newblock {\em Journal of Fluid Mechanics}, 814:361--396, Feb. 2017.

\bibitem{Kunisch2008}
K.~Kunisch and S.~Volkwein.
\newblock Proper orthogonal decomposition for optimality systems.
\newblock {\em {ESAIM}: Mathematical Modelling and Numerical Analysis},
  42(1):1--23, Jan. 2008.

\bibitem{kutz2016dynamic}
J.~N. Kutz, S.~L. Brunton, B.~W. Brunton, and J.~L. Proctor.
\newblock {\em {Dynamic Mode Decomposition: Data-Driven Modeling of Complex
  Systems}}, volume 149.
\newblock SIAM, 2016.

\bibitem{kutz2016multiresolution}
J.~N. Kutz, X.~Fu, and S.~L. Brunton.
\newblock {Multiresolution Dynamic Mode Decomposition}.
\newblock {\em SIAM Journal on Applied Dynamical Systems}, 15(2):713--735,
  2016.

\bibitem{le2017higher}
S.~Le~Clainche and J.~M. Vega.
\newblock {Higher Order Dynamic Mode Decomposition}.
\newblock {\em SIAM Journal on Applied Dynamical Systems}, 16(2):882--925,
  2017.

\bibitem{Lorente2008}
L.~S. Lorente, J.~M. Vega, and A.~Velazquez.
\newblock {Generation of Aerodynamics Databases Using High-Order Singular Value
  Decomposition}.
\newblock {\em Journal of Aircraft}, 45(5):1779--1788, Sept. 2008.

\bibitem{Ly2001}
H.~V. Ly and H.~T. Tran.
\newblock Modeling and control of physical processes using proper orthogonal
  decomposition.
\newblock {\em Mathematical and Computer Modelling}, 33(1-3):223--236, Jan.
  2001.

\bibitem{Mancinelli2017}
M.~Mancinelli, T.~Pagliaroli, R.~Camussi, and T.~Castelain.
\newblock On the hydrodynamic and acoustic nature of pressure proper orthogonal
  decomposition modes in the near field of a compressible jet.
\newblock {\em Journal of Fluid Mechanics}, 836:998--1008, Dec. 2017.

\bibitem{markovsky2019low}
I.~Markovsky.
\newblock {\em {Low-Rank Approximation}}.
\newblock Communications and Control Engineering. Springer International
  Publishing, 2 edition, 2019.

\bibitem{DiMascio2014}
A.~D. Mascio, R.~Muscari, and G.~Dubbioso.
\newblock On the wake dynamics of a propeller operating in drift.
\newblock {\em Journal of Fluid Mechanics}, 754:263--307, July 2014.

\bibitem{MeneveauLC96}
C.~Meneveau, T.~Lund, and W.~Cabot.
\newblock A {L}agrangian dynamic subgrid-scale model of turbulence.
\newblock {\em J Fluid Mech}, 319:353--385, 1996.

\bibitem{MEYER2007}
K.~E. Meyer, J.~M. Pedersen, and O.~\"{O}zcan.
\newblock A turbulent jet in crossflow analysed with proper orthogonal
  decomposition.
\newblock {\em Journal of Fluid Mechanics}, 583:199--227, July 2007.

\bibitem{Mifsud2016}
M.~J. Mifsud, D.~G. MacManus, and S.~Shaw.
\newblock A variable-fidelity aerodynamic model using proper orthogonal
  decomposition.
\newblock {\em International Journal for Numerical Methods in Fluids},
  82(10):646--663, Apr. 2016.

\bibitem{Mifsud2009}
M.~J. Mifsud, S.~T. Shaw, and D.~G. MacManus.
\newblock A high-fidelity low-cost aerodynamic model using proper orthogonal
  decomposition.
\newblock {\em International Journal for Numerical Methods in Fluids}, pages
  n/a--n/a, 2009.

\bibitem{mola2019marine}
A.~Mola, M.~Tezzele, M.~Gadalla, F.~Valdenazzi, D.~Grassi, R.~Padovan, and
  G.~Rozza.
\newblock {Efficient Reduction in Shape Parameter Space Dimension for Ship
  Propeller Blade Design}.
\newblock In R.~Bensow and J.~Ringsberg, editors, {\em Proceedings of MARINE
  2019: VIII International Conference on Computational Methods in Marine
  Engineering}, pages 201--212, 2019.

\bibitem{muld2012flow}
T.~W. Muld, G.~Efraimsson, and D.~S. Henningson.
\newblock {Flow structures around a high-speed train extracted using Proper
  Orthogonal Decomposition and Dynamic Mode Decomposition}.
\newblock {\em Computers \& Fluids}, 57:87--97, 2012.

\bibitem{NajafiBM11}
A.~Najafi-Yazdi, G.A.Bres, and L.~Mongeau.
\newblock An acoustic analogy formulation for moving sources in uniformly
  moving media.
\newblock {\em Proc R Soc Lond}, A467:144--165, 2011.

\bibitem{Nitzkorski2014}
Z.~Nitzkorski and K.~Mahesh.
\newblock A dynamic end cap technique for sound computation using the {Ffowcs}
  {Williams} and {Hawkings} equations.
\newblock {\em Physics of Fluids}, 26(11):115101, Nov. 2014.

\bibitem{PiomelliCS97}
U.~Piomelli, L.~S. Craig, and S.~Sarkar.
\newblock On the computation of sound by large-eddy simulations.
\newblock {\em Journal of Engineering and Mathematics}, 32:217--236, 1997.

\bibitem{PosaBalaras}
A.~Posa, R.~Broglia, M.~Felli, M.~Falchi, and E.~Balaras.
\newblock {Characterization of the wake of a submarine propeller via Large-Eddy
  simulation}.
\newblock {\em Computer and Fluids}, 184:138--152, 2019.

\bibitem{proctor2016dynamic}
J.~L. Proctor, S.~L. Brunton, and J.~N. Kutz.
\newblock {Dynamic Mode Decomposition with Control}.
\newblock {\em SIAM Journal on Applied Dynamical Systems}, 15(1):142--161,
  2016.

\bibitem{quarteroni2015reduced}
A.~Quarteroni, A.~Manzoni, and F.~Negri.
\newblock {\em {Reduced Basis Methods for Partial Differential Equations: An
  Introduction}}, volume~92.
\newblock Springer, 2015.

\bibitem{ripepi2018reduced}
M.~Ripepi, M.~Verveld, N.~Karcher, T.~Franz, M.~Abu-Zurayk, S.~G{\"o}rtz, and
  T.~Kier.
\newblock {Reduced-order models for aerodynamic applications, loads and MDO}.
\newblock {\em CEAS Aeronautical Journal}, 9(1):171--193, 2018.

\bibitem{Rowley2004}
C.~W. Rowley, T.~Colonius, and R.~M. Murray.
\newblock Model reduction for compressible flows using {POD} and {G}alerkin
  projection.
\newblock {\em Physica D: Nonlinear Phenomena}, 189(1-2):115--129, Feb. 2004.

\bibitem{rowley2009spectral}
C.~W. Rowley, I.~Mezi{\'c}, S.~Bagheri, P.~Schlatter, and D.~S. Henningson.
\newblock Spectral analysis of nonlinear flows.
\newblock {\em Journal of fluid mechanics}, 641:115--127, 2009.

\bibitem{morhandbook2019}
G.~Rozza, M.~W. Hess, G.~Stabile, M.~Tezzele, and F.~Ballarin.
\newblock {Basic Ideas and Tools for Projection-Based Model Reduction of
  Parametric Partial Differential Equations }.
\newblock In P.~Benner, S.~Grivet-Talocia, A.~Quarteroni, G.~Rozza, W.~H.~A.
  Schilders, and L.~M. Silveira, editors, {\em Handbook on Model Order
  Reduction}, volume~2, chapter~1. De Gruyter, In Press, 2020.

\bibitem{rozza2018advances}
G.~Rozza, M.~H. Malik, N.~Demo, M.~Tezzele, M.~Girfoglio, G.~Stabile, and
  A.~Mola.
\newblock {Advances in Reduced Order Methods for Parametric Industrial Problems
  in Computational Fluid Dynamics}.
\newblock In R.~Owen, R.~de~Borst, J.~Reese, and P.~Chris, editors, {\em
  ECCOMAS ECFD 7 - Proceedings of 6th European Conference on Computational
  Mechanics (ECCM 6) and 7th European Conference on Computational Fluid
  Dynamics (ECFD 7)}, pages 59--76, Glasgow, UK, 2018.

\bibitem{rudy2018deep}
S.~H. Rudy, J.~N. Kutz, and S.~L. Brunton.
\newblock Deep learning of dynamics and signal-noise decomposition with
  time-stepping constraints.
\newblock {\em J. Comput. Physics}, 396:483--506, 2018.

\bibitem{SakamotoVortShedding}
H.~Sakamoto and H.~Haniu.
\newblock {A Study on Vortex Shedding From Spheres in a Uniform Flow}.
\newblock {\em Journal of Fluids Engineering}, 112(4):386--392, 12 1990.

\bibitem{salmoiraghi2016advances}
F.~Salmoiraghi, F.~Ballarin, G.~Corsi, A.~Mola, M.~Tezzele, and G.~Rozza.
\newblock Advances in geometrical parametrization and reduced order models and
  methods for computational fluid dynamics problems in applied sciences and
  engineering: overview and perspectives.
\newblock {\em ECCOMAS Congress 2016 - Proceedings of the 7th European Congress
  on Computational Methods in Applied Sciences and Engineering}, 1:1013--1031,
  2016.

\bibitem{salmoiraghi2018free}
F.~Salmoiraghi, A.~Scardigli, H.~Telib, and G.~Rozza.
\newblock Free-form deformation, mesh morphing and reduced-order methods:
  enablers for efficient aerodynamic shape optimisation.
\newblock {\em International Journal of Computational Fluid Dynamics},
  32(4-5):233--247, 2018.

\bibitem{SANDBERG2008}
R.~D. Sandberg and N.~D. Sandham.
\newblock Direct numerical simulation of turbulent flow past a trailing edge
  and the associated noise generation.
\newblock {\em Journal of Fluid Mechanics}, 596:353--385, Jan. 2008.

\bibitem{schmid2010dynamic}
P.~J. Schmid.
\newblock Dynamic mode decomposition of numerical and experimental data.
\newblock {\em Journal of fluid mechanics}, 656:5--28, 2010.

\bibitem{schmid2011application}
P.~J. Schmid.
\newblock Application of the dynamic mode decomposition to experimental data.
\newblock {\em Experiments in fluids}, 50(4):1123--1130, 2011.

\bibitem{schmid2011applications}
P.~J. Schmid, L.~Li, M.~P. Juniper, and O.~Pust.
\newblock Applications of the dynamic mode decomposition.
\newblock {\em Theoretical and Computational Fluid Dynamics}, 25(1-4):249--259,
  2011.

\bibitem{Seena2011}
A.~Seena and H.~J. Sung.
\newblock Dynamic mode decomposition of turbulent cavity flows for
  self-sustained oscillations.
\newblock {\em International Journal of Heat and Fluid Flow}, 32(6):1098--1110,
  Dec. 2011.

\bibitem{Seo2008}
J.~H. Seo, Y.~J. Moon, and B.~R. Shin.
\newblock Prediction of cavitating flow noise by direct numerical simulation.
\newblock {\em Journal of Computational Physics}, 227(13):6511--6531, June
  2008.

\bibitem{SEOL2002}
H.~Seol, B.~Jung, J.-C. Suh, and S.~Lee.
\newblock Prediction of non-cavitating underwater propeller noise.
\newblock {\em Journal of Sound and Vibration}, 257(1):131--156, Oct. 2002.

\bibitem{Seror2001}
C.~Seror, P.~Sagaut, C.~Bailly, and D.~Juv{\'{e}}.
\newblock On the radiated noise computed by large-eddy simulation.
\newblock {\em Physics of Fluids}, 13(2):476--487, Feb. 2001.

\bibitem{Shen2020}
W.~Shen, T.~K. Patel, and S.~A. Miller.
\newblock {Extraction of Large-Scale Coherent Structures from Large Eddy
  Simulation of Supersonic Jets for Shock-Associated Noise Prediction}.
\newblock In {\em {AIAA} Scitech 2020 Forum}. American Institute of Aeronautics
  and Astronautics, Jan. 2020.

\bibitem{Sirovich1987}
L.~Sirovich.
\newblock Turbulence and the dynamics of coherent structures. {II}. symmetries
  and transformations.
\newblock {\em Quarterly of Applied Mathematics}, 45(3):573--582, Oct. 1987.

\bibitem{StaBaZuRo2019}
G.~Stabile, F.~Ballarin, G.~Zuccarino, and G.~Rozza.
\newblock {A reduced order variational multiscale approach for turbulent
  flows}.
\newblock {\em Advances in Computational Mathematics}, 45:2349--2368, 2019.

\bibitem{StaHiMoLoRo2017}
G.~Stabile, S.~Hijazi, A.~Mola, S.~Lorenzi, and G.~Rozza.
\newblock {POD-Galerkin reduced order methods for CFD using Finite Volume
  Discretisation: vortex shedding around a circular cylinder}.
\newblock {\em Communications in Applied and Industrial Mathematics},
  8(1):210--236, 2017.

\bibitem{stabile_stabilized}
G.~Stabile and G.~Rozza.
\newblock {Finite volume POD-Galerkin stabilised reduced order methods for the
  parametrised incompressible Navier--Stokes equations}.
\newblock {\em Computers {\&} Fluids}, 173:273--284, sep 2018.

\bibitem{stegeman2015proper}
P.~Stegeman, A.~Ooi, and J.~Soria.
\newblock {Proper Orthogonal Decomposition and Dynamic Mode Decomposition of
  Under-Expanded Free-Jets with Varying Nozzle Pressure Ratios}.
\newblock In {\em Instability and Control of Massively Separated Flows}, pages
  85--90. Springer, 2015.

\bibitem{tezzele2018model}
M.~Tezzele, N.~Demo, M.~Gadalla, A.~Mola, and G.~Rozza.
\newblock {Model Order Reduction by Means of Active Subspaces and Dynamic Mode
  Decomposition for Parametric Hull Shape Design Hydrodynamics}.
\newblock In {\em Technology and Science for the Ships of the Future:
  Proceedings of NAV 2018: 19th International Conference on Ship \& Maritime
  Research}, pages 569--576. IOS Press, 2018.

\bibitem{tezzele2018ecmi}
M.~Tezzele, N.~Demo, A.~Mola, and G.~Rozza.
\newblock {An integrated data-driven computational pipeline with model order
  reduction for industrial and applied mathematics}.
\newblock {\em Special Volume ECMI, Springer, In Press}, 2020.

\bibitem{pygem}
M.~Tezzele, N.~Demo, A.~Mola, and G.~Rozza.
\newblock {PyGeM}: Python geometrical morphing.
\newblock {\em Software Impacts}, page 100047, 2020.

\bibitem{tezzele2019marine}
M.~Tezzele, N.~Demo, and G.~Rozza.
\newblock Shape optimization through proper orthogonal decomposition with
  interpolation and dynamic mode decomposition enhanced by active subspaces.
\newblock In R.~Bensow and J.~Ringsberg, editors, {\em Proceedings of MARINE
  2019: VIII International Conference on Computational Methods in Marine
  Engineering}, pages 122--133, 2019.

\bibitem{tezzele2019mortech}
M.~Tezzele, N.~Demo, G.~Stabile, A.~Mola, and G.~Rozza.
\newblock {Enhancing CFD predictions in shape design problems by model and
  parameter space reduction}.
\newblock {\em Advanced Modeling and Simulation in Engineering Sciences},
  7(40), 2020.

\bibitem{Tezzele2018Dimension}
M.~Tezzele, F.~Salmoiraghi, A.~Mola, and G.~Rozza.
\newblock Dimension reduction in heterogeneous parametric spaces with
  application to naval engineering shape design problems.
\newblock {\em Advanced Modeling and Simulation in Engineering Sciences}, 5(1),
  2018.

\bibitem{valdenazzi2019marine}
F.~Valdenazzi, F.~Conti, S.~Gaggero, D.~Grassi, C.~Vaccaro, and D.~Villa.
\newblock A practical tool for the hydro-acoustic optimization of naval
  propellers.
\newblock In R.~Bensow and J.~Ringsberg, editors, {\em Proceedings of MARINE
  2019: VIII International Conference on Computational Methods in Marine
  Engineering}, pages 296--308, 2019.

\bibitem{of}
H.~G. Weller, G.~Tabor, H.~Jasak, and C.~Fureby.
\newblock A tensorial approach to computational continuum mechanics using
  object-oriented techniques.
\newblock {\em Computers in physics}, 12(6):620--631, 1998.

\bibitem{williams2015data}
M.~O. Williams, I.~G. Kevrekidis, and C.~W. Rowley.
\newblock {A Data--Driven Approximation of the Koopman Operator: Extending
  Dynamic Mode Decomposition}.
\newblock {\em Journal of Nonlinear Science}, 25(6):1307--1346, 2015.

\bibitem{Xiao2015}
D.~Xiao, F.~Fang, C.~Pain, and G.~Hu.
\newblock Non-intrusive reduced-order modelling of the {N}avier-{S}tokes
  equations based on {RBF} interpolation.
\newblock {\em International Journal for Numerical Methods in Fluids},
  79(11):580--595, July 2015.

\bibitem{zhang2014identification}
Q.~Zhang, Y.~Liu, and S.~Wang.
\newblock The identification of coherent structures using proper orthogonal
  decomposition and dynamic mode decomposition.
\newblock {\em Journal of Fluids and Structures}, 49:53--72, 2014.

\end{thebibliography}

\end{document}